\documentclass[twoside]{article}

\usepackage[accepted]{aistats2022}
\usepackage[round]{natbib}

\usepackage[utf8]{inputenc} 
\usepackage{hyperref}       
\usepackage{url}            
\usepackage{booktabs}       
\usepackage{amsfonts}       
\usepackage{nicefrac}       
\usepackage{microtype}      
\usepackage{diagbox}
\usepackage{layout}
\usepackage{multirow}

\usepackage{amsfonts}
\usepackage{amsmath}
\usepackage{amsthm}
\usepackage{amsbsy}
\usepackage{amssymb} 

\usepackage{nicefrac}
\usepackage{mathtools}

\usepackage{mdframed} 
\usepackage{thmtools}

\definecolor{shadecolor}{gray}{0.90}
\declaretheoremstyle[
headfont=\normalfont\bfseries,
notefont=\mdseries, notebraces={(}{)},
bodyfont=\normalfont,
postheadspace=0.5em,
spaceabove=0.5em,
spacebelow=0.5em,
mdframed={
  skipabove=8pt,
  skipbelow=8pt,
  hidealllines=true,
  backgroundcolor={shadecolor},
  innerleftmargin=4pt,
  innerrightmargin=4pt}
]{shaded}

\declaretheorem[style=shaded,within=section]{definition}
\declaretheorem[style=shaded,sibling=definition]{theorem}
\declaretheorem[style=shaded,sibling=definition]{proposition}
\declaretheorem[style=shaded,sibling=definition]{assumption}

\declaretheorem[style=shaded,sibling=definition]{lemma}
\declaretheorem[style=shaded,sibling=definition]{remark}
\declaretheorem[style=shaded,sibling=definition]{example}
\declaretheorem[style=shaded]{objective}

\usepackage[table]{xcolor}
\usepackage{color}
\usepackage{graphicx}
\graphicspath{{./figures/}}  

\usepackage{hyperref}
\usepackage{textcomp}
\usepackage{cases}  

\newcommand{\R}{\mathbb{R}} 
\newcommand{\N}{\mathbb{N}} 

\newcommand{\cC}{{\cal C}}
\newcommand{\cD}{{\cal D}}

\newcommand{\cO}{{\cal O}}

\newcommand{\mo}{{\bf 0}}
\newcommand{\mA}{{\bf A}}

\newcommand{\mD}{{\bf D}}

\newcommand{\mE}{{\bf E}}

\newcommand{\mG}{{\bf G}}
\newcommand{\mH}{{\bf H}}
\newcommand{\mI}{{\bf I}}

\newcommand{\mM}{{\bf M}}

\newcommand{\mS}{{\bf S}}

\newcommand{\mW}{{\bf W}}

\usepackage[colorinlistoftodos,bordercolor=orange,backgroundcolor=orange!20,linecolor=orange]{todonotes}

\newcommand{\eqdef}{\overset{\text{def}}{=}} 

\newcommand{\dotprod}[1]{\left< #1\right>} 
\newcommand{\norm}[1]{ \left\| #1 \right\|}      

\newcommand{\Prob}[1]{\mathbb{P}[#1]}
\providecommand{\Null}[1]{{\bf Null}\left( #1\right)}
\providecommand{\Image}[1]{{\bf Im}\left( #1\right)}

\DeclareMathOperator{\argmin}{argmin}        

\newcommand{\Diag}[1]{\mathbf{Diag}\left( #1\right)}

\providecommand{\trace}[1]{{\bf Trace}\left( #1\right)}

\newcommand{\E}[1]{\mathbb{E}\left[#1\right] } 
\newcommand{\EE}[2]{\mathbb{E}_{#1}\left[#2\right] }

\usepackage{appendix}
\usepackage{wrapfig}
\usepackage{xspace}
\usepackage{float}
\usepackage{algorithmicx}
\usepackage[ruled, lined, linesnumbered, commentsnumbered, longend, vlined]{algorithm2e}
\usepackage{subfigure}

\usepackage{enumitem}
\setlist[itemize]{leftmargin=*} 

\definecolor{pearThree}{HTML}{E74C3C}
\definecolor{pearcomp}{HTML}{B97E29}
\definecolor{pearDark}{HTML}{2980B9}
\definecolor{pearDarker}{HTML}{1D2DEC}
\hypersetup{
  colorlinks,
  citecolor=pearDark,
  linkcolor=pearThree,
  breaklinks=true,
  urlcolor=pearDarker}

\SetKwBlock{With}{With probability $\pi$ update:}{}
\SetKwBlock{Otherwise}{Otherwise with probability $(1-\pi)$:}{}

\usepackage{minitoc}
\usepackage[hang,flushmargin]{footmisc}
\newcommand\affiliations[1]{%
  \begingroup
  \renewcommand\thefootnote{}\footnote{#1}%
  \addtocounter{footnote}{-1}%
  \endgroup
}

\begin{document}

%
\runningtitle{SAN: Stochastic Average Newton Algorithm for Minimizing Finite Sums}

%
\runningauthor{Jiabin Chen , Rui Yuan , Guillaume Garrigos , Robert M.~Gower}

\twocolumn[

\aistatstitle{SAN: Stochastic Average Newton Algorithm for Minimizing Finite Sums}

\aistatsauthor{%
\hspace{-.3in} Jiabin Chen$^{1,5}$
\And 
\hspace{-.3in} Rui Yuan$^{2,5}$ 
\And 
\hspace{-.3in} Guillaume Garrigos$^3$
\And 
\hspace{-.3in} Robert M.~Gower$^4$ }

\vspace*{3em}
]

\affiliations{%
$^1$Baidu Inc. \'Ecole Polytechnique.
$^2$Meta AI. LTCI, T\'el\'ecom Paris, Institut Polytechnique de Paris.
$^3$Universit\'e de Paris and Sorbonne Université, CNRS, Laboratoire de Probabilités, Statistique et Modélisation.
$^4$CCM, Flatiron Institute. LTCI, T\'el\'ecom Paris, Institut Polytechnique de Paris. 
$^5$Equal contribution.
\vspace*{-0.25em}
}


\doparttoc 
\faketableofcontents 

\begin{abstract}
  We present a principled approach for designing stochastic Newton methods for solving finite sum optimization problems. Our approach has two steps. First, we re-write the stationarity conditions as a system of nonlinear equations that associates each data point to a new row. Second, we apply a Subsampled Newton Raphson method to solve this system of nonlinear equations. 
Using our approach, we develop a new Stochastic Average Newton (SAN) method, which is incremental by design, in that it requires only a single data point per iteration.
It is also cheap to implement when solving regularized generalized linear models, with a cost per iteration of the order of the number of the parameters.
We show through  numerical experiments that SAN requires no  knowledge about the problem, neither
 parameter tuning, while remaining competitive as compared to classical variance reduced gradient methods (e.g. SAG and SVRG), incremental Newton and quasi-Newton methods (e.g. SNM, IQN).
\end{abstract}

\section{Introduction}

Consider the problem of minimizing a sum of terms 
\begin{equation}\label{eq:finite_sum}
w^* \in \underset{w \in \R^d}{\rm{argmin}} \ \frac{1}{n} \sum_{i=1}^n f_i(w) \eqdef f(w),
\end{equation}
where  $f_i$ is a convex twice differentiable loss over a given $i$-th data point. 
When the number of \emph{data points} $n$ and \emph{features} $d$ are large, first order methods such as SGD (Stochastic Gradient Descent), SAG~\citep{SAG}, SVRG~\citep{Johnson2013} and ADAM~\citep{ADAM} are the methods of choice for solving~\eqref{eq:finite_sum} because of their low cost per iteration.
The issue with first order methods is that they can require extensive parameter tuning, and/or knowledge of the parameters of the problem.
Consequently, to make a first order method work well
requires careful  tweaking and tuning from an expert, and a careful choice of the model itself. Indeed, neural networks have evolved in such a way that allows for SGD to converge, such as the introduction of batch norm~\citep{batchnorm} and
the push for more over-parametrized networks which greatly speed-up the convergence of SGD~\citep{MaBB18,VaswaniBS19,SGDstruct}.
 Thus the reliance on first order methods ultimately restricts the choice and development of alternative models.

There is now a concerted effort to develop efficient stochastic second order methods that can exploit the sum of terms structured in~\eqref{eq:finite_sum}. The hope for second order methods for solving~\eqref{eq:finite_sum} is  that they require less parameter tuning and converge for wider variety of models and datasets. In particular, here we set out to develop stochastic second order methods that achieve the following objective.
\begin{objective}{} 
Develop a second order method  for solving~\eqref{eq:finite_sum} that is \emph{incremental},  \emph{efficient}, scales well with the dimension $d$, and that requires no \emph{knowledge from the problem}, neither \emph{parameter tuning}.
\end{objective}

Most stochastic second order methods are not incremental, and thus fall short of our first criteria. 
This is due to the fact that most of these methods are only guaranteed to work in a large mini-batch size regime, and not with a single sample.
For instance, the subsampled Newton methods~\citep{Roosta-Khorasani2016,Bollapragada2018,Newton-MR,Erdogdu2015nips,KohlerL17} 
require potentially large mini-batch sizes in order to guarantee that the subsampled Newton direction closely matches the full Newton direction in high probability. 
Stochastic quasi-Newton methods~\citep{Byrd2011,Mokhtari2014,moritz2016linearly,GowerGold2016},   
SDNA~\citep{Qu2015}, the Newton sketch~\citep{Pilanci2015a} and Lissa~\citep{Lissa},
suffer from the same drawback: the need for large mini-batches or full gradient evaluation to work, which makes them all not incremental.
%
 
The two existing methods that we are aware of that are truly incremental are \emph{IQN} (Incremental Quasi-Newton)~\citep{mokhtari2018iqn,Gao2020IncrementalGB} and \emph{SNM} (Stochastic Newton Method)~\citep{SNM,pmlr-v48-rodomanov16}. 
Both methods also enjoy a fast local convergence rate.
Their only drawback is their computational and memory costs per iteration are at least $\cO(d^2)$ (see Table~\ref{tab:complexity} and Section~\ref{sec:exp_extra} for more details).
This is prohibitive in a setting where the number of parameters for the model is large.
Our goal is to develop a method that is not only incremental, but also has a cost per iteration of  $\cO(d)$, as is the case for first-order methods like SGD.

In this paper we develop two new Newton methods for solving~\eqref{eq:stationarity} that effectively make use of second order information, are incremental, and  are governed by a single global convergence theory. 
Our starting point for developing these methods is to re-write the stationarity conditions\\[-0.2cm]
\begin{equation}\label{eq:stationarity}
\frac{1}{n} \sum_{i=1}^n \nabla f_i(w)  = 0.
\end{equation}
At this point, we could apply Newton's method for solving nonlinear equations, otherwise known as the Newton Raphson method. However, this approach would ultimately require a full pass over the data at each iteration.

To avoid taking full passes over the data,  we re-write~\eqref{eq:stationarity} by introducing $n$ auxiliary variables
$\alpha_i \in \R^d$ and solving instead the nonlinear system given by
\begin{align}
\frac{1}{n}\sum_{i=1}^n \alpha_i & = 0, \label{eq:alphaieqzero} \\
\alpha_i & = \nabla f_i(w), \quad\forall i \in \{1,\ldots, n\}. \label{eq:alphainablafi}
\end{align}
Clearly~(\ref{eq:alphaieqzero}--\ref{eq:alphainablafi}) have the same solutions in $w$. The advantage of~(\ref{eq:alphaieqzero}--\ref{eq:alphainablafi}) 
 is that each gradient lies on a separate row. Consequently, applying a \emph{subsampled} Newton Raphson method, that is sampling a row and then applying Newton Raphson, to (\ref{eq:alphaieqzero}--\ref{eq:alphainablafi}) will result in an incremental method. We refer to~(\ref{eq:alphaieqzero}--\ref{eq:alphainablafi}) as the \emph{function splitting} formulation, since it splits the gradient across different rows. 

To solve~(\ref{eq:alphaieqzero}--\ref{eq:alphainablafi}) efficiently, we propose SAN (\emph{Stochastic Average Newton}) in Section~\ref{sec:SAN}. 
SAN is a  subsampled Newton Raphson method that is based on a new variable metric extension of \emph{SNR} (Sketch Newton Raphson Method)~\citep{Yuan2020sketched} that we present in Section~\ref{sec:SNR}, which is itself a nonlinear extension of the Sketch-and-Project method for solving linear systems~\citep{Gower2015}. 
By using a different subsampling of the rows~(\ref{eq:alphaieqzero}--\ref{eq:alphainablafi}), we also derive SANA in Section~\ref{sec:SANA}, which is a variant of SAN that uses unbiased estimates of the gradient.

Note that the idea of applying a subsampled Newton method to a well-chosen system of optimality conditions is not new.
Indeed, it was recently shown in~\citep{Yuan2020sketched} that the SNM method~\citep{SNM} can be seen as the application of a subsampled Newton method to the equations\\[-0.55cm]
\begin{align}
 \frac{1}{n}\sum_{i=1}^n \nabla f_i(\alpha_i) &= 0, \nonumber \\
 w 
& = \alpha_i, \quad \mbox{ for } i=1,\ldots, n,\label{eq:weqai}
\end{align}
which clearly are equivalent to \eqref{eq:stationarity}.
Consequently, the two methods SAN and SNM are both subsampled Newton Raphson methods applied to either a function or a \emph{variable splitting} formulation of~\eqref{eq:stationarity}.

The contributions of our paper are the following:
\begin{itemize}[parsep=0em, topsep=0em]
\item We propose  combining the function splitting reformulation~(\ref{eq:alphaieqzero}--\ref{eq:alphainablafi}) with a subsampled Newton Raphson as a tool for designing stochastic Newton methods, all of which are \emph{variance reducing} in the sense that they are incremental and they converge with
 a constant step size.
\item We introduce SAN (Stochastic Average Newton method) by using this tool, which is incremental and parameter-free, in that, SAN works well with a step size $\gamma =1$ independently of the underlying dataset or the objective function.
\item By specializing  to GLMs (Generalized Linear Models), we develop an efficient implementation of SAN that has the same cost per iteration as the first-order methods. We perform extensive numerical experiments and show that SAN is competitive as compared to SAG and SVRG.
\item To provide a convergence theory of our methods, we extend the class of Sketched Newton Raphson~\citep{Yuan2020sketched} methods to allow for a variable metric that includes SAN as a special case.
\end{itemize}

In Section~\ref{sec:ROW}, we show how to derive the SAN and  SANA methods.
 We then present two different experimental settings comparing  SAN to variance reduced gradients methods in Section~\ref{sec:exp}.
In Section~\ref{sec:SNR}, we study SAN/SANA as instantiations of a new \emph{variable metric} Sketch Newton Raphson  method and present a convergence theory for this class of method.


\text{The following will be assumed throughout the paper.}\vspace*{-1em}
\begin{assumption}\label{Ass:strict convexity of problem}
For all $i\in \{1, \dots, n \}$, the function $f_i : \mathbb{R}^d \longrightarrow \mathbb{R}$ is of class $C^2$ and verifies $\nabla^2 f_i(w) \succ 0$ for every $w \in \mathbb{R}^d$. 
\end{assumption}

\section{Function splitting methods} \label{sec:ROW}

The advantage of the function splitting formulation given by~(\ref{eq:alphaieqzero}) and (\ref{eq:alphainablafi}) is that
there is a separate row for each data point. We will now take advantage of this, and develop new incremental Newton methods based on subsampling the rows of~(\ref{eq:alphaieqzero}--\ref{eq:alphainablafi}).

The reformulation given in~(\ref{eq:alphaieqzero}--\ref{eq:alphainablafi}) is a large system of nonlinear equations. 
For brevity, let $p:=(n+1)d$ and 
 $x = \begin{bmatrix}
w \ ;
\alpha_1 \ ;
\cdots \ ;
\alpha_n
 \end{bmatrix}\in \R^{p}$ 
 be the stacking\footnote{In this paper vectors are columns by default, and given $x_1, \dots, x_n \in \mathbb{R}^q$ we note $[x_1 ; \dots ; x_n]\in \mathbb{R}^{qn}$ the (column) vector stacking the $x_i$'s on top of each other.} of the $w$ and $\alpha_i$ variables.
Thus solving~(\ref{eq:alphaieqzero}--\ref{eq:alphainablafi}) is equivalent to solve $F(x)=0$, where
\begin{align}\label{eq: def_F}
F : & \ \R^{p} \rightarrow  \R^{p} \\
\nonumber & \ x \mapsto  \begin{bmatrix}
\frac{1}{n}\sum\alpha_i ;
\nabla f_1(w) - \alpha_1 ;
\cdots ;
\nabla f_n(w) - \alpha_n
\end{bmatrix}.
\end{align}

Solving nonlinear equations has long been one of the core problems in numerical analysis, with variants of the Newton Raphson method~\citep{Ortega:2000} being  one of the core techniques.
From a given iterate $x^k \in \R^{p}$, the Newton Raphson method computes the next iterate $x^{k+1}$ by linearizing $F$ around $x^k$ and solving the \emph{Newton system}
\begin{equation} \label{eq: newton system}
\nabla F(x^k)^\top (x^{k+1}-x^k)  = -F(x^k). 
\end{equation}
Here $\nabla F(x) \in  \R^{p \times p}$ denotes the Jacobian matrix of $F$ at $x$, and it is assumed that \eqref{eq: newton system} has a solution.
The least norm solution of
the Newton system
is given by
\begin{equation}\label{eq:newton}
x^{k+1} = x^k -  \nabla F(x^k)^{\top \dagger} F(x^k), 
\end{equation}
where $^\dagger$ denotes
the Moore-Penrose pseudoinverse.

This update can also be written as a projection  step:
\begin{align}\label{D:SNR projection+relaxation}
x^{k+1} = \ &  \argmin \Vert x - x^k \Vert^2 \nonumber \\ 
&\text{ s.t. } \nabla F(x^k)^\top (x-x^k) = - F(x^k).
\end{align}
%
 In our setting,~\eqref{eq:newton} is prohibitively expensive because it requires access to all of the data at each step and the solution of a large $(n+1)d \times (n+1)d$ linear system. To bring down the cost of each iteration, and to have a resulting incremental method, at each iteration we will \textit{subsample} the rows of the Newton system before taking a projection step. Next, we present two methods based on subsampling. Later on Section~\ref{sec:SNR}, we generalize this subsampling approach to make use of \emph{sketches} of the system.


\subsection{SAN: the Stochastic Average Newton method} 

\label{sec:SAN}

The SAN method is a subsampled Newton Raphson method that alternates between sampling  equation~\eqref{eq:alphaieqzero} or sampling one of the equations in~\eqref{eq:alphainablafi}. 
After sampling, we then apply a step of Newton Raphson to the sampled equation.

To detail the SAN method, let $\pi \in (0 , 1)$ be a fixed probability, and let  
$x^k=[w^k;\alpha_1^k; \dots; \alpha_n^k] \in \R^p$
be a given $k$-th iterate.
 With probability $\pi$
the SAN method samples equation~\eqref{eq:alphaieqzero} and focuses on finding a solution to this equation. Since~\eqref{eq:alphaieqzero} is a linear equation, it is equal to its own Newton equation.
Furthermore, this linear equation~\eqref{eq:alphaieqzero}
has  $n$ variables and only one equation, thus it has infinite solutions.  We choose a single one of these infinite solution by using a projection step
\begin{align}
 \alpha^{k+1}_1, \dots, \alpha^{k+1}_n = & \underset{\alpha_1, \dots, \alpha_n \in \R^{d}}{\argmin} \textstyle\sum_{i=1}^n \norm{\alpha_i - \alpha_i^k}^2
 \nonumber \\ 
&\mbox{ s.t. }\textstyle\frac{1}{n}\sum_{i=1}^n \alpha_i = 0. \label{eq:SAN implicit linearized 0}
\end{align}
The solution to this projection is given in line~\ref{ln:rsn_mean} in Algorithm~\ref{algo:SAN} when $\gamma =1$. We have added the  step size  $\gamma \in (0 , \; 1]$ to act as relaxation.


Alternatively, with probability $(1-\pi)$ the SAN method then samples the $j$-th equation   in~\eqref{eq:alphainablafi} uniformly among the $n$ equations.
To get the Newton system of $\nabla f_j(w) =\alpha_j$, we linearize around  $w^k \in \R^d$ and 
$\alpha_j^k \in \R^d$  and set the linearization to zero giving
\[  \nabla f_j(w^k)+\nabla^2 f_j(w^k)(w-w^k) = \alpha_j .\]
This linear equation has $2d$ unknowns, and thus also has infinite solutions. Again, we use a projection step to pick a unique solution as follows
\begin{align}
\alpha^{k+1}_j, w^{k+1} &= \underset{\alpha_j, w \in \R^d}{\argmin} \norm{\alpha_j - \alpha_j^k}^2 + \norm{w-w^k}_{\nabla^2 f_j(w^k)}^2 \label{eq:SAN implicit linearized 1..n} \nonumber  \\  
&\mbox{ s.t. } \nabla f_j(w^k)+\nabla^2 f_j(w^k)(w-w^k) = \alpha_j .
\end{align}
Here we have introduced a projection under a norm 
$ \norm{w}_{\nabla^2 f_j(w^k)}~\eqdef~\dotprod{\nabla^2 f_j(w^k)w,w}$ which is based on the Hessian matrix $\nabla^2 f_j(w^k).$
Performing a projection step with respect to the metric induced by the Hessian is often used in Newton type methods such as interior point methods~\citep{Renegar:2001} and quasi-Newton methods~\citep{Goldfarb1970}.
Moreover, we observed that this choice of metric  resulted in a much faster algorithm (see Section~\ref{sec:SANvsSANiI} for experiments that highlight this). 
The closed form solution to the above is given in lines~\ref{ln:rsn_d}-\ref{ln:rsn_alpha} in Algorithm~\ref{algo:SAN} when $\gamma=1$ (see Lemma \ref{L:SAN implicit to explicit} for the details). 

We gather all these updates in Algorithm~\ref{algo:SAN} and call the resulting method the \emph{Stochastic Average Newton} method, or \emph{SAN} for short. 

\begin{algorithm}
\caption{SAN: Stochastic Average Newton}
\label{algo:SAN}
\KwIn{$\{f_i\}_{i=1}^n$, step size $\gamma \in (0,1]$, probability $\pi \in (0,1)$, max iteration $T$} 
Initialize $\alpha^0_1, \cdots, \alpha^0_n, w^0 \in \R^d$ 

\For{$k=1, \dots, T$}{
    \With{
    $ \displaystyle \alpha_i^{k+1} \;=\; \alpha_i^k - \frac{\gamma}{n}\sum_{j=1}^n\alpha_j^k, \  \forall i \in \{1, \cdots, n\} $ 
    \label{ln:rsn_mean}
    } 
    \Otherwise{
    Sample uniformly $j \in \{1, \cdots, n\}$

    $\mH_k \;=\; \mI_d + \nabla^2f_j(w^k)$\label{ln:rsn_H}

    $d^{k} \;=\;  -  \mH_k^{-1}\left(\nabla f_j(w^k) - \alpha_j^k \right) $\label{ln:rsn_d} 

    $ \displaystyle w^{k+1} \;=\; w^k + \gamma d^k$ \label{ln:rsn_x}

    $ \displaystyle\alpha_j^{k+1} \;=\; \alpha_j^k - \gamma d^k $ \label{ln:rsn_alpha}
    }
}
\KwOut{Last iterate $w_{T+1}$}
\end{algorithm}

The SAN method is incremental, since it can be applied with as little as one data point per iteration. 
SAN can also be implemented in such a way that the cost per iteration is $\cO(d)$ in expectation.
 Indeed, 
the \emph{averaging step} on line~\ref{ln:rsn_mean} contributes with a  $\pi \times \cO( nd)$ cost to the total cost in expectation, since all of the vectors $\alpha_i\in\R^d$ for $i=1,\ldots, n,$ are updated. But 
as we
 found through expensive testing in Section~\ref{sec:gridsearch}, SAN converges quickly if $\pi$ is of the order of $\cO(1/n)$, reducing the cost in expectation to $\cO(d)$.
Further, the average of the $\alpha_i$'s can be efficiently implemented by maintaining and updating a variable $\overline{\alpha}^k  =\frac{1}{n}\sum_{j=1}^n\alpha_j^k$.
The main cost for SAN is in solving the linear system $(\mI_d + \nabla^2f_j(w^k)) d = \alpha_j^k - \nabla f_j(w^k).$ 
Solving this system with a direct solver would cost $\cO(d^3)$. Alternatively, the solution can be approximated using an iterative Krylov method for which  each iteration costs $\cO(d)$ by using 
 backpropagation~\citep{Freund1992,Christianson:1992} to compute Hessian-vector products.
For regularized generalized linear models (GLMs), the total cost of this matrix inversion is only $\cO(d)$ operations, as we show next.

\paragraph{Generalized Linear Models.}  Regularized GLMs are models for which we have
\begin{eqnarray} \label{eq:glm}
f_i(w) &=& \phi_i(a_i^\top w) + \lambda R(w),
\end{eqnarray}
where $\phi_i:\R \rightarrow \R$ is a loss function associated with the $i$-th data point $a_i \in \R^d$, $\lambda > 0$ is a regularization parameter and $R$ is a regularizer that is twice differentiable and separable, i.e. 
 $R(w) = \sum_{i=1}^d R_i(w_i)$ with $R_i:\R \rightarrow \R$. 
The inversion on line~\ref{ln:rsn_d} of Algorithm~\ref{algo:SAN} can be efficiently computed using the Woodbury identity because the Hessian $\nabla^2 f_j(w) = \phi_j'' (a_j^\top w)(a_j a_j^\top ) + \lambda \nabla^2R(w)$ is a rank-one perturbation of a diagonal matrix,  which costs $\cO(d)$ to invert (see Lemma~\ref{L:GLM hessians} for an explicit formula).
\begin{remark}[SAN vs. SNM for GLMs]
SAN can be implemented efficiently for all GLMs with separable regularizers.
This is not the case for SNM~\citep{SNM}, which can only be implemented efficiently when the regularizer is the L2 norm.
 For other separable regularizers, 
the cost per iteration for SNM is to $\cO(d^3)$ instead of $\cO(d^2)$.
See Section~\ref{sec:exp_extra} in the supp. material for details.
\end{remark}

\subsection{SANA: alternative with simultaneous projections}
\label{sec:SANA}
Here we present SANA, an alternative version of the SAN method. Instead of alternating between projecting onto linearizations of~\eqref{eq:alphaieqzero} and~\eqref{eq:alphainablafi}, the SANA method  projects onto the  intersection of~\eqref{eq:alphaieqzero} and the linearization of a subsampled equation~\eqref{eq:alphainablafi}.
In other words, the next iterate $x^{k+1} = [w^{k+1}; \alpha_1^{k+1}; \dots; \alpha_n^{k+1}]$ is defined as the unique solution of
\begin{align}
\underset{w, \alpha_1, \dots, \alpha_n \in \R^d}{\argmin}& \sum\limits_{i=1}^n \norm{\alpha_i - \alpha_i^k}^2 + \norm{w-w^k}_{\nabla^2 f_j(w^k)}^2, \nonumber \\
\mbox{ s.t. } &\nabla f_j(w^k)+\nabla^2 f_j(w^k)(w-w^k) = \alpha_j,  \notag \\
\phantom{\mbox{ s.t. }}&\textstyle\frac{1}{n}\sum_{i=1}^n \alpha_i \; =\; 0.  \label{eq:SANA implicit linearized}
\end{align}

The closed form solution  of~\eqref{eq:SANA implicit linearized} corresponds to lines~\ref{ln:SANA0}--\ref{ln:SANA4} in Algorithm~\ref{algo:SANA} when the relaxation parameter is $\gamma = 1$ (see Lemma \ref{L:SANA implicit to explicit} for a proof). 
%


\begin{algorithm}
\caption{SANA}
 \label{algo:SANA}
\KwIn{$\{f_i\}_{i=1}^n$, step size $\gamma \in (0,1]$, max iteration $T$} 
    Initialize $w^0 , \alpha^0_1, \cdots, \alpha^0_n \in \R^d$ s.t. $\textstyle\sum_{i=1}^n \alpha_i^0 = 0$ 

    \For{$k=1, \dots, T$}{
    Sample uniformly $j \in \{1,\dots,n\}$ 

    $\mH_k \;=\; (1-\tfrac{1}{n}) \mI_d +  \nabla^2 f_j(w^k)$\label{ln:SANA0}

    \text{$\textstyle d^k \;=\; -\mH_k^{\text{\tiny $-1$}} (  \nabla f_j(w^k)-\alpha_j^k)$} \label{ln:SANA1}

    $ \displaystyle w^{k+1} \;=\; w^k +  \gamma d^k $ \label{ln:SANA2}

    $ \alpha_j^{k+1} \;=\; \alpha_j^k - \gamma (1-\tfrac{1}{n})d^k $ \label{ln:SANA3}

    $ \alpha_i^{k+1} \;=\; \alpha_i^k + \tfrac{\gamma}{n}d^k, \quad \text{ for } i \neq j$ \label{ln:SANA4}
    }
\KwOut{Last iterate $w_{T+1}$}
\end{algorithm}



Computing one step of this method requires  access to only one function $f_j$ (through its gradient and Hessian evaluated at $w^k$).
In terms of computational cost, each step requires  inverting the $d \times d$ matrix  $(1 - \frac{1}{n})I_d + \nabla^2 f_j(w^k)$.
As with SAN, this cost reduces to $\cO(d)$ in the context of generalized linear models. See Algorithm~\ref{algo:SANA_GLM}  in the appendix for the resulting implementation for GLMs. 
Yet even in the case of GLMs, the SANA method costs $\cO(nd)$ per iteration because it updates all  the $\alpha_i$ vectors at every iteration. Thus the SANA method has complexity which is $\cO(n)$ times larger than SAN in expectation.

Both SAN and SANA can be interpreted as a stochastic relaxed Newton method that uses  estimates of the gradient.  Indeed, computing $d^k$  in  Algorithms~\ref{algo:SAN} and~\ref{algo:SANA} requires solving a relaxed Newton system 
\begin{equation}\label{eq:saeske9ks}
 \left(\delta \mI_d +  \nabla^2 f_j(w^k) \right) d^k = \alpha_j^k -\nabla f_j(w^k),
\end{equation}
where $\delta =1$ and $\delta = 1-\frac{1}{n},$ respectively.
The right hand side of this Newton system is a biased estimate of the gradient for SAN and an unbiased estimate for SANA. To see this, for simplicity, let $\gamma =1$. Taking expectation conditioned on time $k$ over the right-hand side  of~\eqref{eq:saeske9ks} gives
\[\E{ \alpha_j^k -\nabla f_j(w^k)} = \frac{1}{n} \sum_{i=1}^n \alpha_i^k - \nabla f(w^k). \]
For SANA, because the averaging constraint is always enforced in~\eqref{eq:SANA implicit linearized}, we have that $ \frac{1}{n} \sum_{i=1}^n \alpha_i^k = 0$, thus the right hand side of~\eqref{eq:saeske9ks} is always an unbiased estimate of the negative gradient. As for SAN, the averaging constraint  is only enforced every so often with the update on line~\ref{ln:rsn_mean} in Algorithm~\ref{algo:SAN}. Thus for SAN, the right hand is a biased estimate of the negative gradient until the averaging constraint is enforced. In this sense, SAN and SANA are analogous to SAG~\citep{SAG} and SAGA~\citep{SAGA_Nips}.
 We found in practice that this biased estimate of the gradient did not hurt the empirical performance of SAN, and thus we focus on experiments on SAN  in Section~\ref{sec:exp}.
 

%

\begin{figure*}[t]
\centering
\includegraphics[width=.24\linewidth]{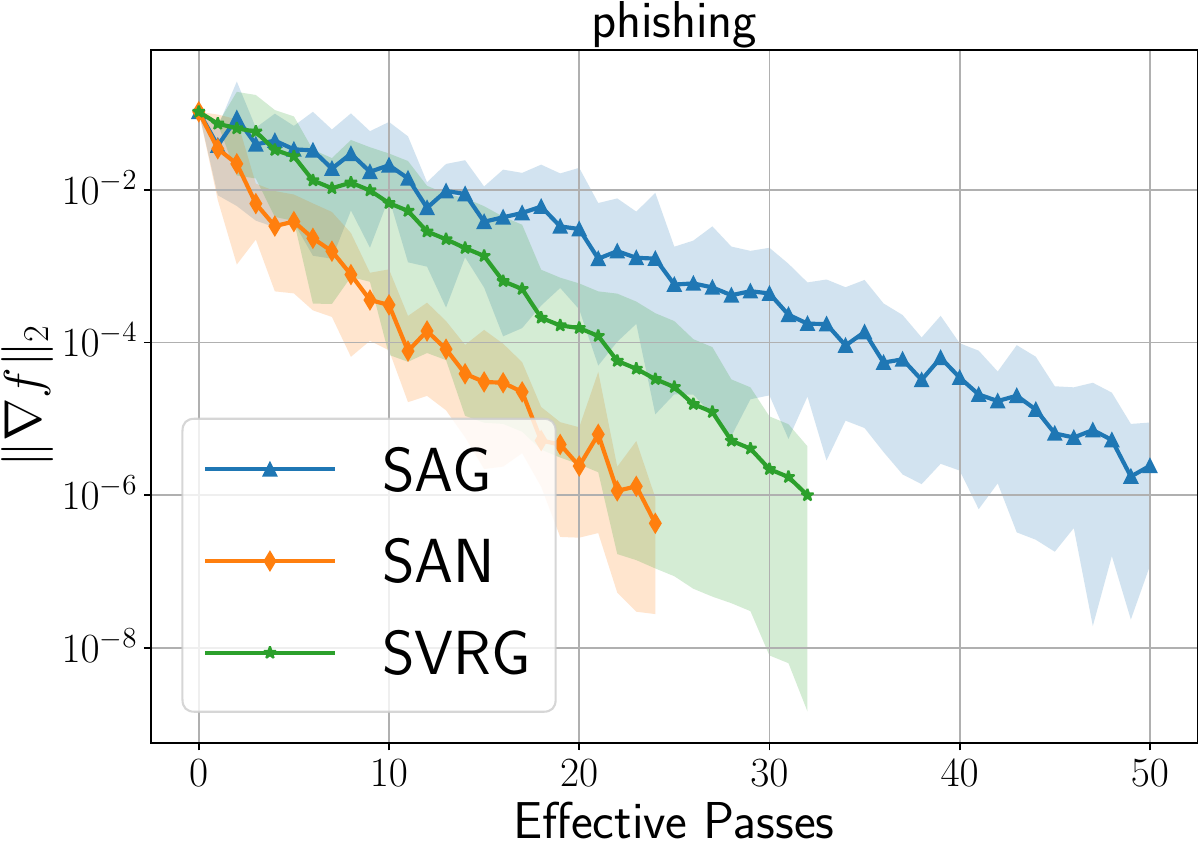}
\includegraphics[width=.24\linewidth]{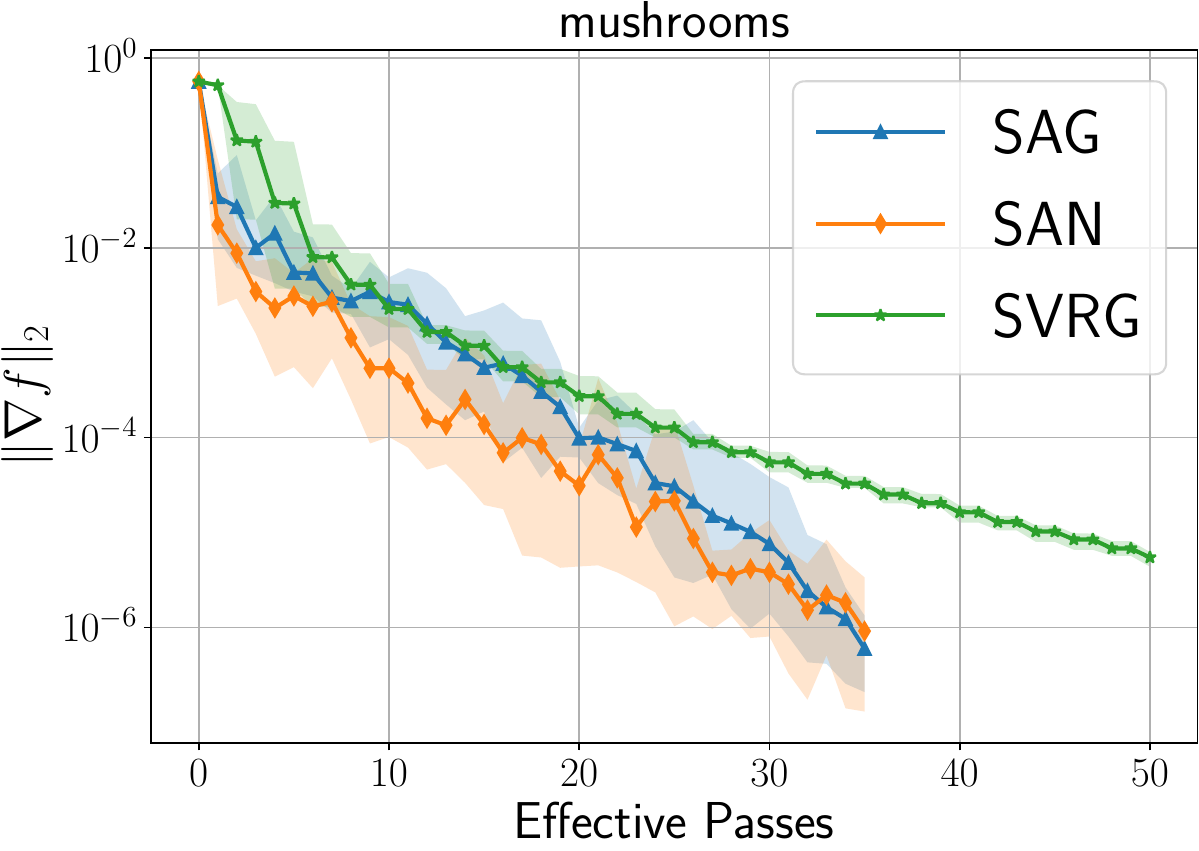}
\includegraphics[width=.24\linewidth]{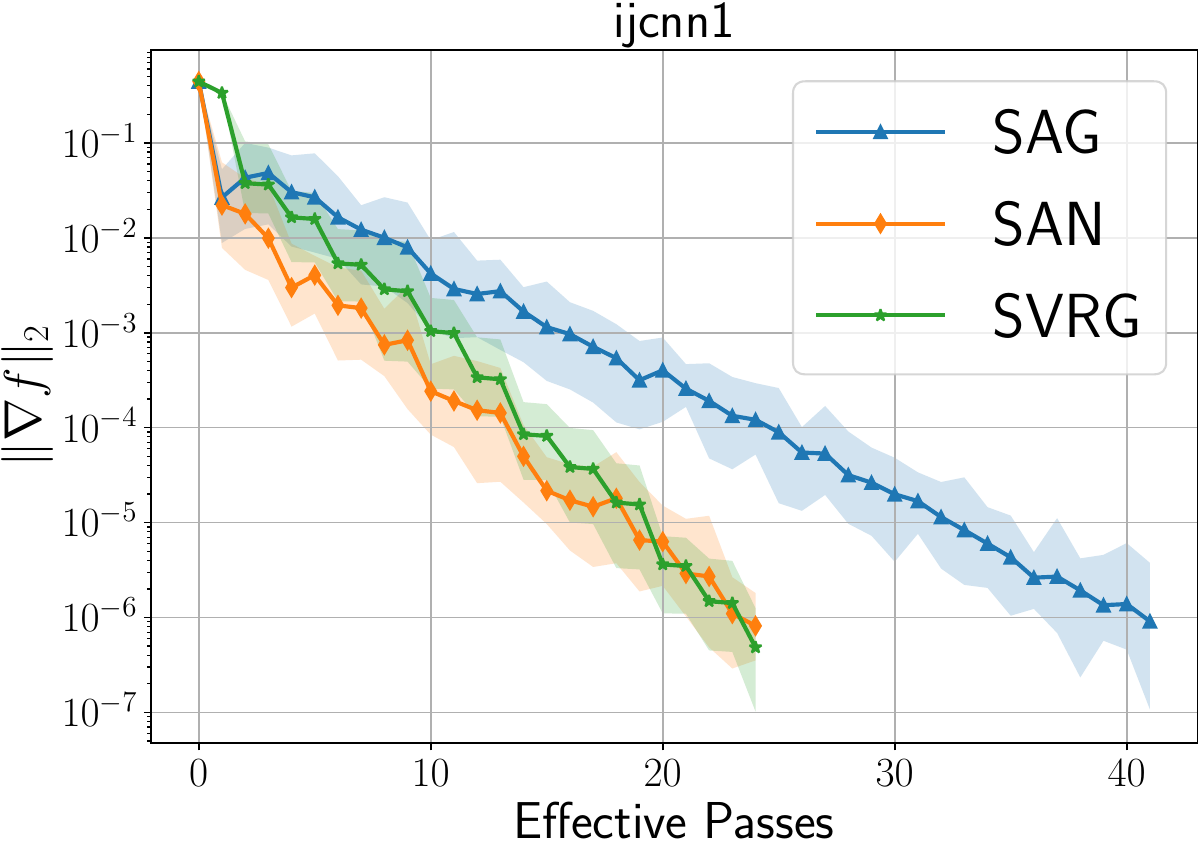}
\includegraphics[width=.24\linewidth]{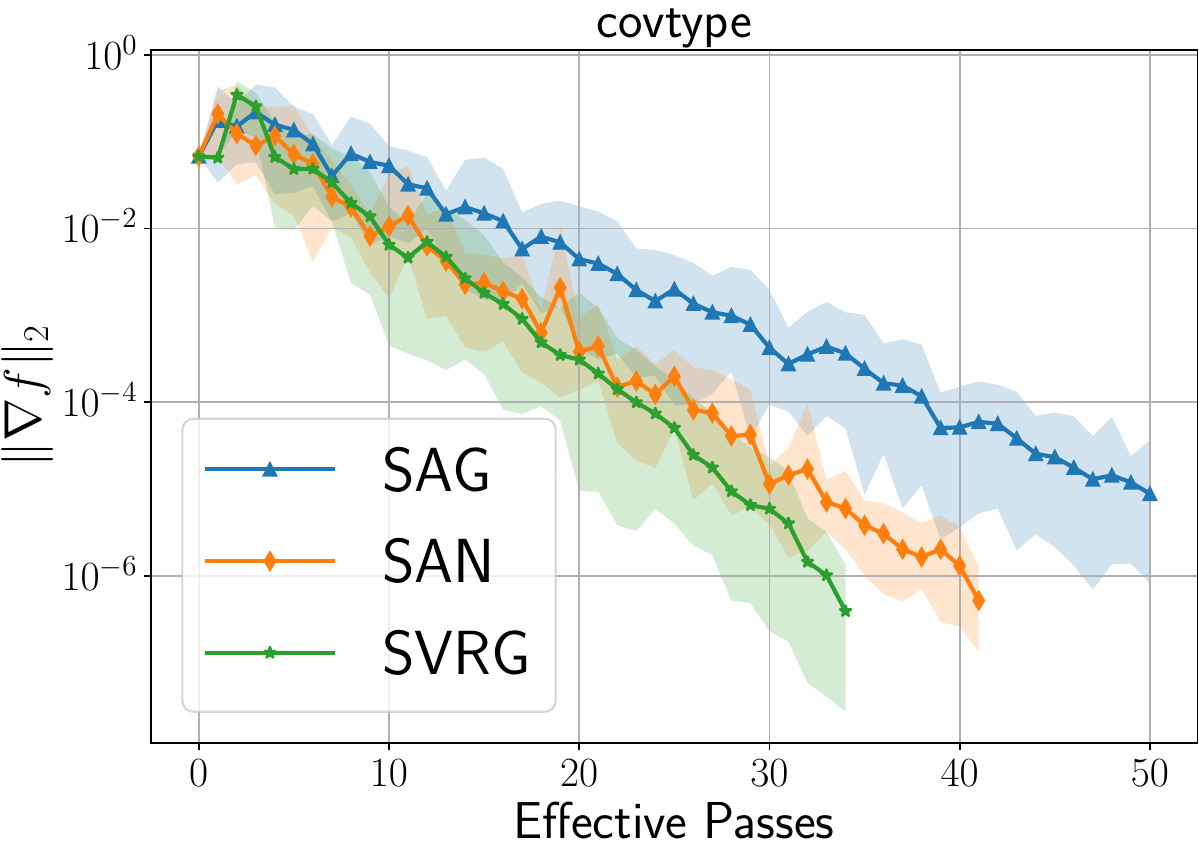} 
\includegraphics[width=.24\linewidth]{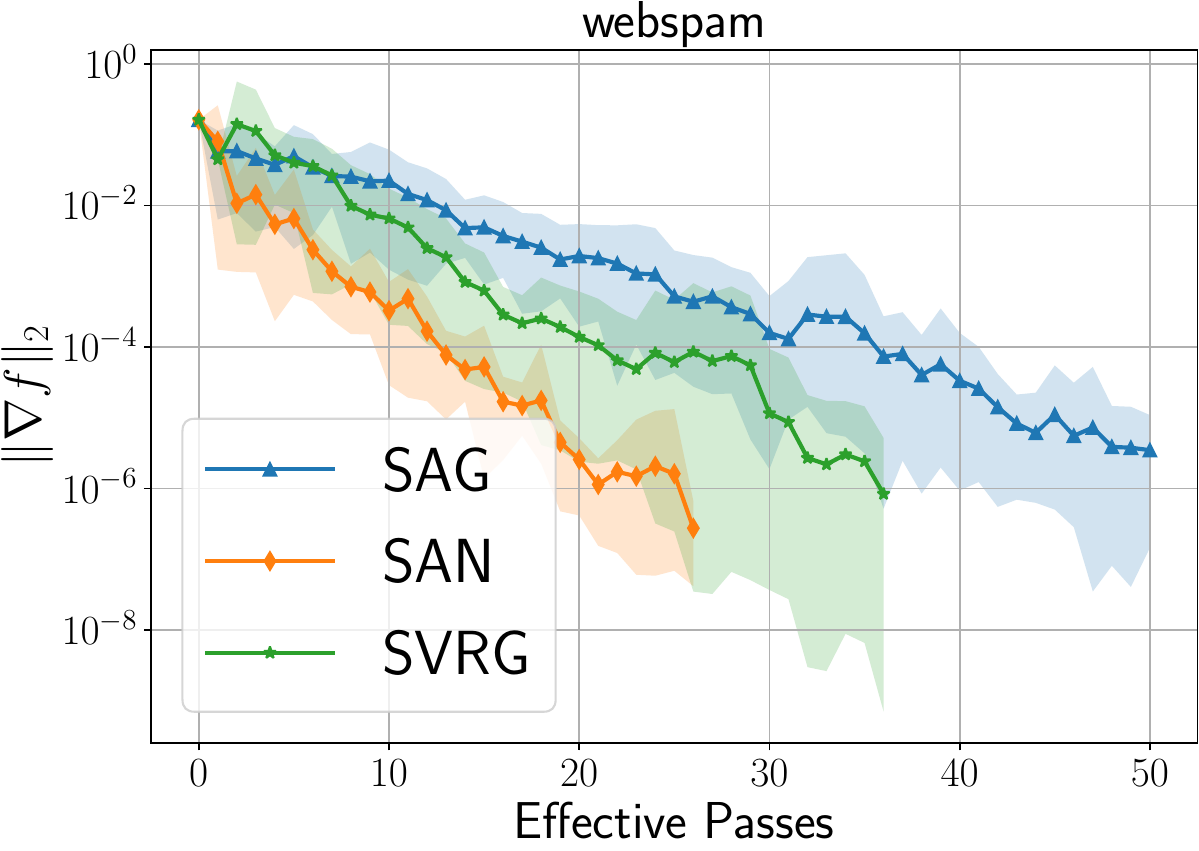}
\includegraphics[width=.24\linewidth]{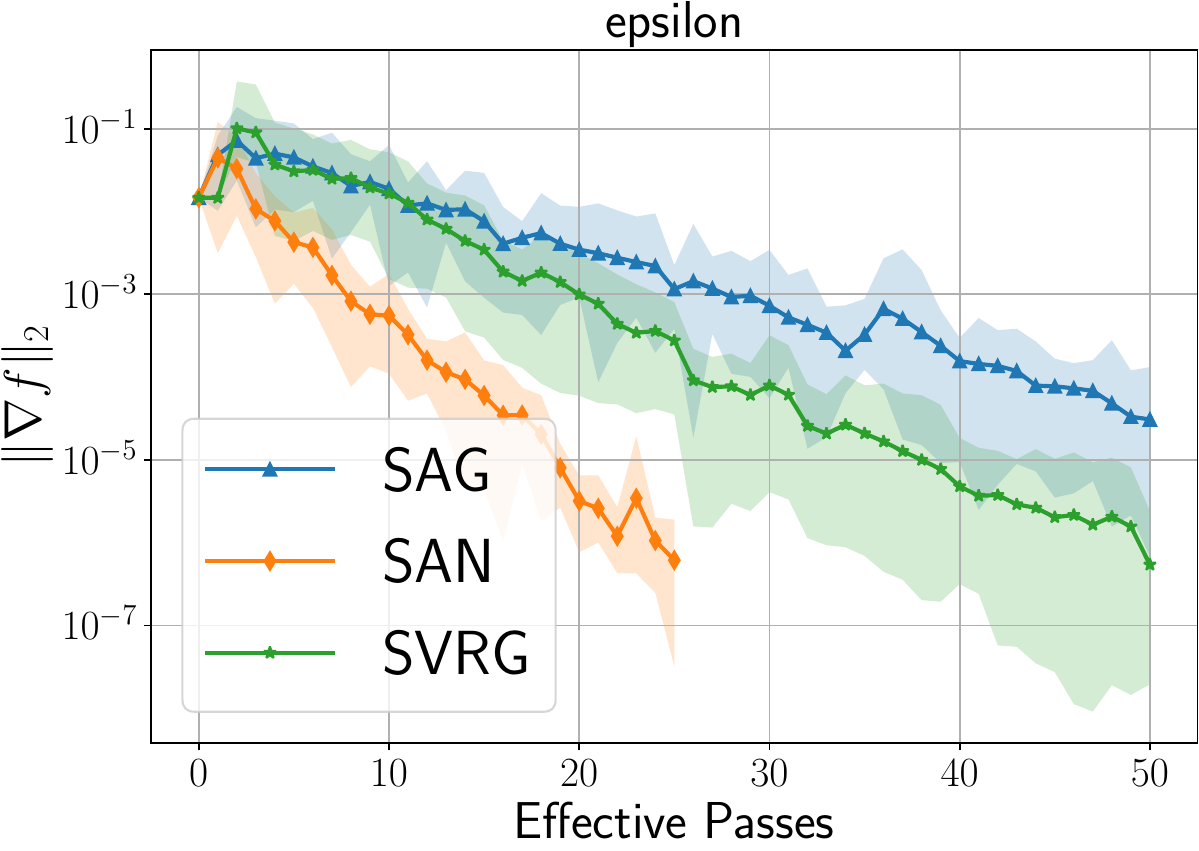}
\includegraphics[width=.24\linewidth]{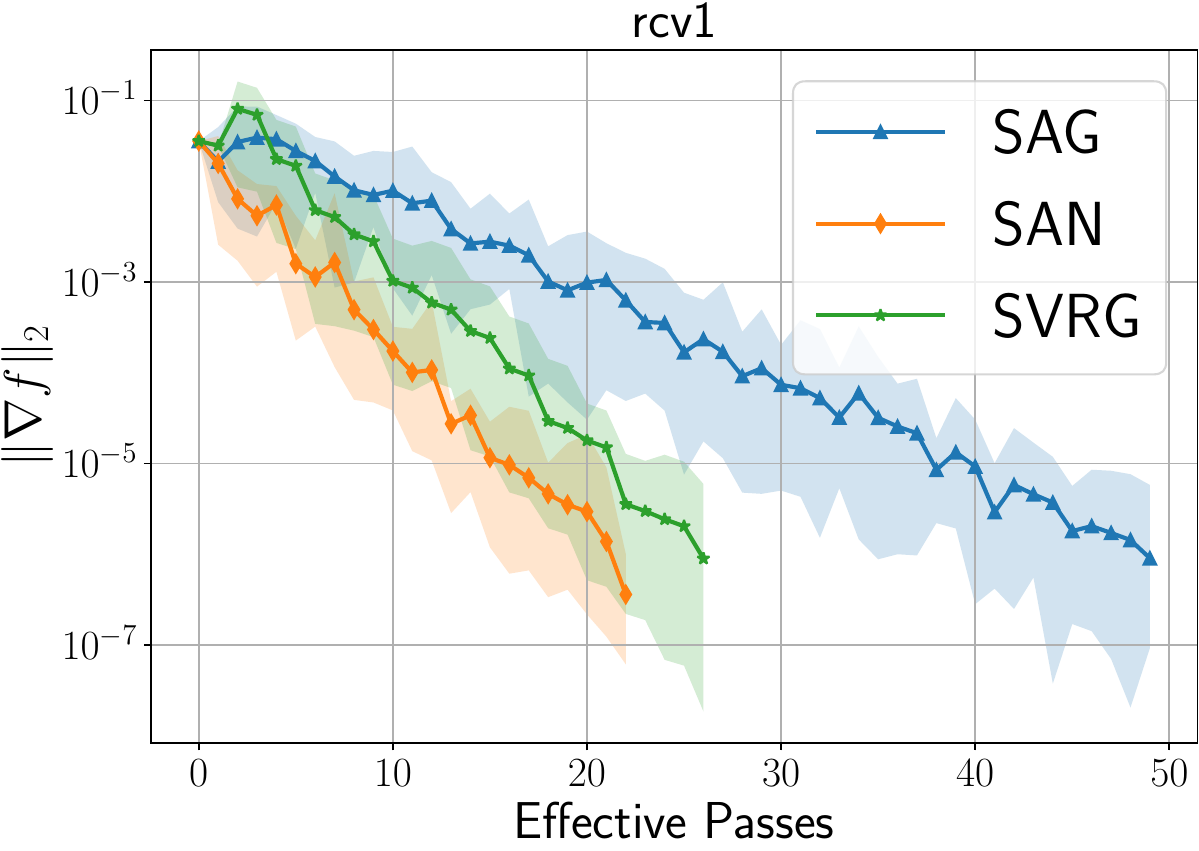}
\includegraphics[width=.24\linewidth]{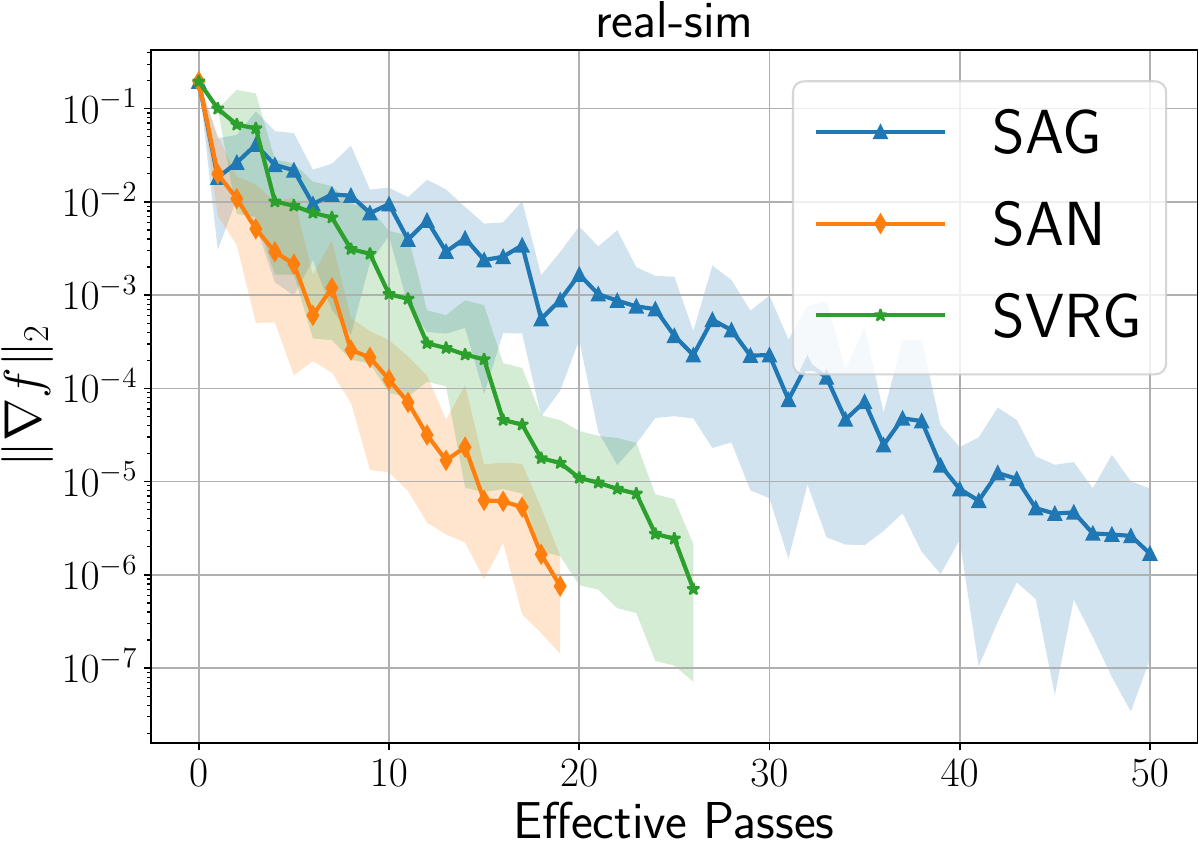}
\caption{Logistic regression with L2 regularization.}
\label{fig:L2SAN}
\end{figure*}

\begin{figure*}[t]
\centering
\includegraphics[width=.24\linewidth]{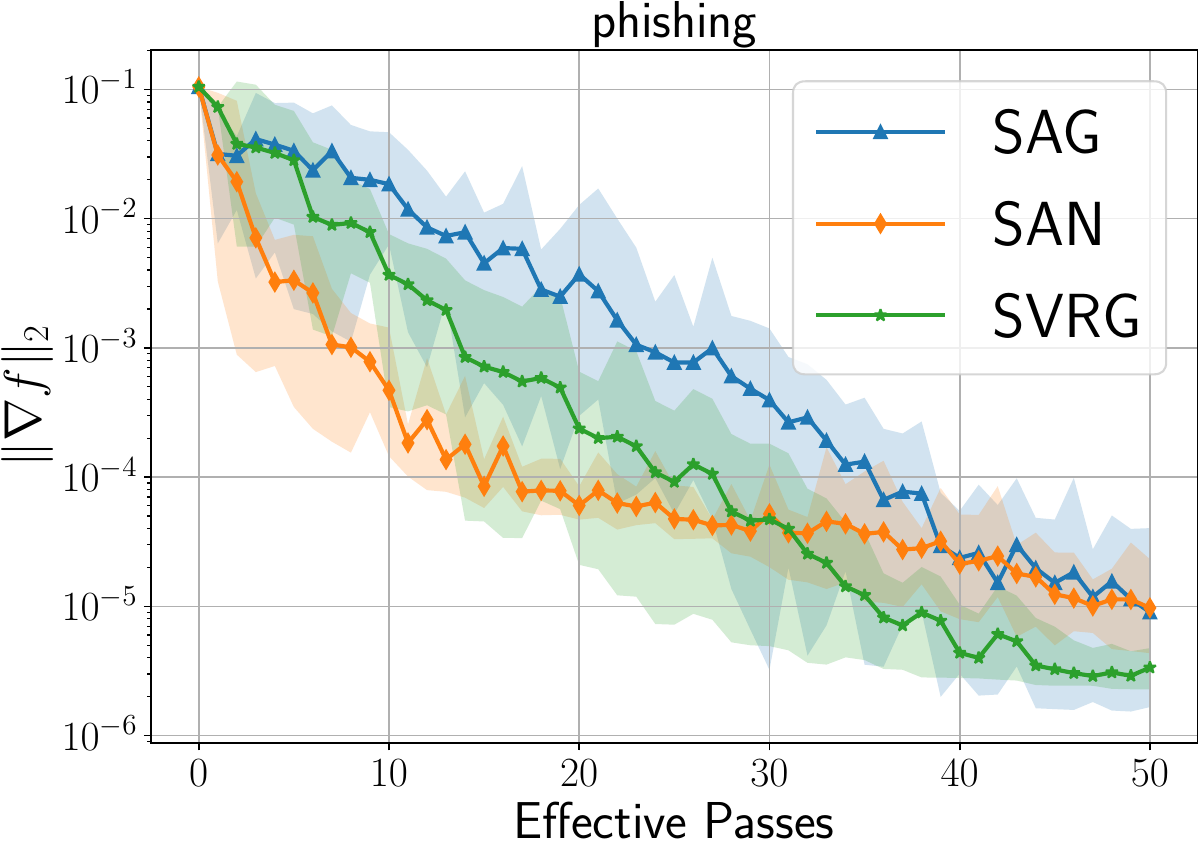} 
\includegraphics[width=.24\linewidth]{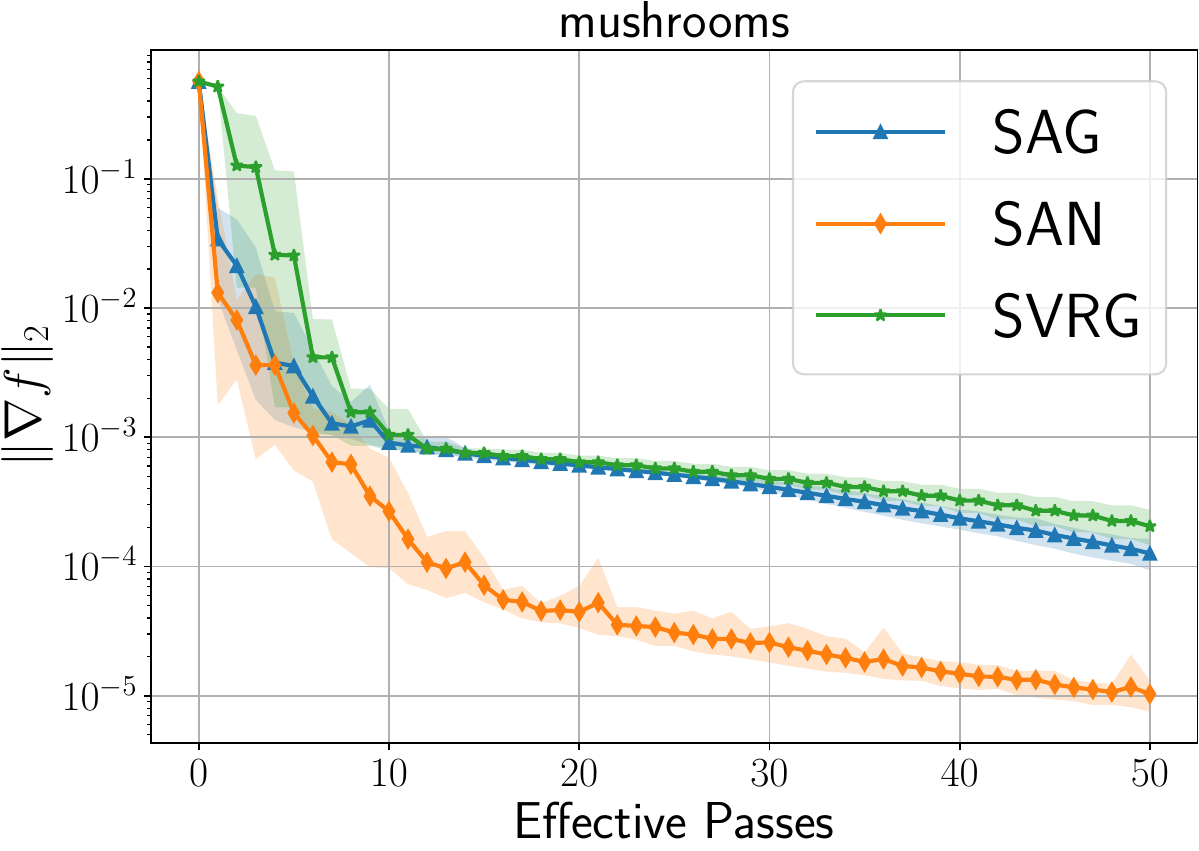}
\includegraphics[width=.24\linewidth]{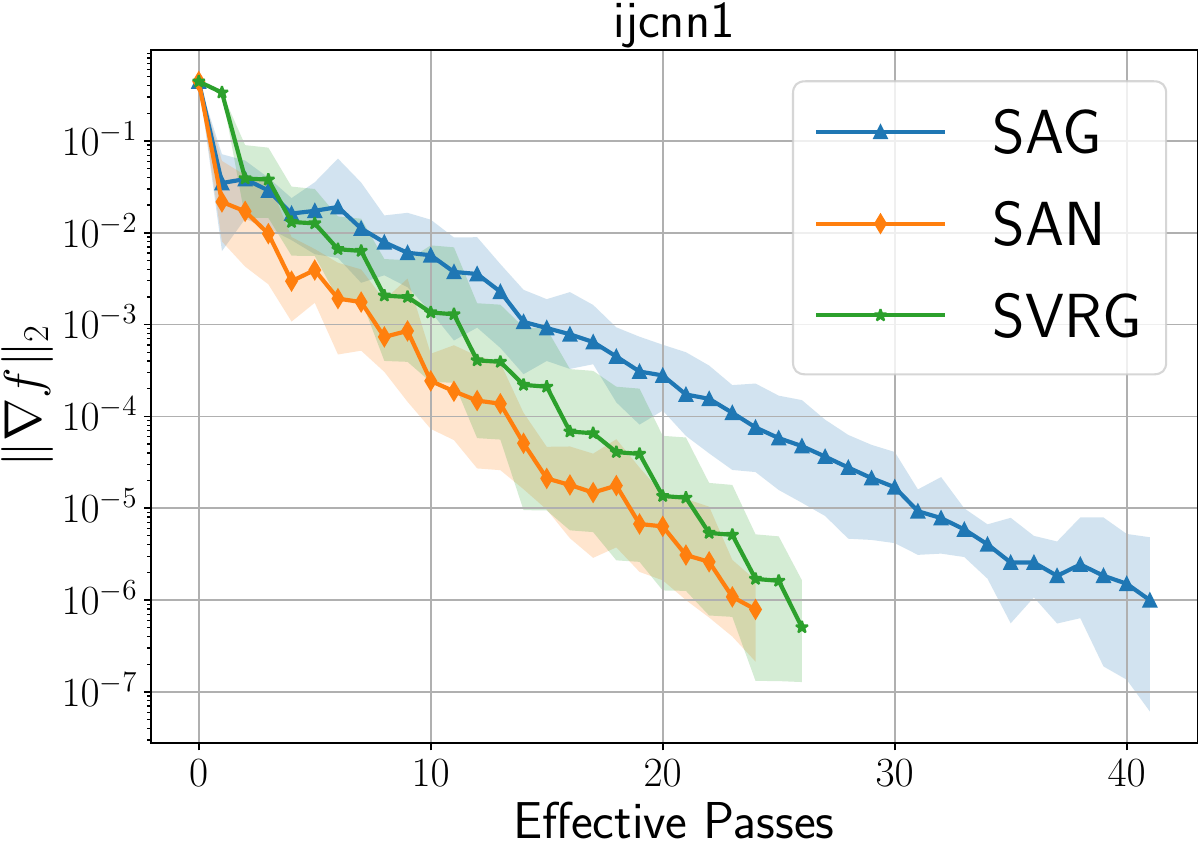}
\includegraphics[width=.24\linewidth]{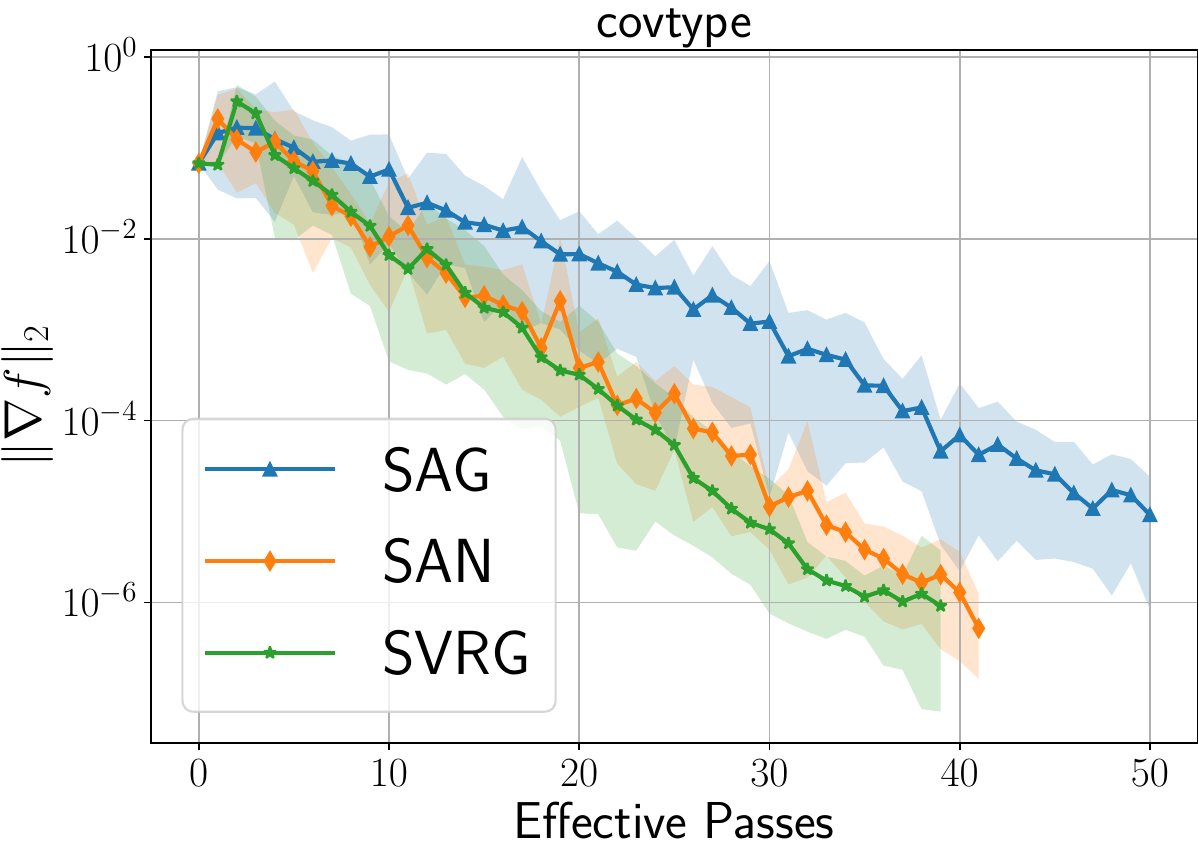} 
\includegraphics[width=.24\linewidth]{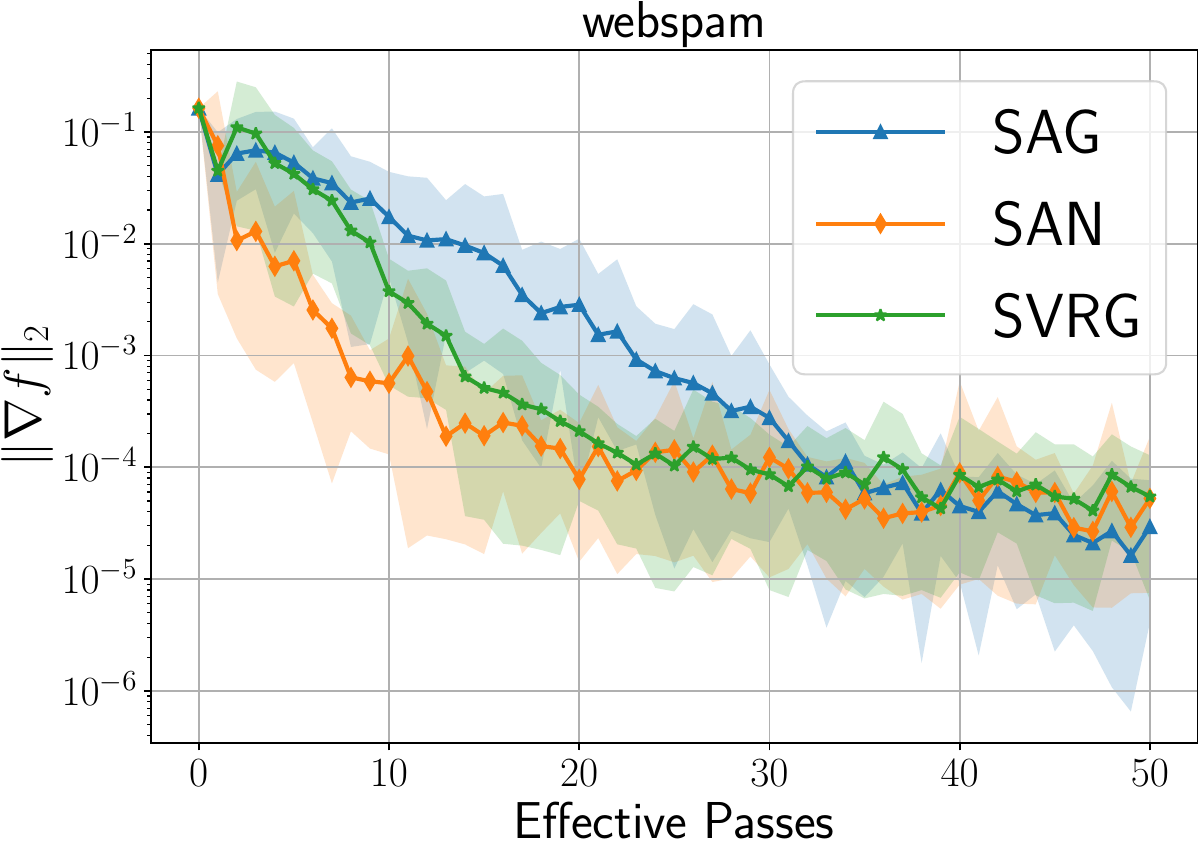}
\includegraphics[width=.24\linewidth]{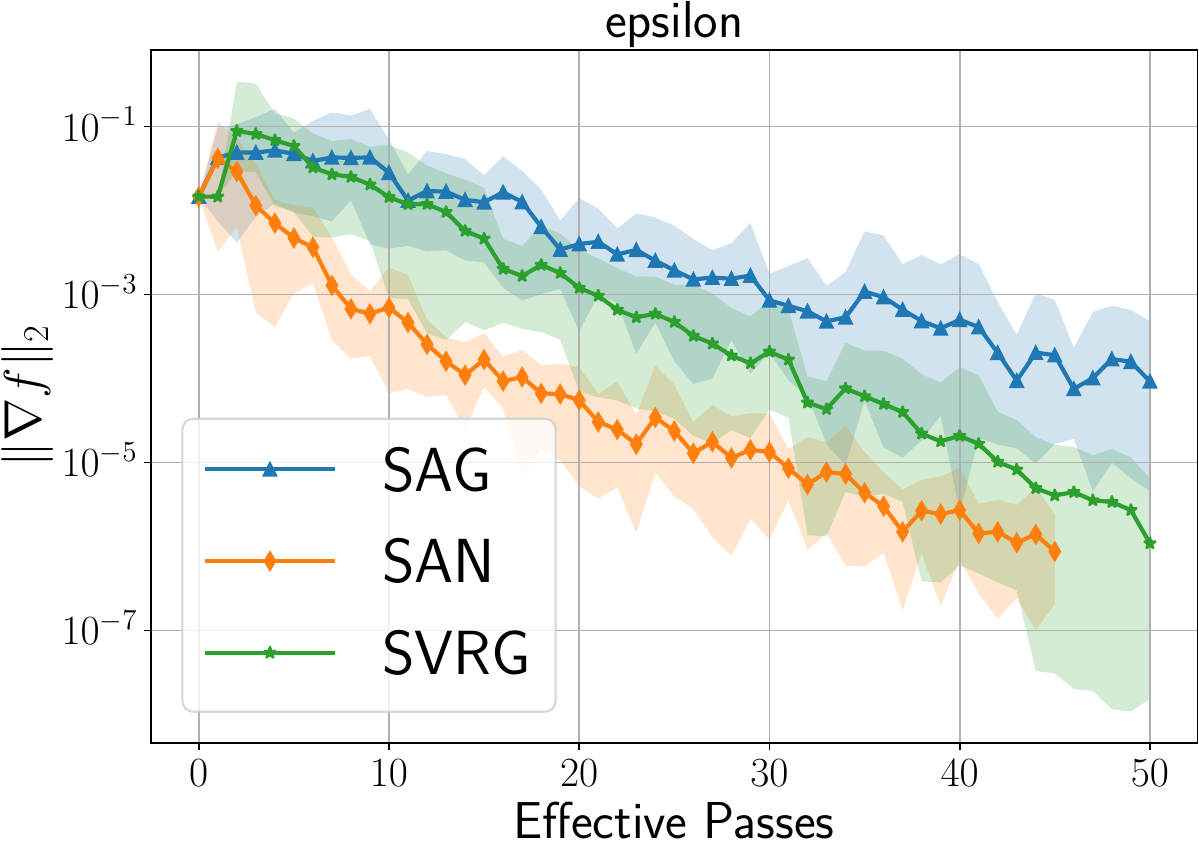}
\includegraphics[width=.24\linewidth]{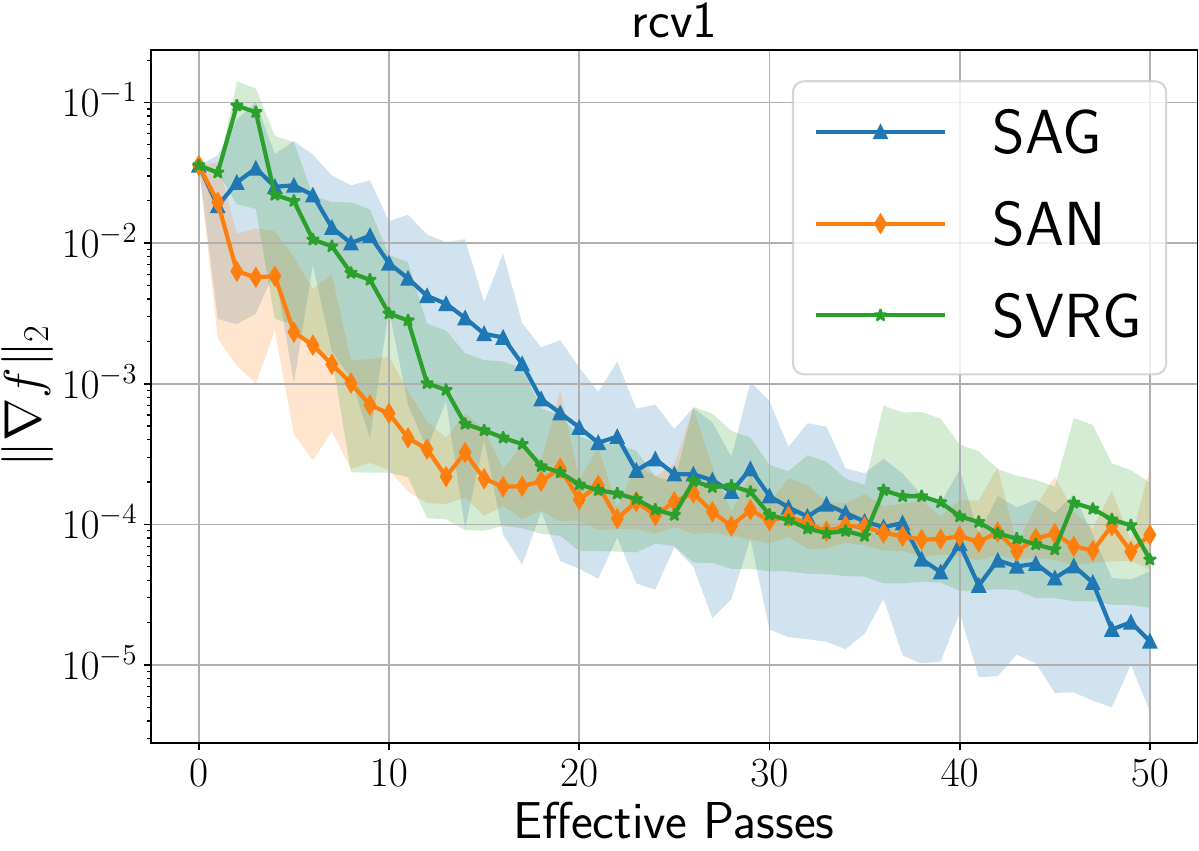}
\includegraphics[width=.24\linewidth]{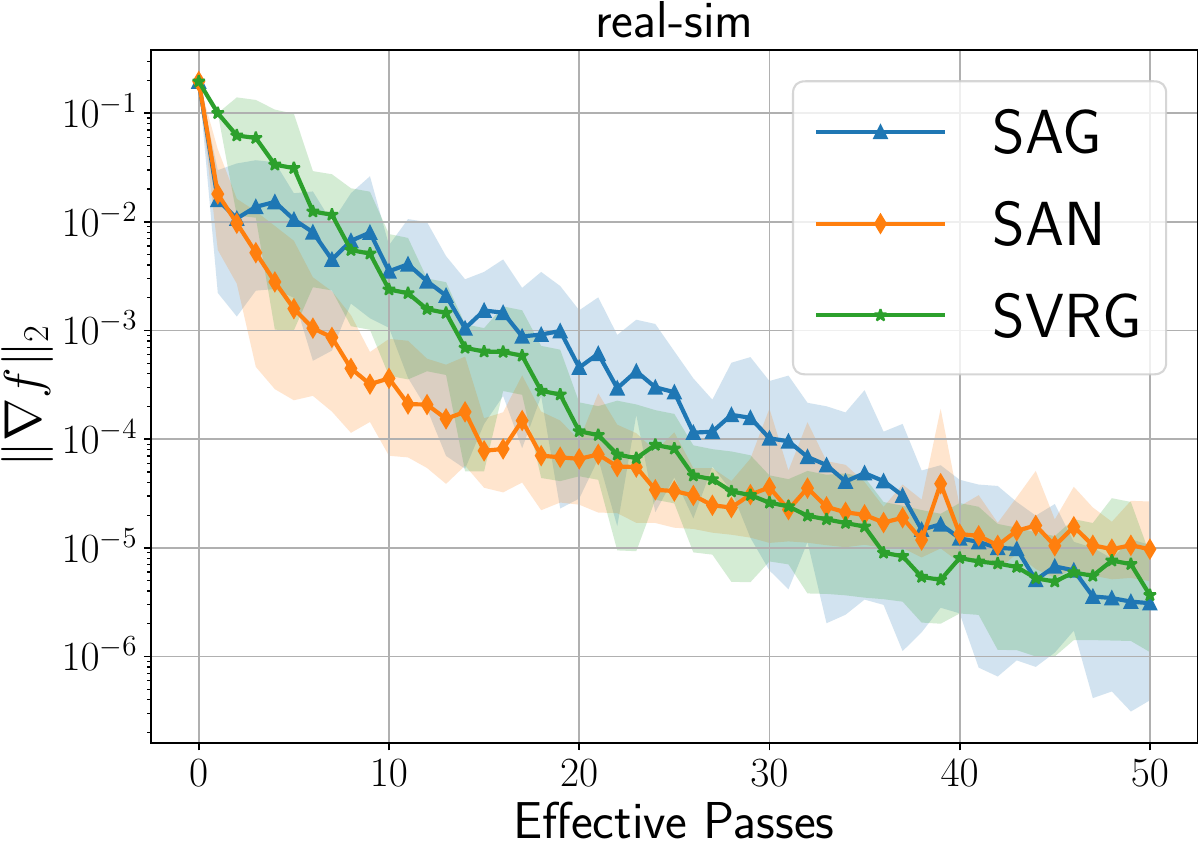}
\caption{Logistic regression with pseudo-Huber regularization.}
\label{fig:pseudo-HuberSAN}
\end{figure*}

\begin{table*}[t]
 \caption{Average cost of one iteration of various stochastic methods applied to GLMs. }
 \label{tab:complexity}
 \centering
 \begin{tabular}{lllll}
  \toprule
       & memory  & memory access & data access & computational cost  \\
  \midrule
  SAN     & $\cO(nd)$               & $\boldsymbol{\cO(d)}$ & $\boldsymbol{\cO(1)}$   & $\boldsymbol{\cO(d)}^{*}$ \\
  SANA    & $\cO(nd)$               & $\cO(nd)$             & $\boldsymbol{\cO(1)}$   & $\cO(nd)$             \\
  SAG     & $\cO(nd)$               & $\boldsymbol{\cO(d)}$ & $\boldsymbol{\cO(1)}$   & $\boldsymbol{\cO(d)}$ \\
  SVRG    & $\boldsymbol{\cO(d)}$   & $\boldsymbol{\cO(d)}$ & $\boldsymbol{\cO(1)}$   & $\boldsymbol{\cO(d)}$ \\
  SNM     & $\cO(n+d^2)$            & $\cO(d^2)$            & $\boldsymbol{\cO(1)}$   & $\cO(d^2)^{**}$        \\
  IQN     & $\cO(nd^2)$            & $\cO(d^2)$            & $\boldsymbol{\cO(1)}$   & $\cO(d^2)$        \\
  \bottomrule
 \end{tabular} \\
 {$^{*}$\footnotesize For SAN this $\cO(d)$ computational cost is derived when $\pi = \cO(1/n)$.} \\ 
 {$^{**}$\footnotesize For SNM this $\cO(d^2)$ computational cost only holds for a L2 regularizer.}
\end{table*}

\section{Experiments for SAN applied for GLMs} \label{sec:exp}

Here we compare  SAN in Algorithm~\ref{algo:SAN} against two variance reduced gradient methods SAG~\citep{SAG} and SVRG~\citep{Johnson2013} for
 solving regularized GLMs \eqref{eq:glm}, 
where  $\phi_i(t) = \log\left(1 + e^{-y_i t} \right)$ is the logistic loss, $y_i \in \{-1, 1\}$ is the $i$-th target value, and $R$ is the regularizer. 

We use eight datasets in our experiments taken from LibSVM~\citep{libsvm},\footnote{All datasets can be found downloaded on \url{https://www.csie.ntu.edu.tw/~cjlin/libsvmtools/datasets/}. Some of the datasets can be found originally in~\citep{phishing,ijcnn1,covtype,webspam,rcv1,UCI}.} with disparate properties (see details of the datasets in Table~\ref{tab:binary_datasets}). 
We fixed an initial random seed, evaluated each method $10$ times, and stopped when the gradient norm was below $10^{-6}$ or a maximum of $50$ effective passes over data had been reached.
In all of our experiments, we plot  effective data passes\footnote{By effective data passes we mean the number of data access divided by $n$.} vs gradient norm, and plot the central tendency as a solid line and all other executions as a shaded region. Plots with function sub-optimality are also provided in Figure~\ref{fig:sub} Section~\ref{sec:sub_opt}, and show much the same relative rankings amongst the methods as the gradient norm plots.

For all methods, we used the default step size. For instance,  for SAG and SVRG we use the step size $\frac{1}{L_{\max}}$ where $L_{\max}$ is the largest smoothness constant of $f_i$, for $i=1,\ldots, n.$ This step size is significantly larger than what has been proven to work for SAG and SVRG.\footnote{SAG has been proven to converge with a step size of $1/16L_{\max}$~\citep{SAG} and  SVRG provably converges with a step size of $1/10L_{\max}$ and loop size of $m=10 L_{\max}/\mu$ where $\mu$ is the strong convexity parameter $f(w)$~\citep{Johnson2013}.} Yet despite this, it is the default setting in sklearn's logistic regression solver~\citep{scikit-learn}, and we also found that it worked well in practice.
The other hyperparameter of SVRG is the inner loop size which is set to $n$ throughout our experiments.
As for SAN, we set the probability $\pi = \frac{1}{n+1}$ and step size $\gamma = 1$.
 More details of the experiments are in Section~\ref{sec:expdetails}.
 \footnote{The code is available on \url{https://github.com/nathansiae/Stochastic-Average-Newton}.}


\paragraph*{Logistic regression with L2 regularization.}
We consider L2-regularized logistic regression, i.e. the regularizer is $R(w) = \frac{\lambda}{2}\norm{w}^2$ with $\lambda =1/n.$ From Figure~\ref{fig:L2SAN}, SAN outperforms SAG and SVRG in all eight datasets, except for \texttt{mushrooms}, \texttt{ijcnn1} and \texttt{covtype}, where SAN remains competitive with SAG or SVRG. 
Note as well that for reaching an approximate solution at early stage, SAN outperforms SAG and SVRG in all datasets, except for \texttt{covtype}.
Furthermore, SAN often has a smaller variance compared to SAG and SVRG based on eye-balling the shaded error bars in Figure~\ref{fig:L2SAN} which was produced by multiple executions.
\vspace{-0.3cm}

\paragraph*{Logistic regression with pseudo-Huber regularization.}

We also tested logistic regression with pseudo-Huber regularizer. The pseudo-Huber regularizer is defined as
$R(w) =\lambda \sum_{i=1}^dR_i(w_i)$ with $R_i(w_i) =\delta^2 \left( \sqrt{1 + \left(\frac{w_i}{\delta} \right)^2} - 1 \right)$ and is used to promote the sparsity of the solution~\citep{FountoulakisG16}.
We set $\delta = 1$ and $\lambda=1/n$. See Section~\ref{sec:expdetails} for more properties and interpretations of the pseudo-Huber regularizer. From Figure~\ref{fig:pseudo-HuberSAN}, SAN is competitive with SAG and SVRG.  Changing the regularizer from an L2 to pseudo-Huber has resulted in a slower convergence for all methods, except on the datasets  \texttt{ijcnn1} and \texttt{covtype}.
SVRG is notably slower when using the pseudo-Huber regularizer, while SAG is the least affected, and SAN is in between.
Besides, SAN again 
outperforms SAG and SVRG in all datasets, except for \texttt{covtype}, for reaching an approximate solution at early stage and
has a smaller variance compared to SAG and SVRG.

Overall, these tests confirm that SAN is efficient for a wide variety of datasets and problems. SAN is efficient in terms of effective passes and cost per iteration. It benefits from both using second order information yet still has the 
same cost $\cO(d)$ as the stochastic first order methods.  SNM and IQN,  the only other methods that fit our stated objective,
have a $\cO(d^2)$ cost per iteration for L2-regularized GLMs. 
SNM costs even more for other regularizers and IQN has a $\cO(nd^2)$ memory cost. 
In Section~\ref{sec:exp_extra}, we present experiments comparing SAN/SANA to the SNM and IQN algorithms.

Another  advantage of SAN is that it requires no prior knowledge of the datasets nor tuning of the hyperparameters. To show this, we did a grid search over  $\pi$ and the step size $\gamma$ of SAN, see Tables~\ref{tab: gridsearchsancovtype}  and~\ref{tab: gridsearchsanijcnn1} in Section~\ref{sec:gridsearch}, 
where we found that SAN performs equally well for a wide range of values of $\pi$ and $\gamma$. Thus for simplicity we set $\pi = \frac{1}{n+1}$ and $\gamma = 1$ in the experiments.
In contrast, SAG and SVRG require the computation of $L_{\max}$ and the performance is highly affected by the step size (see Table~\ref{tab: gridsearchsvrg} and~\ref{tab: gridsearchsag} in Section~\ref{sec:gridsearch}). 

However, the downside of SAN is that it  stores $nd$ scalars much like  SAG/SAGA~\citep{SAGA_Nips}. 
See Table~\ref{tab:complexity} the comparison among different algorithms.



\section{Sketched Newton Raphson with a variable metric}
\label{sec:SNR}

\subsection{Presentation of the SNRVM algorithm}

Though our main focus is in solving the function splitting reformulation~\eqref{eq:alphaieqzero} and~\eqref{eq:alphainablafi}, we find that our forthcoming theory holds for a large class of variable metric \emph{Sketched Newton Raphson} methods~\citep{Yuan2020sketched}, of which SAN/SANA are special cases. All proofs are given in Section~\ref{sec:annex sec 4}. 

 In order to design such method, we first reformulated our original problem~\eqref{eq:finite_sum}  as a system of nonlinear equations $F(x) = 0$ for a given choice of a smooth map $F: \R^{p} \longrightarrow  \R^m$ and where $p,m\in \N$  are appropriately chosen dimensions, e.g., $p = m = (n+1)d$ in~\eqref{eq: def_F} for SAN/SANA. 
We then proposed using a subsampled Newton Raphson method for solving these nonlinear equations.
Here we extend this subsampling to make use of any randomized \emph{sketch} of the system. That is, consider a random \emph{sketching} matrix $\mS_k \in \R^{m \times \tau}$ sampled from some distribution, where $\tau \in \N$ is significantly smaller than $p$ or $m$. We use this sketching matrix to compress the rows of the  Newton system~\eqref{eq: newton system} at each iteration by left multiplying as follows
\begin{equation}\label{eq:linearFxsketched}
 \mS_k^\top F(x^k) + \mS_k^\top \nabla F(x^k)^\top (x-x^k) =0.
\end{equation}
The resulting system has $\tau$ rows and is under-determined. 
To pick a solution, we use the projection
\begin{align}\label{eq:linearFxsketchedproj}
x^{k+1} &= \arg\min_{x\in \R^p} \norm{x- x^k}^2 \\ 
&\mbox{ s.t. } \mS_k^\top F(x^k) + \mS_k^\top \nabla F(x^k)^\top (x-x^k) =0. \nonumber 
\end{align}
The method in~\eqref{eq:linearFxsketchedproj} is known as the SNR (\emph{Sketched Newton Raphson}) method~\citep{Yuan2020sketched}. The SNR method affords a lot of flexibility by choosing different distributions for the sketching matrices. Yet it is not flexible enough to include SAN/SANA, since these require projections under a variable metric. 
To allow for projections under norms other than the L2 norm, 
we introduce a random positive-definite \emph{metric} matrix $\mW_k \in \R^{p \times p}$ that defines the norm under which we project.
Introducing as well a damping parameter $\gamma >0$, as we did for SAN and SANA, we obtain the following method
\begin{align}\label{algo:SNRVM implicit}
\bar x^{k+1} &= \argmin \Vert x - x^k \Vert_{\mW_k}^2 \nonumber \\ 
&\text{ s.t. } \mS_k^\top \nabla  F(x^k)^\top (x-x^k) = - \mS_k^\top F(x^k), \nonumber \\
x^{k+1} &= (1 - \gamma) x^k + \gamma \bar x^{k+1}.
\end{align}
We call this method the \emph{Sketched Newton Raphson with Variable Metric} (SNRVM for short).
The closed form expression\footnote{See Lemma~\ref{L:SNRVM implicit to explicit} in the supp. material for a proof. } for the iterates~\eqref{algo:SNRVM implicit} is given by
\begin{align} \label{algo:SNRVM}
x^{k+1} &= x^k - \gamma\mW_k^{-1} \nabla F(x^k)\mS_k \\
&\quad \ \cdot \left(\mS_k^\top  \nabla F(x^k)^\top \mW_k^{-1} \nabla F(x^k)\mS_k \right)^\dagger \mS_k^\top F(x^k). \nonumber 
\end{align}
SAN/SANA are both instances of the SNRVM method by choosing $\mS_k$ as a subsampling matrix and $\mW_k$ depending on the stochastic Hessian matrices. In Section~\ref{sec:SNR viewpoint} we provide a detailed derivation of  
SAN/SANA as an instance of the SNRVM method.

We assume that at each iteration, the random matrices $(\mS_k, \mW_k)$ are sampled according to a proper finite distribution $\mathcal{D}_{x^k}$ defined in the following.
\begin{assumption}[Proper finite distribution]\label{Ass:distribution D_x}
For every  $x \in \mathbb{R}^{p}$,
there exists $r \in \N$,  probabilities $\pi_1, \ldots, \pi_r >0$ with $\sum_{i=1}^r \pi_i =1$,
and matrices $\{\mS_i(x), \mW_i(x)\}_{i=1}^r$ s.t. for $i=1, \ldots, r$, we have 
\[\mathbb{P}_{(\mS,\mW) \sim \mathcal{D}_x}\left[(\mS,\mW) = (\mS_i(x), \mW_i(x))\right] \; =\; \pi_i.  \]
\end{assumption}


%



The addition of a variable metric to SNR has proven to be very challenging in terms of establishing a convergence theory. 
The current convergence theory and proofs techniques for SNR in~\citet{Yuan2020sketched} all fail with the addition of a variable metric.
This is not so surprising, considering the historic difficulty in developing theory for variable metric methods such as the quasi-Newton methods. Despite the immense practical success of quasi-Newton methods, a meaningful non-asymptotic convergence rate has eluded the optimization community for 70 years, with the first results having only just appeared last year~\citep{greedyqNRodomanov2020,Rodomanov2021super,rodomanov2021rates,jin2021nonasymptotic}.

In the following sections, we provide a general linear convergence theory for SNRVM in Section~\ref{sec:SNRVM} and a more explicit linear convergence rate for SAN and SANA in Section~\ref{sec:SAN_SANA}.

\subsection{Linear convergence rates for SNRVM}
\label{sec:SNRVM}

We start by introducing a  technical assumption which guarantees
  that~\eqref{algo:SNRVM implicit} is well defined, which is always true for SAN and SANA.
 
\begin{assumption}\label{ass:existence} 
For every $x \in\R^{p}$, 
the matrices $\nabla F(x)^\top \nabla F(x)$ and $\mathbb{E}_{\mathcal{D}_x}[\mS \mS^\top] $ are invertible, and every matrix $\mW \sim \mathcal{D}_x$ is symmetric positive definite.
\end{assumption}

\begin{proposition}\label{P:Ass verified for SAN and SANA}
Assumptions \ref{Ass:distribution D_x} and  \ref{ass:existence} are verified for SAN and SANA, under Assumption \ref{Ass:strict convexity of problem}.
\end{proposition}


Let us now introduce the surrogate function 
\begin{align}\label{eq: def_hat_f}
\hat{f}_k(x) \eqdef \frac{1}{2}\norm{F(x)}_{( \nabla F(x^k)^\top \mW_k^{-1}  \nabla F(x^k))^\dagger}^2,
\end{align}
where $\mW_k \sim \mathcal{D}_{x^k}$.
This function is closely related to the SNRVM algorithm.
Indeed, it is possible to show that $x^{k+1}$ is obtained by minimizing a quadratic approximation of $\hat{f}_k$ along a random subspace (see Lemma~\ref{L:SNRVM as minimize quadratic over random subspace}). More precisely, $x^{k+1}$ is the solution of
\begin{align*}
\textstyle\underset{x \in \R^p}{\argmin} & \ \hat{f}_k(x^k) + \langle \nabla \hat{f}_k(x^k) , x-x^k \rangle + \frac{1}{2\gamma}\norm{x-x^k}^2_{\mW_k} \\
\mbox{ s.t. } & x \in x^k +\Image{ \mW_k^{-1} \nabla F(x^k) \mS_k}
\end{align*}
Our forthcoming Theorem \ref{them: AdapSNR} shows that $\hat f_k(x^k)$ converges linearly to zero in expectation.
To achieve this, we need to make an assumption  that controls the evolution of $\hat{f}_k$ along the iterations.
\begin{assumption}\label{ass:AdapUpperBnd}
There exists $L>0$ such that, for every $k \in \N$ and every $x\in \mathbb{R}^p$:
\resizebox{\hsize}{!}{%
        \begin{math}
\begin{aligned}
\textstyle\hat{f}_{k+1}(x)  \leq \hat{f}_k(x^k) + \langle\nabla \hat{f}_k(x^k) , x-x^k\rangle  + \frac{L}{2}\norm{x-x^k}^2_{\mW_k}.
\end{aligned}
        \end{math}%
    }%
\end{assumption}
We now state our core convergence result, which we prove in Section \ref{asec:proofTheoremAdapSNR} in the Appendix.
\begin{theorem}\label{them: AdapSNR}
Let Assumptions~\ref{Ass:distribution D_x}, \ref{ass:existence} and \ref{ass:AdapUpperBnd} hold, and let $\gamma = \left. 1 \right/L$.
Let $(\mS,\mW)\sim\cD_x$ and let
\begin{align*}
\mH(x) &\eqdef \mathbb{E}\left[\mS\left(\mS^\top \nabla F(x)^\top \mW^{-1} \nabla F(x)\mS\right)^\dagger\mS^\top\right], \\
\rho(x) &\eqdef \min\limits_{i=1, \dots, r}
\lambda^{+}_{\min} \left( \mM_i \mH(x) \mM_i^\top \right),
\end{align*}
where $\mM_i \eqdef \mW_i(x^k)^{-\frac{1}{2}}\nabla F(x)$. Assume that there exists $\rho > 0$ such that $\inf\limits_{k \in \mathbb{N}} \ \rho(x^k) \geq \rho$ almost surely.
 It follows that
   \begin{eqnarray*}\label{eq:RSNlinearconv} 
 \E{\hat{f}_{k}(x^{k})} &\leq&   \left(1  -  \rho\gamma \right)^k \E{\hat{f}_0(x^0)}  \text{ a.s.}
\end{eqnarray*}
\end{theorem}
When the metric is constant along iterations and $F(x)$ is a linear function, or equivalently our original problem~\eqref{eq:finite_sum} is a quadratic problem, then the SNRVM method~\eqref{algo:SNRVM} is known as the sketch-and-project method~\citep{Gower2015}. 
 In Section~\ref{asec:linearexample}, we show that Theorem~\ref{them: AdapSNR} when specialized to this case allows us to recover the well known convergence rates for solving linear systems using sketch-and-project. 

The convergence result of SNRVM in Theorem~\ref{them: AdapSNR} is for the surrogate function $\hat{f}_k(x^k)$ and the convergence rate $\rho$ is not explicit. Next, we develop a linear convergence theory of SNRVM for $\norm{F(x^k)}^2$ with the explicit linear convergence rate $\rho$.

Indeed, the existence of a lower bound $\rho>0$ in Theorem \ref{them: AdapSNR} can be guaranteed, provided that we can uniformly control the matrices $\mS,\mW$ and $\nabla F(x)$.
	Let us make this more precise:
	
	\begin{assumption}\label{Ass:explicit rates for SNRVM}
	Assumption \ref{ass:existence} holds, 
		$m=p$~and $F$ is injective.
	We assume that there exists a 
	set $\Omega \subset \mathbb{R}^p$ 
	and constants $\mu_W,L_W,\bar \mu_S, L_S,\mu_{\nabla F}, L_{\nabla F}$ 
	in $(0,+\infty)$ 
	such that, for all $x \in \Omega$, for all $(\mS,\mW) \sim\mathcal{D}_x$, 
		\begin{align*}
		&{\rm{spec}}(\mW) \subset [\mu_W,L_W],
		\quad 
		\sigma(\nabla F(x)) \subset [\mu_{\nabla F}, L_{\nabla F}], \\
		&\bar \mu_S \leq \lambda_{min}\left( \mathbb{E}\left[ \mS\mS^\top \right] \right),
		\quad
		\Vert \mS\mS^\top \Vert \leq L_S.
		\end{align*}
		where ${\rm{spec}}(M)$ (resp. $\sigma(M)$) denote the set of eigenvalues (resp. singular values) of a square matrix $M$.
	\end{assumption}

This assumption is typically verified on bounded sets, if the matrices $\mW, \mS, \nabla F(x)$ enjoy some sort of continuity with respect to $x$ (we detail this argument for SAN in Section \ref{SS:SAN bounded sequence} in the Appendix).
Now we can state our general linear convergence theory of SNRVM in terms of $F(x^k)$ itself, instead of the surrogate $\hat f_k(x^k)$, with explicit rates.

	\begin{theorem}\label{T:SNRVM explicit rates}
	Let Assumptions~\ref{Ass:distribution D_x}, \ref{ass:existence}, \ref{ass:AdapUpperBnd}, and \ref{Ass:explicit rates for SNRVM} hold, and let $\gamma = \left. 1 \right/L$.
	Let $\{x^k\}_{k \in \mathbb{N}}$ be generated by the SNRVM algorithm, and suppose that $x^k \in \Omega$ 
	almost surely.
	Then for all $k \in \mathbb{N}$ we have that
	\begin{equation*}
	\mathbb{E}\left[ \Vert F(x^k) \Vert^2 \right] \leq C(1-\gamma \rho)^k 
	\quad \text{ almost surely},
	\end{equation*}
	with
	$\rho = \frac{\mu_{\nabla F}^2}{L_{\nabla F}^{2}} 
	\frac{\mu_{W}}{L_{W}}
	\frac{\bar \mu_S}{L_S}$, and
	$C=2
	\E{\hat{f}_0(x^0)}
	\frac{L_{\nabla F}^{2}}{\mu_W}$.
	\end{theorem}

\subsection{Linear convergence rates for SAN and SANA}
\label{sec:SAN_SANA}

Theorem~\ref{T:SNRVM explicit rates} provides a general convergence theory for SNRVM. In this section, when specialized to SAN and SANA, we are able to get more insightful convergence results.

Indeed, Assumption \ref{Ass:explicit rates for SNRVM} can be ensured for SAN and SANA
under reasonable assumptions on the regularity of the functions $f_i$, in the following sense:

\begin{assumption}\label{Ass:explicit rates for SAN}
There exists $0<\mu_f \leq L_f$ such that for all $i\in \{1, \dots, n \}$, the function $f_i : \mathbb{R}^d \longrightarrow \mathbb{R}$ is of class $C^2$, $\mu_f$-strongly convex and has a $L_f$-Lipschitz continuous gradient.
\end{assumption}

The next Theorem, which is our main convergence result for SAN and SANA, describes linear rates for the iterates $(w^k, \alpha_1^k, \dots, \alpha_n^k)$ themselves:

\begin{theorem}\label{T:CV SAN explicit}
Let Assumptions \ref{ass:AdapUpperBnd} and \ref{Ass:explicit rates for SAN} hold.
Let $\{x^k\}_{k \in \mathbb{N}}$ be a sequence generated by SAN with $\pi = 1/(n+1)$, or by SANA, with $\gamma=1/L$.
Let $w^* = {\rm{argmin}}~f$.
Then for every $k \in \mathbb{N}$:
\begin{equation*}
\mathbb{E}\hspace*{-2pt}\left[ \Vert w^k\hspace*{-2pt}-\hspace*{-2pt}w^* \Vert^2  \right]
+ \sum\limits_{i=1}^n \mathbb{E}\hspace*{-2pt}\left[   \Vert \alpha_i^k\hspace*{-2pt}-\hspace*{-2pt}\nabla f_i(w^*) \Vert^2 \right]
\hspace*{-2pt}\leq\hspace*{-2pt}C (1-\gamma \rho)^k
\end{equation*}
holds almost surely, where we can take
\begin{equation*}
\scriptstyle \rho = \frac{\min\{1,\mu_f^3\}}{14n^3(2+L_f^2)^2\max\{1,L_f^3\}},
\
C=
18n^2
\frac{\max\{1,L_f^2\}(2+L_f^2)^2}{\min\{1,\mu_f^3\}}\E{\hat{f}_0(x^0)}.
\end{equation*}
\end{theorem}

The resulting linear rate of convergence $\rho$ in Theorem~\ref{T:CV SAN explicit} depends on $n$ and converges to $1$ as $n$ goes to infinity.
This is not surprising, since it is also the case for variance reduced methods  such as SAGA~\cite{SAGA_Nips} and SVRG~\cite{Johnson2013}.
We also note that our theoretical rate has a much worse dependence in $\mu_f$, $L_f$ and $n$ than the rates of SVRG and SAGA, because of the presence of exponents greater than $1$.
This might suggest that our analysis is not tight: indeed we observed empirically that SAN performs as well as SVRG and SAG,
 even in a regime where $n$ is large and the problem is severely ill-conditioned (see Table \ref{tab:binary_datasets} for more details).

%
%

\if{
The assumption on the existence of a lower bound $\rho >0$ fortunately holds for SAN and SANA so long as the iterates are bounded, which we observe in all our experiments.
%
\begin{proposition}\label{P:rhopositive SAN}
Consider the setting of Theorem \ref{them: AdapSNR}. If $\{x^k\}_{k \in \mathbb{N}}$ is generated by SAN or SANA  and is  bounded almost surely then there exists $\rho>0$ such that $\inf\limits_{k \in \mathbb{N}} \ \rho(x^k) \geq \rho$ almost surely.
\end{proposition}
}\fi

\section{Conclusion}
\label{sec:conclusion}
We introduced the use of a subsampled Newton Raphson method applied to a specific function splitting problem as a tool for designing new incremental Newton methods.
We showcase this by developing SAN, an average Newton method that is empirically highly competitive as compared to variance reduced gradient methods, and does not require parameter tuning.
Further venues of investigation include:
\begin{itemize}[parsep=0em, topsep=0em]
	\item Improving our theoretical analysis, to obtain rates that better matches the ones of usual variance reduced methods, motivated by our numerical experiments.
	\item Leveraging SNRVM's structure to design a more efficient variant of SAN including mini-batching.
	This should be simple thanks to our function splitting point of view, as we would simply sample many rows at once in \eqref{eq:weqai}.
	\item Exploring different sketching techniques together with original and alternative splitting point of views to design methods that have not been discovered yet. 
\end{itemize}



\bibliography{references}
\bibliographystyle{unsrtnat}

\onecolumn



\newpage

\appendix

\part{Appendix}

\parttoc




The Appendix is organized as follows: In Section~\ref{S:closed form SAN}, we carefully derive the closed form updates of SAN and SANA presented in Algorithm~\ref{algo:SAN} and~\ref{algo:SANA}. 
In Section~\ref{sec:GLMs}, we specialize SAN and SANA for the case of regularized generalized linear models and provide more detailed and efficient pseudo-codes for such case. 
In Sections~\ref{sec:exp2}, we give further details on the numerical experiments and provide additional experiments for SANA to compare with SNM and IQN. 
In Section~\ref{sec:SNR viewpoint} and~\ref{sec:annex sec 4}, we provide the proofs for the claims and results in Section~\ref{sec:SNR}.

\section{A closed form expression for SAN and SANA}\label{S:closed form SAN}

In this section,  we  show that the updates of the SAN method given in Algorithm \ref{algo:SAN} are equivalent to the implicit formulation in \eqref{eq:SAN implicit linearized 0}-\eqref{eq:SAN implicit linearized 1..n}.
We then derive the closed form updates of the  SANA method in Section~\ref{asec:closedformSANA}. Finally in Section~\ref{sec:generic projection}, we provide a useful lemma. It provides an alternatively way to directly deduce the closed form updates~\eqref{eq:SAN implicit linearized 0} and~\eqref{eq:SAN implicit linearized 1..n} of SAN.

\subsection{Closed form expression for SAN}
\label{sec:closed form SAN}

We start with the following technical lemma.
\begin{lemma}\label{L:SAN linearized implicit to explicit}
Let $j \in \{1, \dots, n\}$.
Let $\hat w \in \mathbb{R}^d$, and $\hat \alpha_1, \dots, \hat \alpha_n \in \mathbb{R}^d$.
Let $c_j \in  \mathbb{R}^d$, and $\mH_j \in  \R^{d\times d}$ be a symmetric positive definite matrix.
The optimization problem
\begin{eqnarray*}
\min_{w, \alpha_1, \dots, \alpha_n \in \R^d} & \frac{1}{2}\sum\limits_{i=1}^n \norm{\alpha_i - \hat \alpha_i}^2 + \frac{1}{2}\norm{w-\hat w}_{\mH_j}^2 \nonumber \\
 \;\;\mbox{subject to }  &  \mH_j (w-\hat w) - \alpha_j =  c_j, \nonumber 
\end{eqnarray*}
has a unique solution $(w, \alpha_1, \dots, \alpha_n)$ given by
\begin{eqnarray*}
w & = &  \hat w + (\mI_d + \mH_j)^{-1}  (c_j + \hat \alpha_j), \\
\alpha_j & = & \hat \alpha_j - (\mI_d + \mH_j)^{-1}  (c_j + \hat \alpha_j), \\
\alpha_i &=& \hat \alpha_i \text{ for } i \neq j.
\end{eqnarray*}
\end{lemma}

\begin{proof}
Denoting $x = (w,\alpha_1, \dots, \alpha_n)$, let us define 
\begin{equation}\label{e:SAN linearized implicit to explicit:functions}
\Phi(x) \eqdef \frac{1}{2}\sum\limits_{i=1}^n \norm{\alpha_i - \hat \alpha_i}^2 + \frac{1}{2}\norm{w-\hat w}_{{\mH_j}}^2,
\quad \text{ and } \quad 
\Psi_j(x) \eqdef \mH_j(w-\hat w) - \alpha_j -  c_j.
\end{equation}
The fact that ${\mH_j}$ is positive definite implies that $\Phi$ is strongly convex.
Moreover $\Psi_j$ is affine, so we deduce that this problem has a unique solution.
Moreover, this solution, let us call it $x=(w,\alpha_1, \dots, \alpha_n)$, is characterized as the unique vector in $\mathbb{R}^{(n+1)d}$ satisfying the following KKT conditions:
\begin{equation}\label{eq:KKTSAN}
(\exists \beta_j \in \mathbb{R}^{d}) \quad \text{ such that } \quad 
\begin{cases}
\nabla \Phi(x)  + \nabla \Psi_j(x)\beta_j = 0, \\
 \Psi_j(x) = 0.
\end{cases}
\end{equation}
The derivatives in the above KKT conditions are given by
\begin{equation}\label{e:SAN linearized implicit to explicit:gradients}
\nabla  \Phi(x) = 
\begin{bmatrix}
{\mH_j}(w - \hat w) \\
\alpha_1 - \hat \alpha_1 \\
\vdots \\
\alpha_n - \hat \alpha_n
\end{bmatrix} 
\qquad 
\mbox{and}
\qquad
\nabla  \Psi_j(x) =
\begin{array}{cc}
\begin{bmatrix}
\mH_j \\ \mo_d \\ \vdots \\ \mI_d \\ \vdots \\ \mo_d
\end{bmatrix}
&
\begin{array}{c}
 \\
\leftarrow j+1
\end{array}
\end{array}
\end{equation}
Using the expression for these derivatives, we can rewrite the KKT conditions~\eqref{eq:KKTSAN} as
\begin{equation} \label{eq:kktsan}
(\exists \beta_j \in \mathbb{R}^{d}) \quad \text{ such that } \quad 
\begin{cases}
\mH_j(w - \hat w) + \mH_j \beta_j = 0, \\
\alpha_j - \hat \alpha_j  - \beta_j = 0, \\
\alpha_i - \hat \alpha_i  +0 = 0, \text{ for all } i \neq j \\
\mH_j(w-\hat w) - \alpha_j - c_j = 0.
\end{cases} 
\end{equation}
We immediately see that $\alpha_i = \hat \alpha_i$ for $i \neq j$.
Combining the second and fourth equations in \eqref{eq:kktsan}, we obtain
\begin{equation*}
\beta_j = \alpha_j - \hat \alpha_j = \mH_j(w - \hat w) - c_j - \hat \alpha_j.
\end{equation*}
Multiplying this new equality by $\mH_j$ allows us to rewrite the first equation in \eqref{eq:kktsan} as:
\begin{equation*}
\mH_j(w - \hat w) +  \mH_j \beta_j = 0
\Leftrightarrow 
\mH_j(w - \hat w) +  \mH_j^2  (w- \hat w) =  \mH_j (c_j + \hat \alpha_j).
\end{equation*}
Using the fact that $\mH_j$ is invertible, the latter is equivalent to write:
\begin{equation*}
w = \hat w + \left( \mI_d +  \mH_j \right)^{-1}  (c_j + \hat \alpha_j).
\end{equation*}
Moreover, since $\mH_j (\mI_d + \mH_j)^{-1} = \mI_d - (\mI_d + \mH_j)^{-1}$, we can also turn the fourth equation in \eqref{eq:kktsan} into
\begin{equation*}
\alpha_j 
= \mH_j(w- \hat w) - c_j
= \left( \mI_d - (\mI_d + \mH_j)^{-1}  \right) (c_j + \hat \alpha_j) - c_j 
=
\hat \alpha_j - (\mI_d + \mH_j)^{-1}  (c_j + \hat \alpha_j).
\end{equation*}
This proves the claim.
\end{proof}

\begin{lemma}\label{L:SAN implicit to explicit}
Let $\pi \in [0,\,1 ]$ and $\gamma \in (0,\,1]$ a step size.
Algorithm \ref{algo:SAN} (SAN) is equivalent to the following algorithm:  \\
With probability $\pi$, update according to
\begin{equation}\label{e:sanie1}
\begin{cases}
\bar x^{k+1} = {\rm{argmin}}~ \Vert w - w^k \Vert^2 + \sum\limits_{i=1}^n \Vert \alpha_i - \alpha_i^k \Vert^2 \quad  \text{ subject to } \quad 
\frac{1}{n}\sum\limits_{i=1}^n \alpha_i = 0, \\
x^{k+1} = (1 - \gamma) x^k + \gamma \bar x^{k+1},
\end{cases}
\end{equation}
Otherwise with probability $(1-\pi)$, sample $j \sim \{1, \dots, n\}$ uniformly and update according to
\begin{equation}\label{e:sanie2}
\begin{cases}
\bar x^{k+1} = & {\rm{argmin}}~ \Vert w - w^k \Vert_{\nabla^2 f_j(w^k)}^2 + \sum\limits_{i=1}^n \Vert \alpha_i - \alpha_i^k \Vert^2 \\
 &  \text{ subject to } \quad 
\nabla^2 f_j(w^k)(w - w^k) - \alpha_j = - \nabla f_j(w^k), \\
x^{k+1} = & (1 - \gamma) x^k + \gamma \bar x^{k+1}.
\end{cases}
\end{equation}
\end{lemma}

\begin{proof}
Suppose that we are in the case (which holds with probability $\pi$)  given by~\eqref{e:sanie1}.
In the projection step, we see that $w$ is not present in the constraint, which implies that $\bar w^{k+1} = w^k$, and therefore $w^{k+1} = w^k$.
On the other hand, $(\bar \alpha_1, \dots, \bar \alpha_n)$ is the projection of $( \alpha_1, \dots,  \alpha_n)$ onto a simple linear constraint, and can be computed in closed form as
\begin{equation*}
(\forall i \in \{1, \dots, n\}) \quad 
\bar \alpha_i^{k+1} = \alpha_i^k - \frac{1}{n} \sum\limits_{i=1}^n \alpha_i^k.
\end{equation*}
Consequently
\begin{equation*}
(\forall i \in \{1, \dots, n\}) \quad 
\alpha_i^{k+1} = (1- \gamma) \alpha_i^k + \gamma \left(  \alpha_i^k - \frac{1}{n} \sum\limits_{i=1}^n \alpha_i^k \right)
=
\alpha_i^k - \frac{\gamma}{n} \sum\limits_{i=1}^n \alpha_i^k,
\end{equation*}
which gives us exactly the step \ref{ln:rsn_mean} in Algorithm \ref{algo:SAN}.

Let now $j$ be in $\{1, \dots, n\}$ sampled uniformly, and suppose that we are in the case given by~\eqref{e:sanie2}.
Using Lemma \ref{L:SAN linearized implicit to explicit} we can compute an explicit form for $\bar x^{k+1}$ given by
\begin{eqnarray*}
\bar w^{k+1} & = & w^k + \left( \mI_d + \nabla^2 f_j(w^k)  \right)^{-1}  (\alpha_j^k - \nabla f_j(w^k)), \\
\bar \alpha_j^{k+1} & = & \alpha_j^k -  \left( \mI_d + \nabla^2 f_j(w^k)  \right)^{-1}  (\alpha_j^k - \nabla f_j(w^k)), \\
\bar \alpha_i^{k+1} & = & \hat \alpha_i^{k} \text{ for all } i \neq j.
\end{eqnarray*}
Consequently, after applying the relaxation step we have
\begin{eqnarray*}
w^{k+1} & = & w^k + \gamma \left( \mI_d + \nabla^2 f_j(w^k)  \right)^{-1}  (\alpha_j^k - \nabla f_j(w^k)), \\
\alpha_j^{k+1} & = & \alpha_j^k -  \gamma \left( \mI_d + \nabla^2 f_j(w^k)  \right)^{-1}  (\alpha_j^k - \nabla f_j(w^k)), \\
\bar \alpha_i^{k+1} & = & \hat \alpha_i^{k} \text{ for all } i \neq j.
\end{eqnarray*}
which is exactly the steps \ref{ln:rsn_d}-\ref{ln:rsn_alpha} in Algorithm \ref{algo:SAN}.
\end{proof}

\subsection{Closed form expression for SANA}
\label{asec:closedformSANA}

\begin{lemma}\label{L:rsp linear implicit to explicit}
Let $j \in \{1, \dots, n\}$.
Let $c_j \in \mathbb{R}^d$, $\hat w \in \mathbb{R}^d$, and let $\hat \alpha_1, \dots, \hat \alpha_n \in \mathbb{R}^d$ be such that $\sum_{i=1}^n \hat \alpha_i=0$.
Let $\mH_j \in \R^{d\times d}$ be a  positive definite matrix.
Then the optimization problem
\begin{eqnarray}
\min_{w, \alpha_1, \dots, \alpha_n \in \R^d}  \frac{1}{2}\sum\limits_{i=1}^n \norm{\alpha_i - \hat \alpha_i}^2 + \frac{1}{2}\norm{w-\hat w}_{\mH_j}^2, \nonumber \\
 \;\;\mbox{subject to }  \mH_j (w-\hat w) - \alpha_j =  c_j, \nonumber  \\
 \phantom{\mbox{subject to }}  \sum\limits_{i=1}^n \alpha_i  = 0,\label{eq:rowsimrpobapp}
\end{eqnarray}
has a unique solution $(w, \alpha_1, \dots, \alpha_n)$ given by
\begin{eqnarray*}
d & =  &\left(  \frac{n-1}{n}  \mI_d + \mH_j \right)^{-1}  (c_j + \hat \alpha_j),\\
w &=&\hat w + d, \\
\alpha_j &=& \hat \alpha_j - \frac{n-1}{n} d, \\
\alpha_i &=& \hat \alpha_i + \frac{1}{n} d, \quad \text{ for } i \neq j.
\end{eqnarray*}
\end{lemma}

\begin{proof}
Noting $x = (w,\alpha_1, \dots, \alpha_n)$, let us define $\Phi(x)$ and $\Psi_j(x)$ as in \eqref{e:SAN linearized implicit to explicit:functions}, together with
\begin{equation*}
\Psi_0(x) \eqdef \sum\limits_{i=1}^n \alpha_i.
\end{equation*}
The fact that $\mH_j$ is positive definite implies that we are minimizing a strongly convex function over a set of affine equations.
We deduce that this problem has a unique solution.
Moreover, this solution, let us call it $x=(w,\alpha_1, \dots, \alpha_n)$, is characterized as the unique vector in $\mathbb{R}^{(n+1)d}$ satisfying the following KKT conditions:
\begin{equation*}
(\exists \beta_0, \beta_j \in \mathbb{R}^{d}) \quad \text{ such that } \quad 
\begin{cases}
\nabla \Phi(x) + \nabla  \Psi_0(x)\beta_0 + \nabla \Psi_j(x)\beta_j = 0, \\
 \Psi_0(x) = 0, \\
 \Psi_j(x) = 0.
\end{cases}
\end{equation*}
Here, we can compute $\nabla  \Phi(x)$ and $\nabla  \Psi_j(x)$ as in \eqref{e:SAN linearized implicit to explicit:gradients}, together with
\begin{equation*}
\nabla  \Psi_0(x) = 
\begin{bmatrix}
\mo_d \\
\mI_d \\
\vdots \\
\mI_d
\end{bmatrix}.
\end{equation*}
Therefore, we can rewrite the KKT conditions as
\begin{equation} \label{eq:kktrow}
(\exists \beta_0, \beta_j \in \mathbb{R}^{d}) 
\quad \text{ such that } \quad 
\begin{cases}
\mH_j(w - \hat w) + 0 +  \mH_j \beta_j = 0, \\
\alpha_j - \hat \alpha_j + \beta_0 - \beta_j = 0, \\
\alpha_i - \hat \alpha_i + \beta_0 +0 = 0, \text{ for all } i \neq j\\
\sum_{i=1}^n \alpha_i = 0 \\
\mH_j(w-\hat w) - \alpha_j - c_j = 0.
\end{cases} 
\end{equation}
The last equation in \eqref{eq:kktrow} can be rewritten as
\begin{equation}\label{eq:tempalphaoa9j3o8aj3}
\alpha_j = \mH_j(w-\hat w) - c_j.
\end{equation}
Summing the equations involving $\alpha_i$ for $i \neq j$, and using the fact that $\sum\limits_{i=1}^n \alpha_i = \sum\limits_{i=1}^n \hat \alpha_i =0$, together with \eqref{eq:tempalphaoa9j3o8aj3}, we can deduce that
\begin{equation*}
0 
= 
\sum\limits_{i \neq j} (\alpha_i - \hat \alpha_i + \beta_0 )
=
- \alpha_j +  \hat \alpha_j + (n-1) \beta_0
=
-  \mH_j (w-\hat w) + c_j + \hat \alpha_j +(n-1) \beta_0.
\end{equation*}
In other words, we obtain that
\begin{equation}\label{eq:tempbeta1smaoej9}
\beta_0 = \frac{1}{n-1}\left( \mH_j(w-\hat w) -  (c_j + \hat \alpha_j) \right).
\end{equation}
Injecting the above expression into the second equation of~\eqref{eq:kktrow}, and using again \eqref{eq:tempalphaoa9j3o8aj3}, gives
\begin{eqnarray*}\label{eq:beta2prev}
\beta_j 
& = &
\mH_j(w-\hat w) - c_j - \hat \alpha_j + \frac{1}{n-1}\left( \mH_j(w-\hat w) -  (c_j + \hat \alpha_j) \right) \nonumber \\
& = & 
 \frac{n}{n-1} \left(   \mH_j (w - \hat w) - (c_j + \hat \alpha_j)\right) 
\end{eqnarray*}
Combining this expression of $\beta_j$ with the first equation in \eqref{eq:kktrow} leads to
\begin{eqnarray*}
0 & = &
\mH_j(w-\hat w) + \frac{n}{n-1} \mH_j   \left(   \mH_j (w - \hat w) - (c_j + \hat \alpha_j)\right) \\
& = &
\mH_j  \left( 
\left( \mI_d + \frac{n}{n-1} \mH_j \right)
(w-\hat w) 
-\frac{n}{n-1}    (c_j + \hat \alpha_j)
\right)
\end{eqnarray*}
Using the fact that $\mH_j$ is positive definite, we obtain that:
\begin{eqnarray*}
w & = & \hat w + \frac{n}{n-1} \left(   \mI_d + \frac{n}{n-1} \mH_j \right)^{-1}  (c_j + \hat \alpha_j) \\
 & = & \hat w + \left(  \frac{n-1}{n}  \mI_d + \mH_j \right)^{-1}  (c_j + \hat \alpha_j), \\
 & = & \hat w + d,
\end{eqnarray*}
where $d$ is defined as $\left(  \frac{n-1}{n}  \mI_d + \mH_j \right)^{-1}  (c_j + \hat \alpha_j)$.

Going back now to \eqref{eq:tempalphaoa9j3o8aj3} we can write
\begin{eqnarray*}
\alpha_j &=& \mH_j d - c_j
=
\mH_j \left(  \frac{n-1}{n}  \mI_d + \mH_j \right)^{-1}  (c_j + \hat \alpha_j)- c_j \\
& = & 
\left( \mI_d - \frac{n-1}{n}\left(  \frac{n-1}{n}  \mI_d + \mH_j \right)^{-1} \right)(c_j + \hat \alpha_j)- c_j  \\
& = &
\hat \alpha_j - \frac{n-1}{n}\left(  \frac{n-1}{n}  \mI_d + \mH_j \right)^{-1}(c_j + \hat \alpha_j) \\
& = &
\hat \alpha_j - \frac{n-1}{n} d.
\end{eqnarray*}
It remains to compute $\alpha_i$, for $i \neq j$.
Start with the first equation of \eqref{eq:kktrow} and see that $w - \hat w + \beta_j=0$.
This implies that $\beta_j = -d$.
We can therefore use the second equation of \eqref{eq:kktrow} to write that 
\begin{equation*}
\beta_0 = \beta_j - (\alpha_j - \hat \alpha_j) = -d+ \frac{n-1}{n}d = -\frac{1}{n}d.
\end{equation*}
We can finally call the third equation of \eqref{eq:kktrow} and write that $\alpha_i = \hat \alpha_i - \beta_0 = \hat \alpha_i +\frac{1}{n}d$. 
\end{proof}

\begin{lemma}\label{L:SANA implicit to explicit}
Let $\gamma \in (0,\,1]$ be a step size.
Algorithm \ref{algo:SANA} (SANA) is equivalent to the following algorithm: 
update the iterates according to
\begin{equation}\label{e:sanaie}
\begin{cases}
\bar x^{k+1} = & {\rm{argmin}}~ \Vert w - w^k \Vert_{\nabla^2 f_j(w^k)}^2 + \sum\limits_{i=1}^n \Vert \alpha_i - \alpha_i^k \Vert^2 \\
 &  \text{ subject to } 
 \begin{cases}
 \frac{1}{n}\sum\limits_{i=1}^n \alpha_i = 0, \\
 \nabla^2 f_j(w^k)(w - w^k) - \alpha_j = - \nabla f_j(w^k),
 \end{cases}
 \\
x^{k+1} = & (1 - \gamma) x^k + \gamma \bar x^{k+1}.
\end{cases}
\end{equation}
\end{lemma}

\begin{proof}
Consider the iterates defined by \eqref{e:sanaie}.
Using Lemma \ref{L:rsp linear implicit to explicit}, we can compute an explicit form for $\bar x^{k+1}$:
\begin{eqnarray*}
d^k & =  &\left(  \frac{n-1}{n}  \mI_d + \nabla^2 f_j(w^k) \right)^{-1}  (\alpha_j^k - \nabla f_j(w^k)),\\
\bar w^{k+1} &=&w^k + d^k, \\
\bar \alpha_j^{k+1} &=& \alpha_j^k - \frac{n-1}{n} d^k, \\
\bar \alpha_i^{k+1} &=& \alpha_i^k + \frac{1}{n} d^k, \quad \text{ for } i \neq j.
\end{eqnarray*}
After applying the relaxation step $x^{k+1} =  (1 - \gamma) x^k + \gamma \bar x^{k+1}$, we obtain exactly the steps \ref{ln:SANA1}-\ref{ln:SANA4} in Algorithm~\ref{algo:SANA}.
\end{proof}

\subsection{Generic projection onto linear systems}
\label{sec:generic projection}

Here we provide a useful lemma that can directly deduce the closed form updates of~\eqref{eq:SAN implicit linearized 0} and~\eqref{eq:SAN implicit linearized 1..n} of SAN. It will also be used later in the appendix.

\begin{lemma} \label{lem:Hproject}
Let $\mA \in \R^{n\times d}$, $\mS \in \R^{n \times \tau}$, $b\in \Image{\mA}$, and $\mH$ be a symmetric positive definite matrix. The optimization problem
\begin{eqnarray*}
x^* &=& \arg \min_{x\in\R^d} \quad \frac{1}{2}\norm{x}_\mH^2, \\ 
&\quad& \mbox{subject to} \ \mS^\top \mA x = \mS^\top b,
\end{eqnarray*}
has a unique solution, called the weighted sketch-and-project optimal solution:
\begin{equation} \label{eq:Hproject}
x^* = \mH^{-1}\mA^\top \mS(\mS^\top \mA\mH^{-1}\mA^\top \mS)^\dagger \mS^\top b.
\end{equation}
\end{lemma}

\begin{proof}
First, note that this problem is strongly convex because $\mH$ is supposed to be positive definite, and therefore admits a unique solution $x^* \in \mathbb{R}^d$.
Second, since $\mH$ is invertible, we can do the change of variables $y \eqdef \mH^{1/2}x$. 
This allows us to write that  $x^* = \left(\mH\right)^{-1/2} y^*$ where $y^*$ is the unique solution of
\begin{align*}
\arg \min_{y\in\R^d}& \frac{1}{2}\norm{y}^2, \\ 
\mbox{subject to}& \ \mS^\top \mA \mH^{-1/2}y = \mS^\top b.
\end{align*}
The unique solution to the above problem is the minimal-norm solution of the linear system $\mS^\top \mA\left(\mH\right)^{-1/2}y = \mS^\top b$, which can be simply expressed by using the pseudo-inverse \cite[Definition 1]{BenCha63} :
\begin{align*}
y^* = \left( \mS^\top \mA \mH^{-1/2} \right)^\dagger \mS^\top b.
\end{align*}
Using the relation $\mM^\dagger = \mM^\top(\mM\mM^\top)^\dagger$ \cite[Lemma 1 \& Eq. 10]{Pen55}, we obtain
\begin{align*}
y^* = \left(\mH\right)^{-1/2}\mA^\top \mS(\mS^\top \mA\mH^{-1}\mA^\top \mS)^\dagger \mS^\top b.
\end{align*}
Multiplying this equality by $\mH^{-1/2}$ gives us the desired expression for $x^*$.
\end{proof}

This lemma is useful. Later it will be applied in Lemma~\ref{L:SNRVM implicit to explicit} and consequently provide the explicit updates of~\eqref{eq:SAN implicit linearized 0} and~\eqref{eq:SAN implicit linearized 1..n} of SAN in Section~\ref{sec:SNR viewpoint:SAN}. Thus this is a different way to obtain the closed form updates of SAN, compared to Section~\ref{sec:closed form SAN}.

\section{Implementations for regularized GLMs}\label{sec:GLMs}

\subsection{Definition and examples}

Here we specify our algorithms for the case of regularized generalized linear models.
Throughout this section, we assume that our finite sum minimization problem~\eqref{eq:finite_sum} 
is a  GLM (generalized linear model) defined as follows.
\begin{assumption}[Regularized GLM]\label{Ass:GLM regularized}
Our problem~\eqref{eq:finite_sum} writes as
\begin{equation}
\min_{w \in \R^d} \frac{1}{n} \sum_{i=1} f_i(w) \; \eqdef \;  \phi_i(\langle a_i,w \rangle)+R(w),
\end{equation}
where $\{a_i\}_{i=1}^n \subset \mathbb{R}^d$ are  \emph{data points}, $\{\phi_i\}_{i=1}^n$ are twice differentiable real convex \emph{loss} functions with $\phi_i''(t) >0$, and $R$ is a separable \emph{regularizer} with $R(w) = \sum\limits_{j=1}^d R_j(w_j)$ where $R_j$ is a twice differentiable real convex function with $R_j''(t) >0, $ for all $t \in \R.$
\end{assumption}


%
Some classic examples of GLMs include ridge regression where $\phi_i(t) = \frac{1}{2}(t - y_i)^2$ and $R_j(t) = \frac{\lambda}{2}t^2$ where $\lambda >0$ is a regularization parameter. L2-regularized logistic regression,   the example on which we perform most of our experiments, is also a GLM with
\begin{equation}
\phi_i(t) = \ln \left( 1 + e^{-y_i t} \right) \quad \mbox{and} \quad R_j(t) = \frac{\lambda}{2}t^2.
\end{equation}
 We also consider other forms of separable regularizers such as the pseudo-huber regularizer where $R_j(t) = \lambda \delta^2 \left(\sqrt{1+ \left(\frac{t}{\delta}\right)^2} -1 \right)$ where $\delta$ is a parameter.
%
%
%
%

In the next section, we will show that for GLMs, our methods can be efficiently implemented. But first we need the following preliminary results.
\begin{lemma}[Simple computations with Regularized GLMs]\label{L:GLM hessians}
For GLMs (Assumption \ref{Ass:GLM regularized}) 
we have for all $j \in \{1, \dots, n \}$, all $w \in \mathbb{R}^d$ and every $\mu \geq 0$ that
\begin{enumerate}[label=(\arabic*)]
  \item\label{itm:diagonal Hessian} $ \begin{array}{lcl}
   \nabla R(w) &=& [R_1'(w_1) \dots R_d'(w_d)]^\top,\\
   \nabla^2 R(w) &=& \Diag{R_1''(w_1), \dots, R_d''(w_d)}.
   \end{array}$
  \item\label{itm:f Hessian} $\begin{array}{lcl}
   \nabla f_j(w) &=& \nabla R(w) + \phi_j'(\langle a_j, w \rangle) a_j, \\
   \nabla^2 f_j(w) &=& \nabla^2 R(w) + \phi_j''(\langle a_j,w \rangle) a_j a_j^\top .
   \end{array}$
  \item\label{itm:inverse} With $\hat a_j := \left( \mu \mI_d + \nabla^2 R(w_k) \right)^{-1} a_j$, we have
   \begin{equation*}
   \left(  \mu \mI_d + \nabla^2 f_j(w) \right)^{-1} =
   \left( \mu \mI_d + \nabla^2 R(w_k) \right)^{-1} - \frac{\phi_j''(\langle a_j,w \rangle)}{1 + \phi_j''(\langle a_j,w \rangle) \langle \hat a_j, a_j \rangle} \hat a_j \hat a_j^\top.
   \end{equation*}  
  \item\label{itm:L2} If $R(w) = \frac{\lambda}{2}\Vert w \Vert^2$, with $\lambda >0$, then
   \begin{equation*}
   \left(  \mu \mI_d + \nabla^2 f_j(w) \right)^{-1} =
   \frac{1}{\mu +\lambda}\left(\mI_d - \frac{\phi_j''(\langle a_j,w \rangle)}{\mu  + \lambda + \phi_j''(\langle a_j,w \rangle) \Vert a_j \Vert^2}a_j a_j^\top\right).
   \end{equation*} 
\end{enumerate}
\end{lemma}


\begin{proof}~
\ref{itm:diagonal Hessian} and~\ref{itm:f Hessian} are trivial. For~\ref{itm:inverse}, let $\Phi := \phi_j''(\langle a_j,w \rangle)$, which is nonnegative because of the Assumption \ref{Ass:GLM regularized}.
Consider now the Sherman–Morrison formula:
\begin{equation*}
(\mM + uu^\top)^{-1} = \mM^{-1} - \frac{1}{1 + \langle \mM^{-1}u,u \rangle } (\mM^{-1} u)(\mM^{-1} u)^\top.
\end{equation*} 
This allows us to write, for $\mD = (\mu \mI_d +  \nabla^2 R(w))^{-1}$ and $\mM = \Phi^{-1}(\mu \mI_d +  \nabla^2 R(w))$, that 
  \begin{eqnarray*}
  \left(  \mu \mI_d + \nabla^2 f_j(w) \right)^{-1} 
  & =& 
  \left( \mu \mI_d + \nabla^2 R(w) + \Phi a_j a_j^\top \right)^{-1} 
  = 
  \Phi^{-1} \left( \mM + a_j a_j^\top \right)^{-1} \\  
  & =& 
  \Phi^{-1} \left(  \mM^{-1} - \frac{1}{1 + \langle \mM^{-1}a_j,a_j \rangle } (\mM^{-1} a_j)(\mM^{-1} a_j)^\top\right) \\
  & = &
  \mD - \frac{\Phi}{1 + \Phi\langle \mD a_j, a_j \rangle} (\mD a_j)(\mD a_j)^\top.
  \end{eqnarray*} 
\ref{itm:L2} is a direct consequence of the fact that
$\mD = \left( \mu \mI_d + \nabla^2 R(w_k) \right)^{-1} = \frac{1}{\mu + \lambda} \mI_d$.
\end{proof}

\subsection{SAN with GLMs}

Here we give the detailed derivation of our implementation of SAN  for GLMs, see Algorithm~\ref{algo:SAN_GLM}.
Upon examination, we can see that every step of Algorithm~\ref{algo:SAN_GLM} has a cost of $\cO(d)$, except on line~\ref{ln:SANGLMaverage}. As explained in Section~\ref{sec:SAN}, the averaging cost on line~\ref{ln:SANGLMaverage} costs $\cO(d)$ in which $\pi$ is of the order of $\cO(1/n)$. The only step that we have left an implicit computation is  on lines~\ref{ln:SANGLMRinver1} and~\ref{ln:SANGLMRinver2} which require inverting $(\mI_d+ \nabla^2 R(w^k))$. But this to comes at a cost of $\cO(d)$ since in our Assumption~\ref{Ass:GLM regularized} the regularizer is separable, and thus the Hessian is a diagonal matrix whose inversion also costs $\cO(d).$

\begin{algorithm}[H]
\caption{{SAN for regularized GLMs} \label{algo:SAN_GLM}}
\KwIn{%
Data $\{a_i \}_{i=1}^n$, loss functions $\{\phi_i\}_{i=1}^n$, regularizer $R$, 
$\pi \in (0,\,1)$, step size $\gamma\in (0,\,1]$, max iteration $T$}
Initialize $\alpha^0_1, \cdots, \alpha^0_n, w^0 \in \R^d$  and $\overline{\alpha}^0 = \frac{1}{n}  \sum_{i=1}^n \alpha^0_i.$

\For{$k=0, \dots, T-1$}{
    \With{
    $\alpha_i^{k+1} = \alpha_i^k - \gamma \overline{\alpha}^k , \quad \mbox{for all } i \in \{1, \cdots, n\}$ \label{ln:SANGLMaverage}

    $w^{k+1} = w^k$
    }
    \Otherwise{
    Sample $j \in \{1,\ldots, n\}$ uniformly

    $g^k =  \nabla R(w^k) + \phi^{\prime}_j( \langle a_j, w^k \rangle) a_j  - \alpha^k_j$

    $\hat{a}_j = (\mI_d + \nabla^2 R(w^k) )^{-1} a_j$ \label{ln:SANGLMRinver1}

    $d^k = \frac{\phi_j''(\langle a_j,w^k \rangle)\langle \hat a_j, g^k \rangle }{1 + \phi_j''(\langle a_j,w^k \rangle) \langle \hat a_j, a_j \rangle}\hat a_j
    - \left( \mI_d + \nabla^2 R(w^k) \right)^{-1}g^k$ \label{ln:SANGLMRinver2}
    
    $w^{k+1} = w^k + \gamma d^k$
    
    $\alpha^{k+1}_{j} = \alpha^k_j - \gamma d^k$
    
    $\alpha^{k+1}_{i} = \alpha^k_i $ for $i \neq j$
    
    $\overline{\alpha}^{k+1} =  \overline{\alpha}^{k} -  \frac{\gamma}{n}d^k$
    }
}
\KwOut{Last iterate $w^{T}$}
\end{algorithm}

Next we formalize the costs of Algorithm~\ref{algo:SAN_GLM} in the following remark. 
By  computational cost, we refer to the total number of floating point operations, that is  the number of scalar multiplications and additions.
\begin{remark}
The average costs of SAN (Algorithm~\ref{algo:SAN_GLM}) per iteration under Assumption~\ref{Ass:GLM regularized} are:
\begin{itemize}[nosep]
  \item \texttt{ Memory storage} of $\cO(nd)$ scalars.
  \item  \texttt{  Memory access} of $\cO(\pi nd + (1-\pi) d)$ which is $\cO(d)$ when $\pi \simeq 1/n$.
  \item  \texttt{ Data access} of $\cO(1)$ .
  \item \texttt{  Computational cost} of $\cO(\pi dn + (1-\pi)d)$ which is $\cO(d)$ when $\pi \simeq 1/n$.
\end{itemize}
\end{remark}
In calculating the average computational cost per iteration, we used that in expectation the updates on lines~ 4--5 occur with probability $\pi$, while the updates on lines~7--14 occur with probability $(1-\pi)$.

\begin{lemma}\label{L:SAN for GLMs}
The SAN Algorithm \ref{algo:SAN} applied to Regularized GLMs (in the sense of Assumption \ref{Ass:GLM regularized}) is Algorithm \ref{algo:SAN_GLM}.
\end{lemma}

\begin{proof}
Let $k \in \{0, \dots, T-1\}$. 
With probability $\pi$  from Algorithm \ref{algo:SAN} we have
\begin{equation*}
\alpha_i^{k+1} \;=\; \alpha_i^k - \frac{\gamma}{n}\sum_{j=1}^n\alpha_j^k, \quad \mbox{for all } i \in \{1, \cdots, n\}.
\end{equation*}
This can be rewritten as $\alpha_i^{k+1} = \alpha_i^k - \gamma \bar \alpha^k$, which is the update on  line 4 in Algorithm \ref{algo:SAN_GLM}.

With probability $(1-\pi)$  from Algorithm \ref{algo:SAN} we have
\begin{eqnarray*}
d^k &= &  - \left(\mI_d + \nabla^2f_j(w^k)\right)^{-1}\left( \nabla f_j(w^k) - \alpha_j^k\right), \\
w^{k+1} & = & w^k + \gamma d^k , \\
\alpha_j^{k+1} & = & \alpha_j^k - \gamma d^k.
\end{eqnarray*}
Using Lemma \ref{L:GLM hessians}, we see that
\begin{equation*}
g^k := \nabla f_j(w^k) - \alpha_j^k = \nabla R(w^k) + \phi'_j(\langle a_j, w^k \rangle) a_j - \alpha_j^k.
\end{equation*}
Still using Lemma \ref{L:GLM hessians}, and introducing the notation $\hat{a}_j = (\mI_d + \nabla^2 R(w^k) )^{-1} a_j$, we see that
\begin{eqnarray*}
\left(  \mI_d + \nabla^2 f_j(w^k) \right)^{-1} &=& 
\left( \mI_d + \nabla^2 R(w^k) \right)^{-1} - \frac{\phi_j''(\langle a_j,w^k \rangle)}{1 + \phi_j''(\langle a_j,w^k \rangle) \langle \hat a_j, a_j \rangle} \hat a_j \hat a_j^\top.
\end{eqnarray*}  
Therefore,
\begin{equation*}
d^k = \frac{\phi_j''(\langle a_j,w^k \rangle)\langle \hat a_j, g^k \rangle }{1 + \phi_j''(\langle a_j,w^k \rangle) \langle \hat a_j, a_j \rangle}\hat a_j
- \left( \mI_d + \nabla^2 R(w^k) \right)^{-1}g^k,
\end{equation*}
which concludes the proof.
\end{proof}

Finally, when the regularizer is the L2 norm, then we can implement SAN  even more efficiently as follows.
\begin{example}[Ridge regularization]
If $R(w) = \frac{\lambda}{2} \Vert w \Vert^2$ with $\lambda >0$, then the stochastic Newton direction $d^k$ can be computed explicitly (see also Lemma \ref{L:GLM hessians}):
\begin{equation*}
d^k = 
\frac{\phi_j''(r^k)}{1+\lambda}\cdot\frac{\langle a_j, \alpha_j^k \rangle -\phi_j'(r^k)\norm{a_j}^2-\lambda r^k}{1+\lambda+\phi_j''(r^k)\norm{a_j}^2}a_j
-
\frac{1}{1+\lambda}\left(\lambda w^k  + \phi_j'(r^k)a_j - \alpha_j^k\right),
\end{equation*}
where $r^k = \langle a_j, w^k \rangle $.
\end{example}

\subsection{SANA with GLMs}

In Algorithm~\ref{algo:SANA_GLM} we give the specialized implementation of SANA (Algorithm~\ref{algo:SANA}) for GLMs.
\begin{algorithm}[H]
\caption{{SANA for regularized GLMs}}\label{algo:SANA_GLM}
\KwIn{Data $\{a_i \}_{i=1}^n$, loss functions $\{\phi_i\}_{i=1}^n$, regularizer $R$, 
step size $\gamma\in (0,\,1]$, max iteration $T$}

Initialize $\alpha^0_1, \cdots, \alpha^0_n, w^0 \in \R^d$, with $\sum_{i=1}^n \alpha_i^0 =0$;

Pre-compute $\mu = \frac{n-1}{n}$; 

\For{$k=0, \dots, T-1$}{
    Sample $j \in \{1,\ldots, n\}$ uniformly;

    $g^k  = \nabla R(w^k) + \phi'_j(\langle a_j, w^k \rangle) a_j - \alpha_j^k$

    $\hat{a}_j = (\mu \mI_d + \nabla^2 R(w^k))^{-1} a_j$

    $d^k = \frac{\phi_j''(\langle a_j,w^k \rangle)\langle \hat a_j, g^k \rangle }{1 + \phi_j''(\langle a_j,w^k \rangle) \langle \hat a_j, a_j \rangle}\hat a_j
    - \left( \mu \mI_d + \nabla^2 R(w^k) \right)^{-1}g^k$
    
    $w^{k+1} = w^k + \gamma d^k$
    
    $\alpha^{k+1}_{j} = \alpha^k_j + \gamma \mu d^k$
    
    $\alpha_i^{k+1}  =  \alpha_i^k -  \frac{\gamma}{n} d^k, \quad \text{ for } i \neq j$
}
\KwOut{Last iterate $w^{T}$}
\end{algorithm}

Next we formalize the costs of Algorithm~\ref{algo:SANA_GLM} in the following remark. 
\begin{remark}
The costs of SANA (Algorithm~\ref{algo:SANA_GLM}) per iteration under Assumption~\ref{Ass:GLM regularized} are:
\begin{itemize}[nosep]
  \item \texttt{ Memory storage} of $\cO(nd)$ scalars.
  \item \texttt{  Memory access} of $\cO(nd)$. 
  \item \texttt{ Data access} of $\cO(1)$  .
  \item \texttt{  Computational cost} of $\cO(nd)$.
\end{itemize}
%
\end{remark}
%
%
%
%

\begin{lemma}\label{L:SANA for GLMs}
The SANA Algorithm \ref{algo:SANA} applied to Regularized GLMs (in the sense of Assumption \ref{Ass:GLM regularized}) is Algorithm \ref{algo:SANA_GLM}.
\end{lemma}
  
\begin{proof}
Let $k \in \{0, \dots, T-1\}$, and $\mu := 1 - n^{-1}$.
Let $j$ be sampled over $\{1,\ldots, n\}$ uniformly. From Algorithm \ref{algo:SANA} we have
\begin{eqnarray*}
g^k & = & \nabla f_j(w^k)-\alpha_j^k, \\
d^k & = & -\left(\mu \mI_d +  \nabla^2 f_j(w^k) \right)^{-1} g^k, \\ 
w^{k+1} & = & w^k +  \gamma d^k,\\
\alpha_j^{k+1} & = & \alpha_j^k - \gamma \mu d^k, \\
\alpha_i^{k+1} & = & \alpha_i^k +  \tfrac{\gamma}{n} d^k, \quad \text{ for } i \neq j.
\end{eqnarray*}
Using Lemma \ref{L:GLM hessians}, we see that
\begin{equation*}
g^k  = \nabla R(w^k) + \phi'_j(\langle a_j, w^k \rangle) a_j - \alpha_j^k.
\end{equation*}
Still using Lemma \ref{L:GLM hessians}, and introducing the notation $\hat{a}_j = (\mu \mI_d + \nabla^2 R(w^k) )^{-1} a_j$, we have that
\begin{eqnarray*}
\left(  \mu \mI_d + \nabla^2 f_j(w^k) \right)^{-1} &=& 
\left( \mu \mI_d + \nabla^2 R(w^k) \right)^{-1} - \frac{\phi_j''(\langle a_j,w^k \rangle)}{1 + \phi_j''(\langle a_j,w^k \rangle) \langle \hat a_j, a_j \rangle} \hat a_j \hat a_j^\top.
\end{eqnarray*}  
Therefore,
\begin{equation*}
d^k = \frac{\phi_j''(\langle a_j,w^k \rangle)\langle \hat a_j, g^k \rangle }{1 + \phi_j''(\langle a_j,w^k \rangle) \langle \hat a_j, a_j \rangle}\hat a_j
- \left( \mu \mI_d + \nabla^2 R(w^k) \right)^{-1}g^k,
\end{equation*}
which concludes the proof.

\end{proof}

\section{Experimental details in Section~\ref{sec:exp} and additional experiments} 
\label{sec:exp2}

We present the details of the experiments in Section~\ref{sec:exp},
in order to guide readers to reproduce the exact same results in Figure~\ref{fig:L2SAN} and Figure~\ref{fig:pseudo-HuberSAN}. 
We also explain some grid search results about sensitivity of hyperparameters in Section~\ref{sec:gridsearch}, showing in particular that SAN does not require parameter tuning. 
Then we provide additional experiments for SANA, SNM and IQN in Section~\ref{sec:exp_extra} which are not included in our main paper. These experimental results support that SANA introduced in Section~\ref{sec:SANA} is also a reasonable method. Finally, we provide experiments to compare SAN and SAN without variable metric in Section~\ref{sec:SANvsSANiI} to illustrate the importance of such variable metric.



\subsection{Experimental details in Section~\ref{sec:exp}}
\label{sec:expdetails}

All experiments in Section~\ref{sec:exp} were run in Python 3.7.7 on a laptop with an Intel Core i9-9980HK CPU and 32 Gigabyte of DDR4 RAM running OSX 11.3.1. 

All datasets were taken directly from LibSVM~\citep{libsvm} on \url{https://www.csie.ntu.edu.tw/~cjlin/libsvmtools/datasets/} and the scaled versions were used if provided. 
All datasets were preprocessed by adding an intercept, i.e. a constant feature one. 
For the datasets whose binary labels are not in $\{-1, 1\}$, e.g., \texttt{phishing}, \texttt{mushrooms} and \texttt{covtype},
we assigned $-1$ to the smallest labels and $+1$ to those largest ones.
All learnable parameters were initialized by zeros, e.g., $w^0 = 0 \in \R^d$ and $\alpha_i^0 = 0 \in \R^d$ for $i = 1, \cdots, n$ for SAN.

Table~\ref{tab:binary_datasets} provides the details of the datasets we used in Section~\ref{sec:exp}, including the condition number and $L_{\max}$. 
For a given dataset, let $\mA = [a_1 \cdots a_n] \in \R^{d \times n}$ be the data matrix,
the condition number in Table~\ref{tab:binary_datasets} is computed by
\[\mbox{condition number of the dataset} \; \eqdef \; \sqrt{\frac{\lambda_{\max}\left(\mA\mA^\top\right)}{\lambda_{\min}^{+}\left(\mA\mA^\top\right)}},\]
where  $\lambda_{\max}$ and $\lambda_{\min}^{+}$ are the largest and smallest non-zero eigenvalue operators respectively. 
$L_{\max}$ is defined as $L_{\max} = \max_{i = 1, \dots, n } L_i$, where $L_i = \frac{1}{4} \norm{a_i}^2 + \lambda$
is the smoothness constant
of the regularized logistic regression $f_i$. Notice that the step size's choice for SAG and SVRG is of the order of $\cO(1/L_{\max})$. 

From Table~\ref{tab:binary_datasets}, note that we have datasets that are middle scale (top row of Figure~\ref{fig:L2SAN}) and large scale (bottom row of Figure~\ref{fig:L2SAN}), well conditioned (\texttt{ijcnn1}) and ill conditioned (\texttt{webspam} and \texttt{rcv1}), sparse (\texttt{rcv1} and \texttt{real-sim}) and dense (\texttt{epsilon}), under-parametrized (\texttt{phishing}, \texttt{mushrooms}, \texttt{ijcnn1}, \texttt{covtype}, \texttt{webspam}, \texttt{epsilon} and \texttt{real-sim}) and over-parametrized (\texttt{rcv1}).

\begin{table}
\caption{Details of the binary data sets used in the logistic regression experiments}
\label{tab:binary_datasets}
\centering
\begin{tabular}{ llllll }
   \toprule
   dataset   & dimension ($d$)    & samples ($n$)  & $L_{\max}$            & sparsity & condition number \\
 \midrule
 phishing    & $68+1$             & $11055$        & $0.5001$              & $0.5588$ & $4.1065 \times 10^{18}$ \\
 mushrooms   & $112+1$            & $8124$         & $5.5001$              & $0.8125$ & $1.3095 \times 10^{19}$ \\
 ijcnn1      & $22+1$             & $49990$        & $1.2342$              & $0.4091$ & $25.6587$               \\
 covtype     & $54+1$             & $581012$       & $2.154$               & $0.7788$ & $9.6926 \times 10^{17}$ \\
 webspam     & $254+1$            & $350000$       & $0.5$                 & $0.6648$ & $6.9973\times 10^{255}$ \\
 epsilon     & $2000+1$           & $400000$       & $0.5$                 & $0.0$    & $3.2110 \times 10^{10}$ \\
 rcv1        & $47236+1$          & $20242$        & $0.5$                 & $0.9984$ & $5.3915 \times 10^{25}$ \\
 real-sim    & $20958+1$          & $72309$        & $0.5$                 & $0.9976$ & $1.3987 \times 10^{20}$ \\
 \bottomrule
\end{tabular}
\end{table}

\paragraph*{Pseudo-Huber function.} Recall the definition of the pseudo-Huber function used as the regularizer in our experiments in Figure~\ref{fig:pseudo-HuberSAN}: 
$R(w) = \sum_{i=1}^dR_i(w_i)$ with
\[R_i(w_i) =\delta^2 \left( \sqrt{1 + \left(\frac{w_i}{\delta} \right)^2} - 1 \right).\]
When $w_i$ is large, $R_i(w_i) \rightarrow \delta|w_i|$ for all $i = 1, \cdots, n$, i.e. $R(w)$ approximates L1 loss with a factor $\delta$; when $w_i$ is closed to zero$, R_i(w_i) \rightarrow \frac{1}{2}w_i^2$ for all $i = 1, \cdots, n$, i.e. $R(w)$ approximates L2 loss. This function can be served as a regularizer to promote the sparsity of the solution~\citep{FountoulakisG16}. 



Besides, the pseudo-Huber is $\cC^{\infty}$. The gradient of the pseudo-Huber is given by
\[\nabla R(w) \; = \; \left[\frac{w_1}{\sqrt{1+\left(\frac{w_1}{\delta}\right)^2}} \ \cdots \ \frac{w_d}{\sqrt{1+\left(\frac{w_d}{\delta}\right)^2}}\right]^\top \in \R^d\]
and the Hessian is given by
\[\nabla^2R(w) \; = \; \Diag{\left(1+\left(\frac{w_1}{\delta}\right)^2\right)^{-3/2}, \cdots, \left(1+\left(\frac{w_d}{\delta}\right)^2\right)^{-3/2}} \; \leq \; \mI_d.\]
Thus the pseudo-Huber is $1$-smooth which is the same as L2 regularizer. Consequently, $L_{\max}$ for the pseudo-Huber regularized logistic regression is the same as the L2-regularized one.

\subsection{Function sub-optimality plots}
\label{sec:sub_opt}

\begin{figure}
\centering
\includegraphics[width=.24\linewidth]{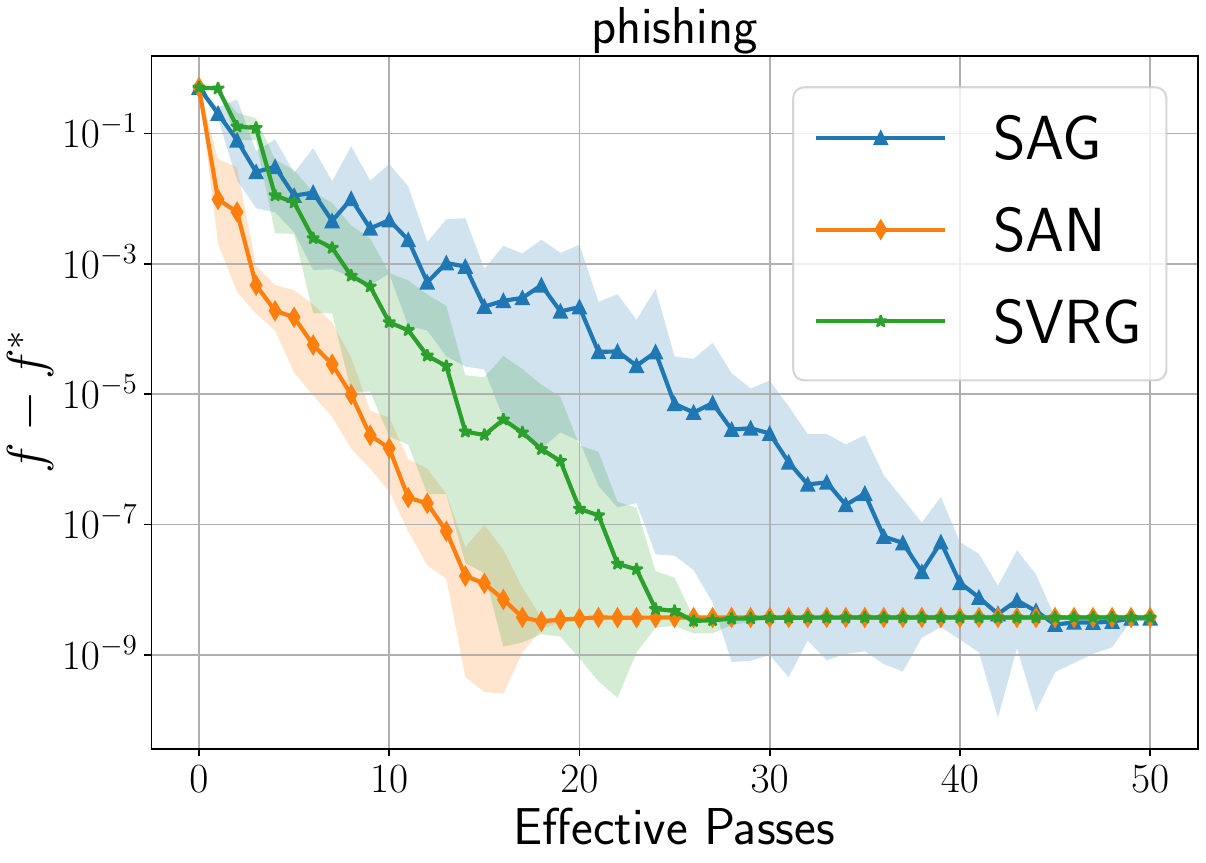}
\includegraphics[width=.24\linewidth]{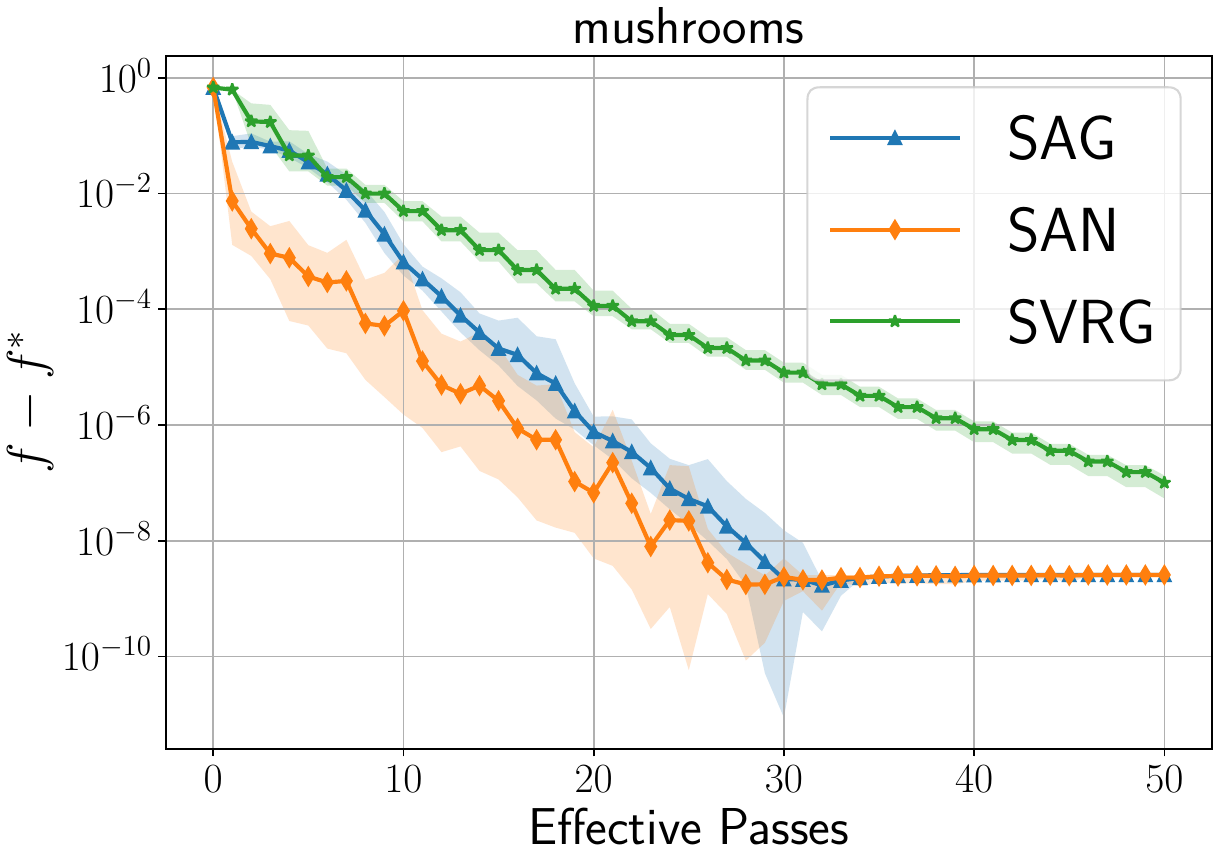}
\includegraphics[width=.24\linewidth]{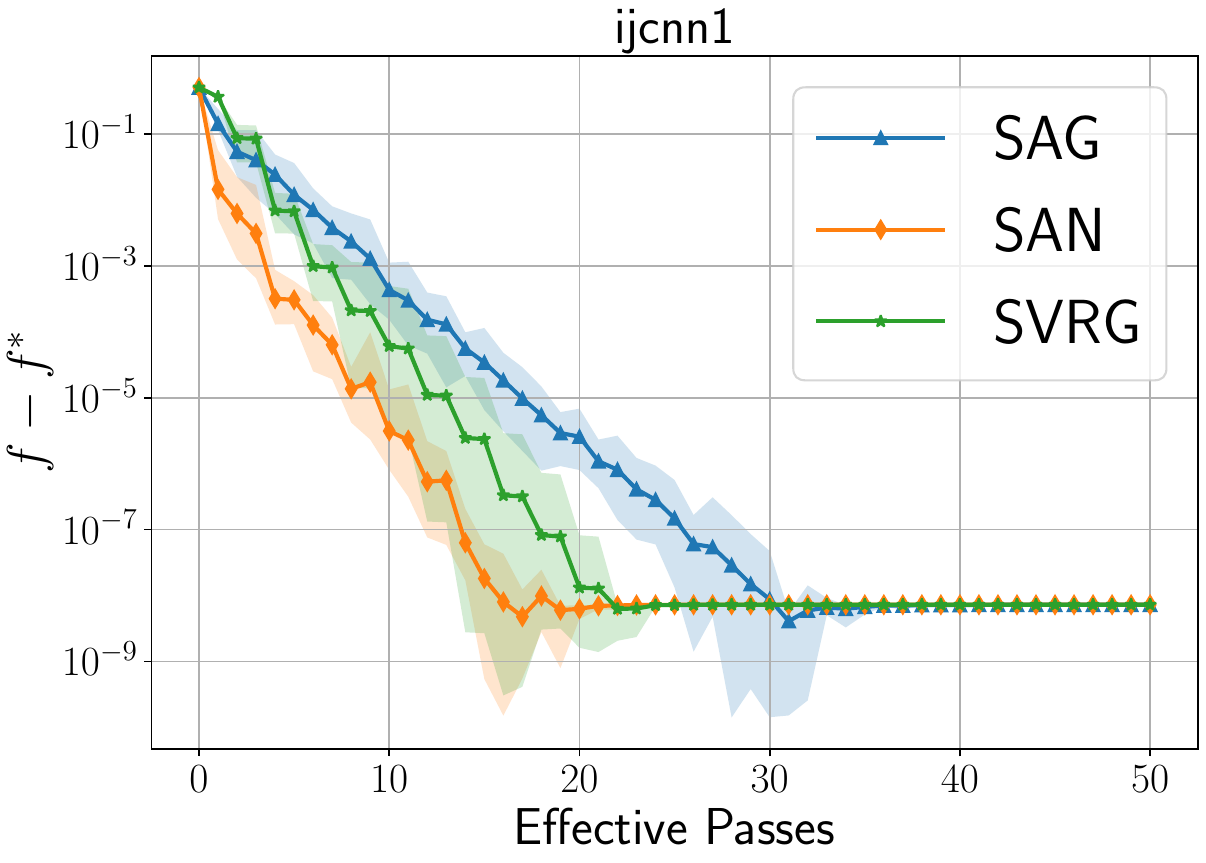}
\includegraphics[width=.24\linewidth]{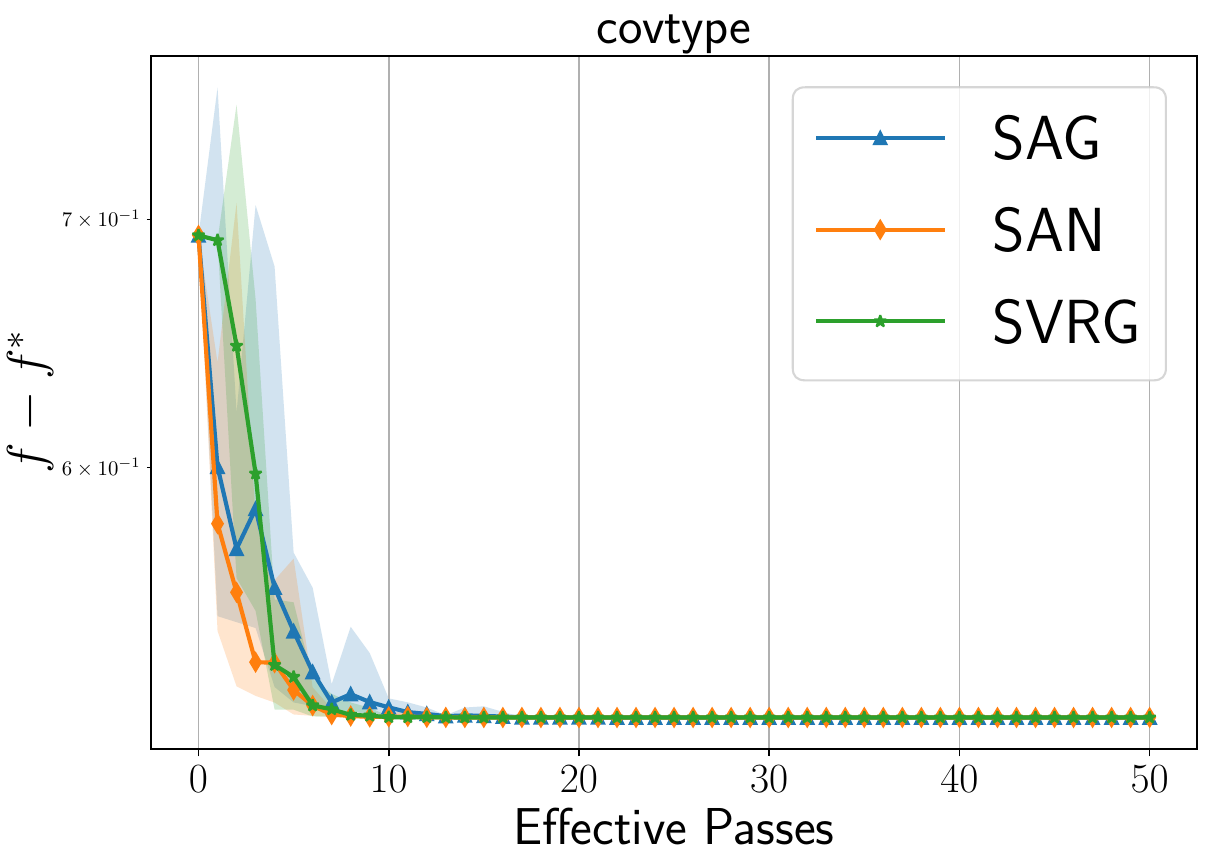} 
\caption{Function sub-optimality of logistic regression with L2 regularization.}
\label{fig:sub}
\end{figure}

The performance of an algorithm for solving a convex problem is usually done by measuring one of the following quantities: the solution gap $\Vert w^k - w^* \Vert$ where $w^*$ is the solution of the problem, the optimization gap $f(w^k) - \inf f$, and the stationarity gap $\Vert \nabla f(w^k)\Vert$.
In this paper we choose to measure and compare performance of algorithms in terms of $\Vert \nabla f(w^k) \Vert$.
The main reason for this is that the solution gap $\Vert w^k - w^* \Vert$ and the optimization gap $f(w^k) - \inf f$ both require to compute the solution of the problem to a high precision.
While this is possible to do for small problems, it quickly becomes intractable for large problems (see Figure \ref{fig:sub} for \texttt{covtype}), which we want to address in this paper (see Table \ref{tab:binary_datasets} for more bigger datasets than \texttt{covtype}).
The flat curves appeared in Figure \ref{fig:sub} after certain effective passes, especially for \texttt{covtype} dataset, are due to the imprecise computation of $\inf f$ from the solver \verb|scipy.optimize.fmin_l_bfgs_b|. Indeed, the curves in Figure \ref{fig:sub} are in logarithmic scale. When $f(w^k) - \widehat{\inf f} < 0$ with $\widehat{\inf f}$ the tentative solution of the problem computed by the solver, it means that the solution $\widehat{\inf f}$ is imprecise, i.e. the solver performs worse than the tested algorithms. In this case, Figure \ref{fig:sub} plots $\left|f(w^k) - \widehat{\inf f}\right| =  \widehat{\inf f} - f(w^k) > 0$ where the curves remain flat in logarithmic scale.

We argue that the quantity $\Vert \nabla f(w^k) \Vert^2$ is a fair and good proxy for the more classical optimization gap $f(w^k) - \inf f$.
Our argument for this is twofold.
First, we observe empirically on small problems (for which we can compute $\inf f$ with precision) that the curves for $\Vert \nabla f(w^k) \Vert^2$ and $f(w^k) - \inf f$ behave the same (see Figure \ref{fig:sub}). 
Second, we verify theoretically that $\Vert \nabla f(w^k) \Vert^2$ and $f(w^k) - \inf f$ are of the same order.
Indeed, Assumption \ref{Ass:strict convexity of problem} implies that $f$ is strongly convex on every compact.
In particular, it verifies on every compact a Lojasiewicz inequality:
\begin{equation*}
(\forall R>0)(\exists \mu >0)(\forall w \in \mathbb{B}(0,R)) 
\quad
f(w) - \inf f \leq \frac{1}{2 \mu} \Vert \nabla f(w) \Vert^2.
\end{equation*}
Moreover, $f$ is convex, so if we assume that $f$ has a $L$-Lipschitz gradient, we obtain the following inequality:
\begin{equation*}
(\forall w \in \mathbb{R}^d) \quad \frac{1}{2L}\Vert \nabla f(w) \Vert^2 \leq f(w) - \inf f.
\end{equation*}
Note that this assumption is verified for the functions considered in our experiments (see Section \ref{sec:expdetails}).

\subsection{Effect of hyperparameters}\label{sec:gridsearch}

As we discussed in Section~\ref{sec:exp}, SAN involves neither  prior knowledge of the datasets (e.g., $L_{\max}$), 
nor the hyperparameter tuning, while both SAG~\citep{SAG} and SVRG~\citep{Johnson2013} do. 
To support this conclusion, under different 
hyperparameters setting, we measure the performance of the given algorithm 
by monitoring the number of effective passes over the data required to reach below 
a threshold (e.g., $10^{-4}$ in our case) of $\norm{\nabla f}$.
We repeat this procedure $5$ times and report the average results. 

\paragraph*{Grid search for SAN.}

SAN has two hyperparameters: the probability $\pi$ doing the averaging step in Algorithm~\ref{algo:SAN} 
and the step size $\gamma$. We searched $\pi$ in a wide range, 
$\pi \in \{\frac{1}{2n},  \frac{1}{n}, \frac{10}{n}, \frac{100}{n}, \frac{1000}{n} \}$; as for $\gamma$, 
through our extensive experiments, we observed that SAN works consistently well 
when $\gamma$ is around one as we expected for second order methods,
thus we tried $\gamma \in \{0.7, 0.8, 0.9, 1.0, 1.1, 1.2, 1.3 \}$.

\begin{table}
  \caption{\texttt{covtype} dataset: grid search of $\pi$ and $\gamma$ for SAN}
  \label{tab: gridsearchsancovtype}
  \centering
    \begin{tabular}{llllllll}
    \toprule
    \diagbox{$\pi$}{$\gamma$} &  $0.7$  &  $0.8$ & $0.9$ & $1.0$ & $1.1$ & $1.2$ & $1.3$ \\
    \midrule
    $1/2n$ & $27$ & $25$ & $23$ & $21$ & $21$ & $22$ & $24$ \\
    $1/n$ & $26$ & $26$ & $25$ & $22$ & $22$ & $23$ & $24$ \\
    $10/n$ & $28$ & $24$ & $24$ & $23$ & $23$ & $22$ & $22$ \\
    $100/n$ & $28$ & $26$ & $23$ & $23$ & $22$ & $23$ & $24$ \\
    $1000/n$ & $28$ & $26$ & $27$ & $27$ & $25$ & $24$ & $26$ \\
    \bottomrule
    \end{tabular}
\end{table}

\begin{table}
  \caption{\texttt{ijcnn1} dataset: grid search of $\pi$ and $\gamma$ for SAN}
  \label{tab: gridsearchsanijcnn1}
  \centering
    \begin{tabular}{llllllll}
    \toprule
    \diagbox{$\pi$}{$\gamma$} &  $0.7$  &  $0.8$ & $0.9$ & $1.0$ & $1.1$ & $1.2$ & $1.3$ \\
    \midrule
    $1/2n$ & $13$ & $13$ & $14$ & $12$ & $11$ & $12$ & $13$ \\
    $1/n$ & $14$ & $13$ & $12$ & $12$ & $12$ & $12$ & $13$ \\
    $10/n$ & $13$ & $13$ & $12$ & $12$ & $12$ & $12$ & $13$ \\
    $100/n$ & $13$ & $11$ & $12$ & $13$ & $11$ & $12$ & $14$ \\
    $1000/n$ & $16$ & $13$ & $13$ & $13$ & $14$ & $14$ & $14$ \\
    \bottomrule
    \end{tabular}
\end{table}


Table~\ref{tab: gridsearchsancovtype} and~\ref{tab: gridsearchsanijcnn1} show the grid search results on
datasets \texttt{covtype} and \texttt{ijcnn1}. We see that the average effective data passes required to reach
the threshold is stable. It means that SAN is not sensitive to these hyperparameters. This 
advantage allows us to use $\pi = \frac{1}{n+1}$ and $\gamma = 1.0$ as default choice in our experiments
shown in Section~\ref{sec:exp}.

\paragraph*{Grid search for SAG and SVRG.}
Additionally we evaluated the effect of step size $\gamma$ which is a crucial hyperparameter for first order methods.
Let $f_i$ be $L_i$-smooth for all $i \in \{1, \dots, n \}$ and $L_{\max} = \max_{i \in \{1, \dots, n \}} L_i$.
As $\gamma = \frac{1}{L_{\max}}$ is thought as the rule of thumb choice in practice~\citep{scikit-learn} for SAG and SVRG, we searched over the values given by 
\[\gamma \in \left\{ \frac{1}{10L_{\max}}, \frac{1}{5L_{\max}}, \frac{1}{3L_{\max}}, \frac{1}{2L_{\max}}, 
\frac{1}{L_{\max}}, \frac{2}{L_{\max}}, \frac{5}{L_{\max}} \right\}\] on different datasets. 

\begin{table}[!h]
  \caption{Grid search of the step size $\gamma$ for SVRG on four datasets}
  \label{tab: gridsearchsvrg}
  \centering
    \begin{tabular}{llllllll}
    \toprule
    \diagbox{Datasets}{$\gamma$} &  $\frac{1}{10L_{\max}}$  &  $\frac{1}{5L_{\max}}$ & $\frac{1}{3L_{\max}}$ & 
    $\frac{1}{2L_{\max}}$ & $\frac{1}{L_{\max}}$ & $\frac{2}{L_{\max}}$ & $\frac{5}{L_{\max}}$  \\
    \midrule
    \texttt{covtype} & $44$ & $24$ & $18$ & $\boldsymbol{14}$ & $20$ & $\times$ & $\times$ \\
    \texttt{ijcnn1} & $22$ & $12$ & $\boldsymbol{10}$ & $\boldsymbol{10}$ & $15$ & $25$ & $\times$ \\
    \texttt{phishing} & $14$ & $11$ & $\boldsymbol{10}$ & $14$ & $18$ & $44$ & $\times$  \\
    \texttt{mushrooms} & $\times$ & $\times$ & $44$ & $36$ & $28$ &$\boldsymbol{20}$ & $46$ \\
    \bottomrule
    \end{tabular}
\end{table}

\begin{table}[!h]
  \caption{Grid search of the step size $\gamma$ for SAG on four datasets}
  \label{tab: gridsearchsag}
  \centering
    \begin{tabular}{llllllll}
    \toprule
    \diagbox{Datasets}{$\gamma$} &  $\frac{1}{10L_{\max}}$  &  $\frac{1}{5L_{\max}}$ & $\frac{1}{3L_{\max}}$ & 
    $\frac{1}{2L_{\max}}$ & $\frac{1}{L_{\max}}$ & $\frac{2}{L_{\max}}$ & $\frac{5}{L_{\max}}$  \\
    \midrule
    \texttt{covtype} & $21$ & $\boldsymbol{19}$ & $23$ & $24$ & $40$ & $\times$ & $\times$ \\
    \texttt{ijcnn1} & $\boldsymbol{14}$ & $16$ & $17$ & $17$ & $22$ & $34$ & $\times$ \\
    \texttt{phishing} & $\boldsymbol{14}$ & $17$ & $21$ & $21$ & $30$ & $48$ & $\times$ \\
    \texttt{mushrooms} & $\times$ & $47$ & $32$ & $24$ & $\boldsymbol{18}$ & $25$ & $\times$ \\
    \bottomrule
    \end{tabular}
\end{table}

From our observations to Table~\ref{tab: gridsearchsvrg} and~\ref{tab: gridsearchsag} ,\footnote{The symbol $\times$ 
in these tables means that the algorithm can not reach below the threshold $10^{-4}$ after $50$ data passes.} we can draw
the conclusions that compared to SAN, there is no universal step size choice for SAG and SVRG which gives a consistent good performance on different datasets. This point is one of our original motivations to develop a second order method
that requires neither prior knowledge from datasets nor the hyperparameter tuning.



\subsection{Additional experiments for SANA, SNM and IQN applied for L2 logistic regression}
\label{sec:exp_extra}

We present some additional results of SANA, SNM~\citep{SNM} and IQN~\citep{mokhtari2018iqn} compared to SAN on L2 logistic regression scenario. 

\begin{figure}[!h]
\centering
\includegraphics[width=.35\linewidth]{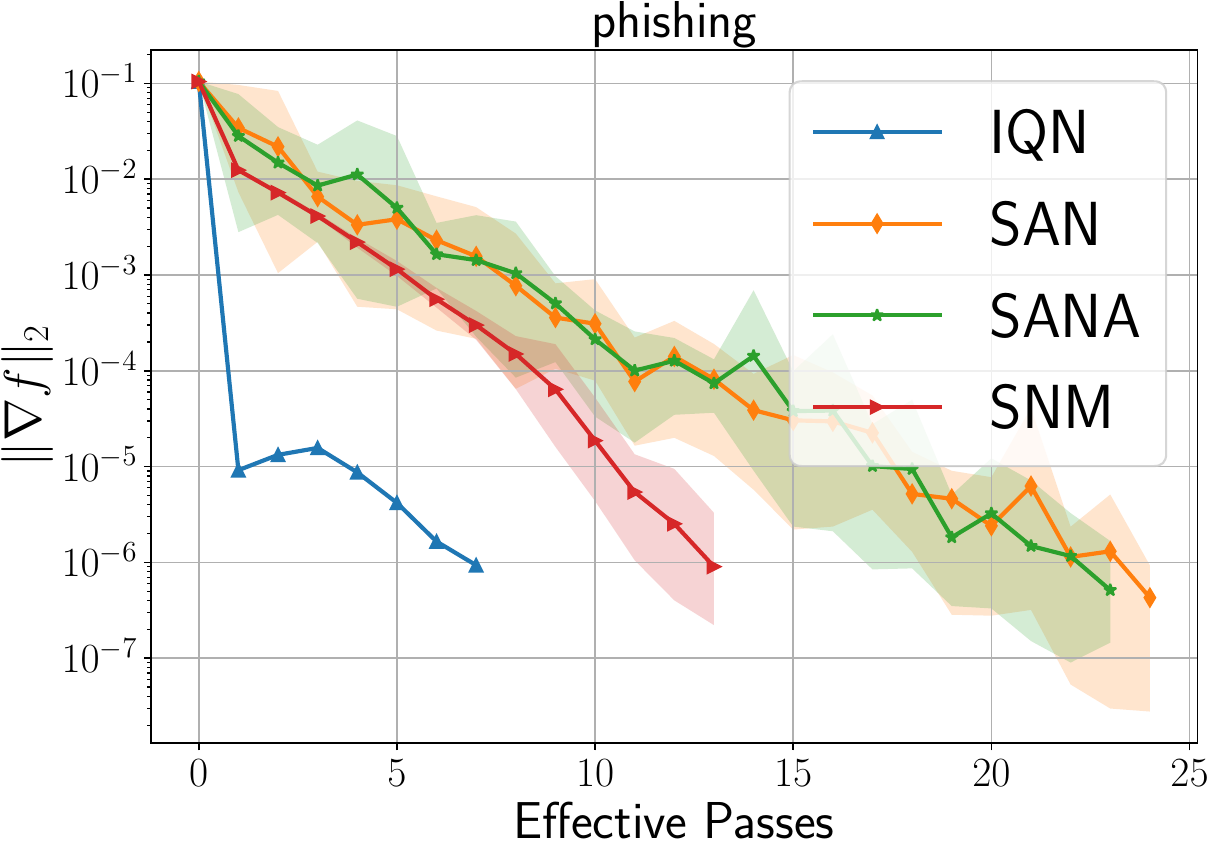}
\includegraphics[width=.35\linewidth]{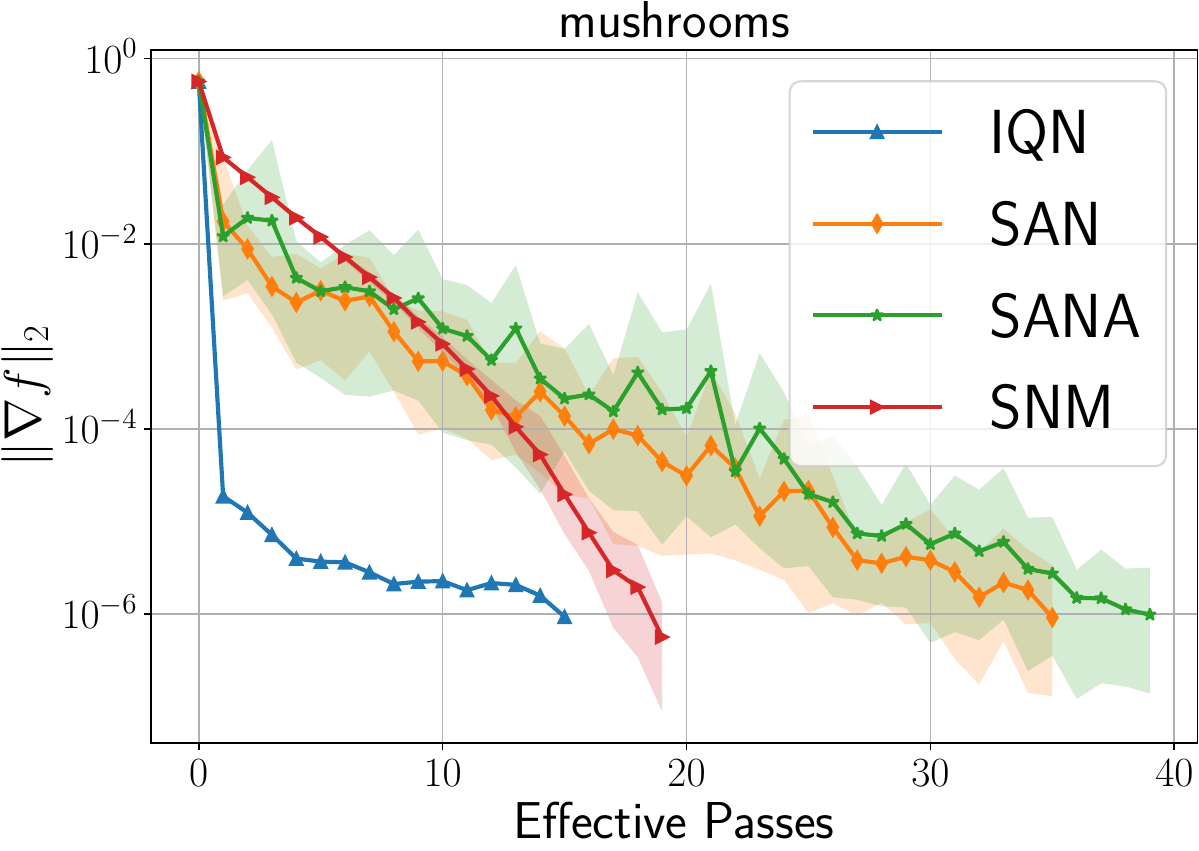} \\
\includegraphics[width=.35\linewidth]{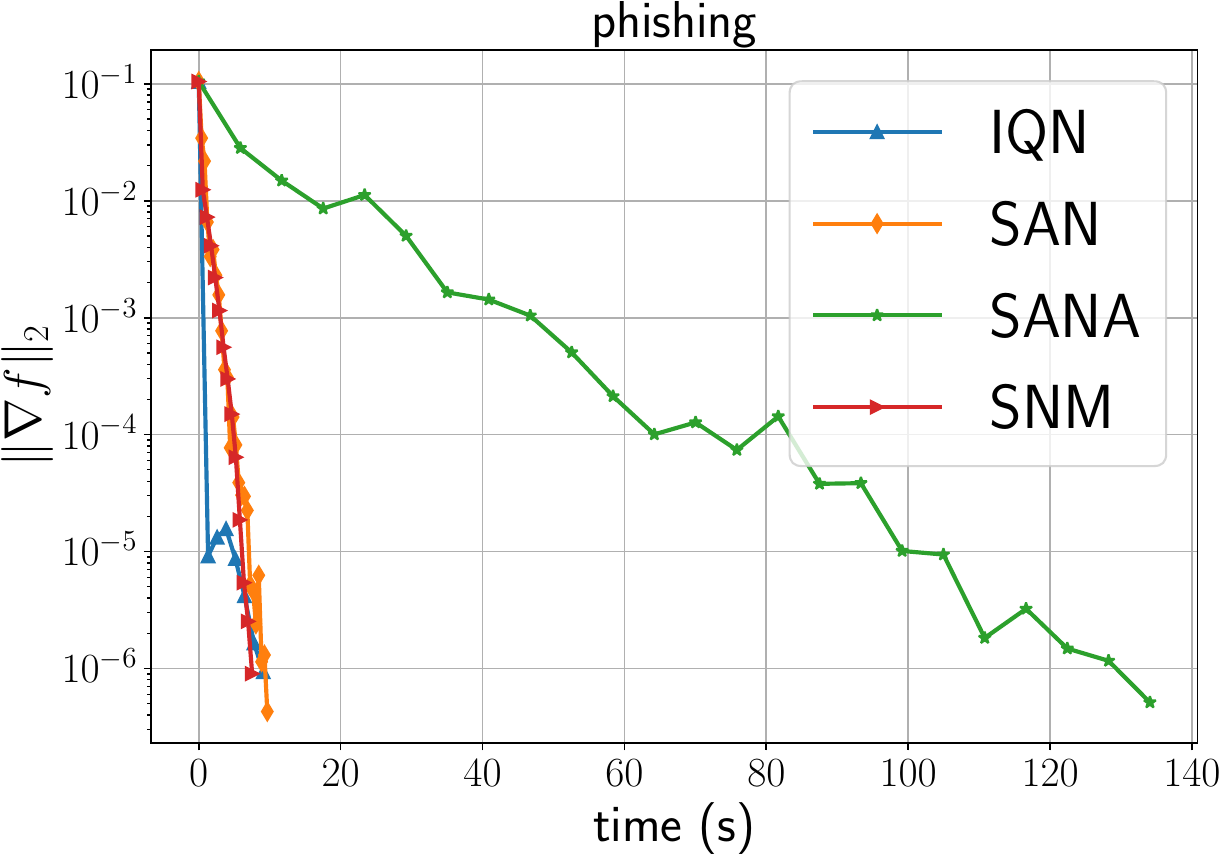}
\includegraphics[width=.35\linewidth]{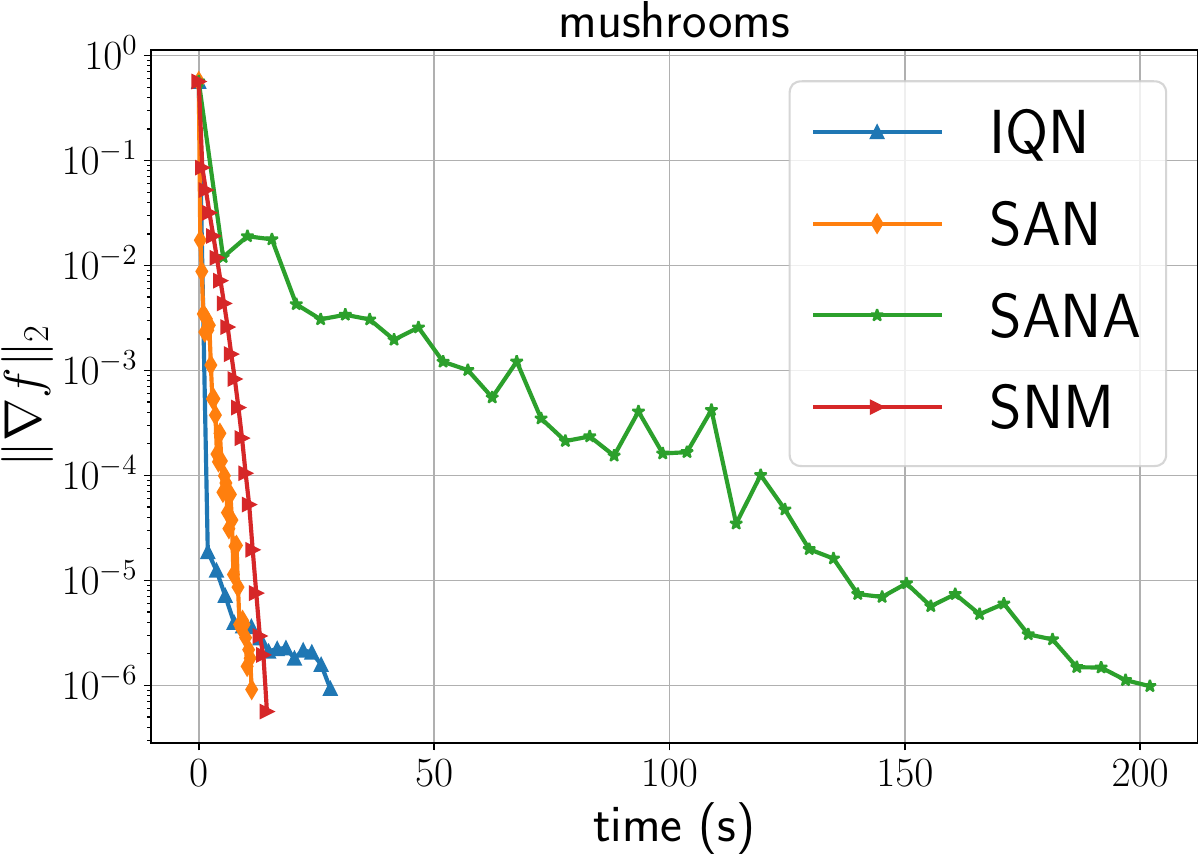}
\caption{L2-regularized logistic regression for SAN, SANA, SNM and IQN on middle size datasets. Top row is evaluated
in terms of effective data passes and bottom row is evaluated in terms of computational time.}
\label{fig:compare SANA L2}
\end{figure}

First, we show the results on middle size datasets, \texttt{phishing} and \texttt{mushrooms} in Figure~\ref{fig:compare SANA L2}.
On the one hand, in terms of effective passes of data, SANA has a similar performance as SAN despite the fact that SANA is unbiased and SAN is a biased estimate.
Both methods are less efficient than SNM and IQN. 
Notice that the initialization process of SNM is expensive, as it requires a computation of the full Newton system. Such process is not counted into the effective passes.
On the other hand, in terms of computational time, we observe that SAN does as well as IQN and SNM: SAN's cheap iteration cost compensates for its slower convergence rate.
On the other hand, we observe for SANA that it is not
competitive with respect to the other methods in terms of time taken. This is coherent in a regime where $d \ll n$ since  SANA has a computation cost of $\mathcal{O}(nd)$ per
iteration (see Table~\ref{tab:complexity}), while the cost for SAN and SNM, IQN is respectively $\cO(d)$ and $\cO(d^2)$, 
However, it shows that SANA is still 
a meaningful incremental second order method that satisfies our objective statement. This supports our general approach to design algorithms via function splitting.

\begin{figure}[!h]
\centering
\includegraphics[width=.35\linewidth]{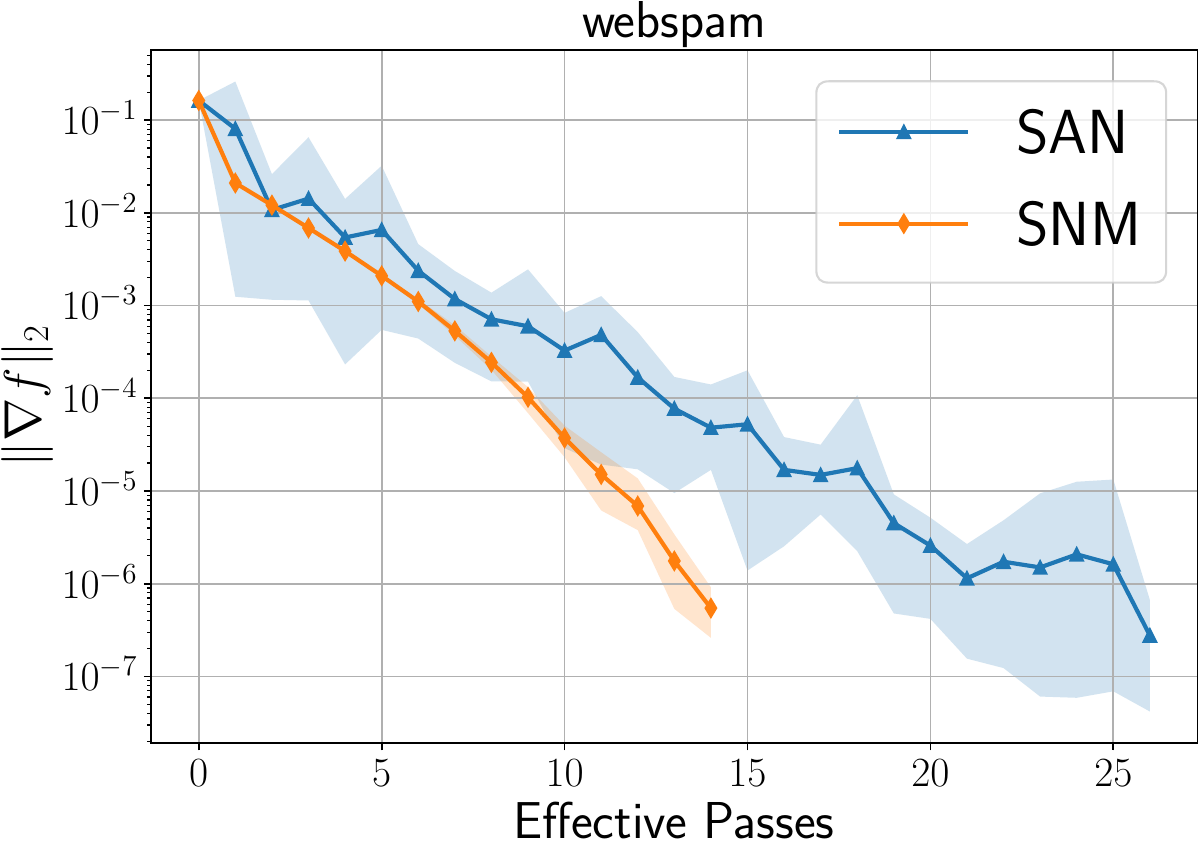} 
\includegraphics[width=.35\linewidth]{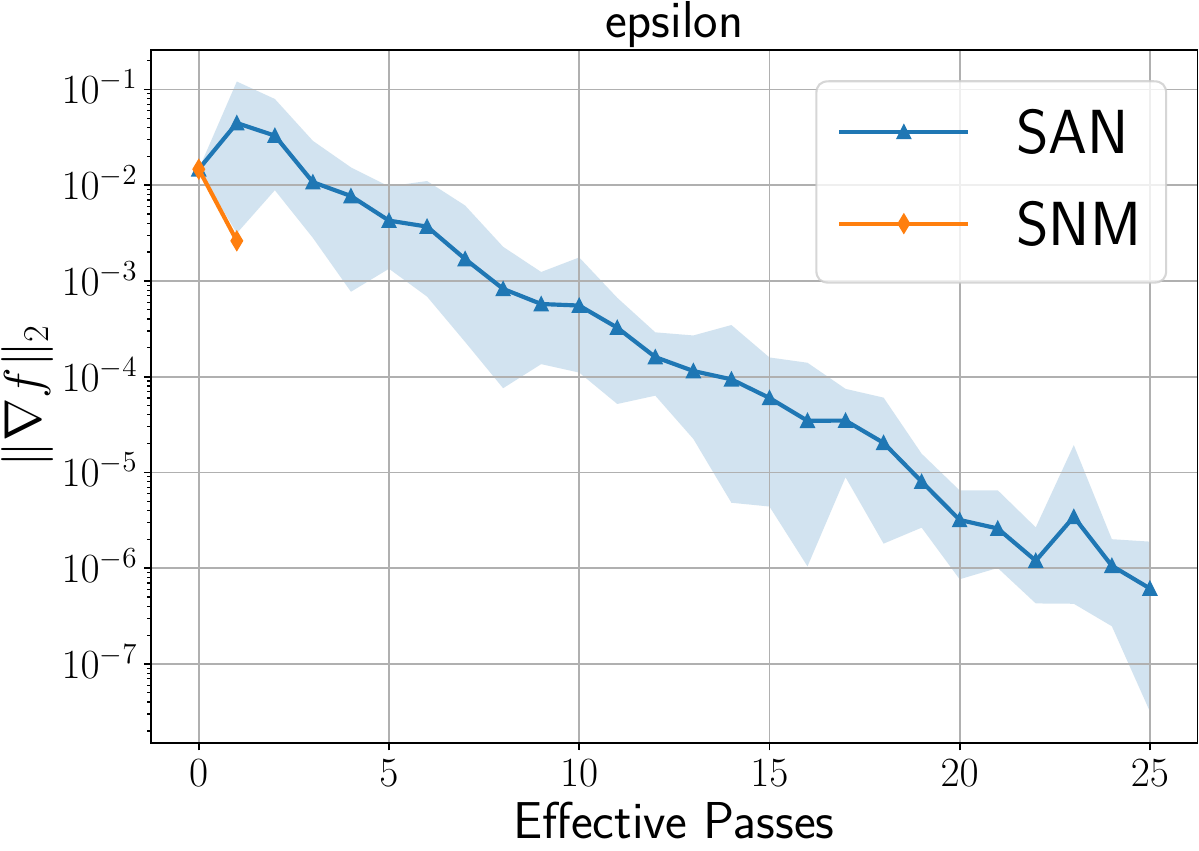}\\
\includegraphics[width=.35\linewidth]{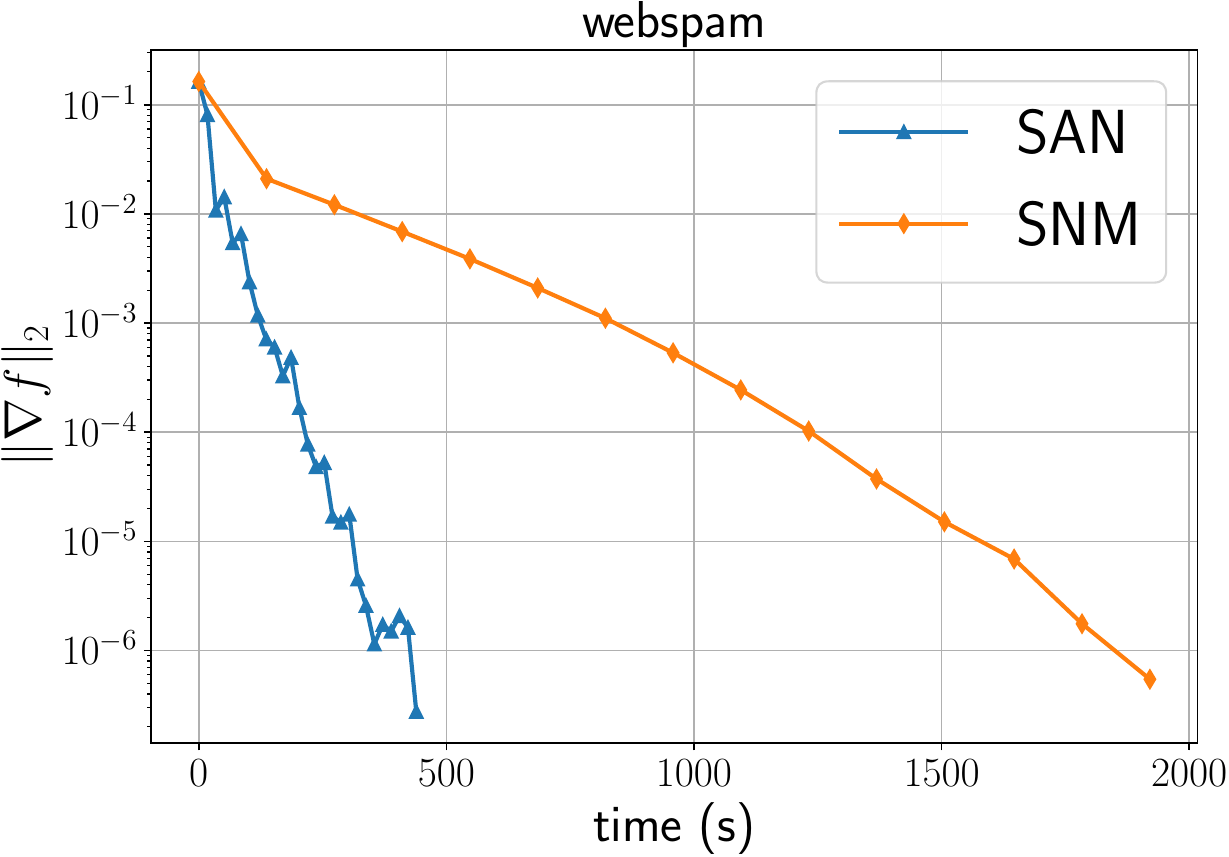}
\includegraphics[width=.35\linewidth]{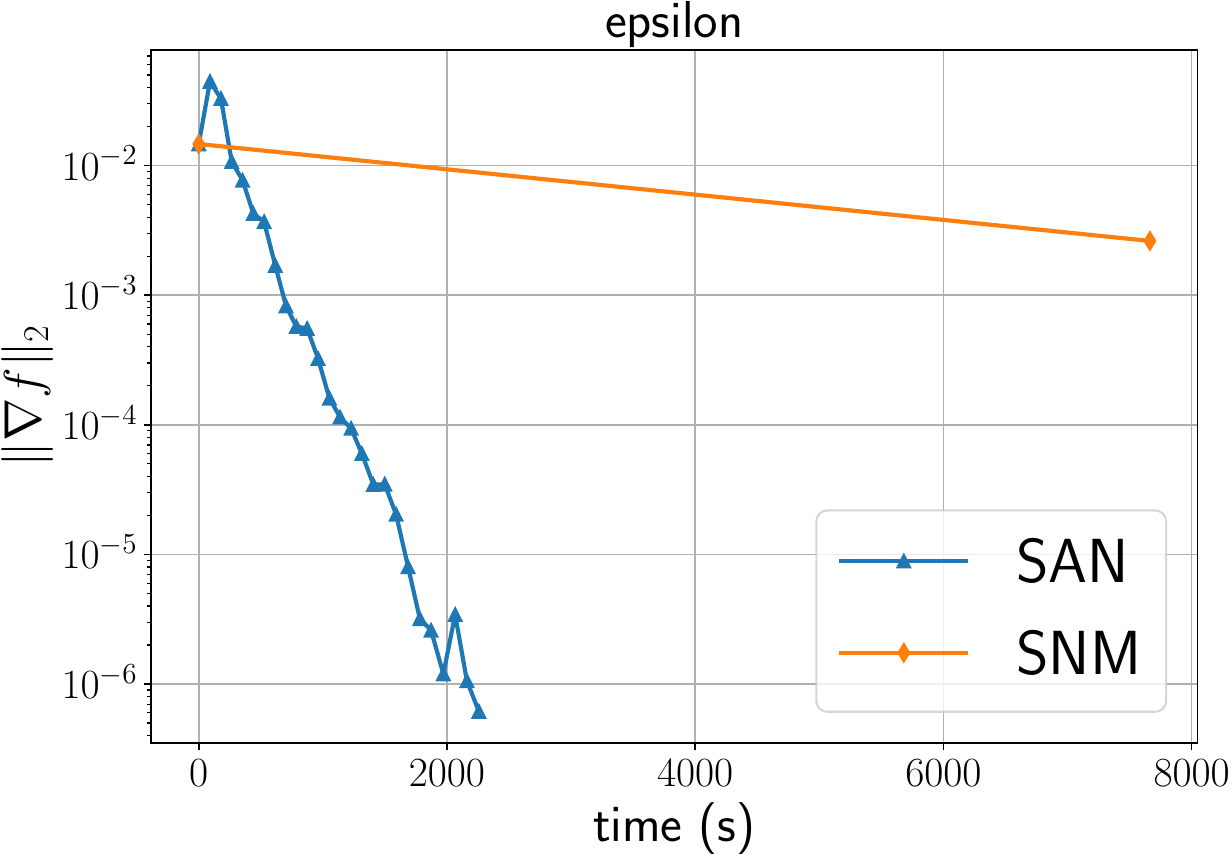}
\caption{L2-regularized logistic regression for SAN and SNM on large size datasets. Top row is evaluated in terms of effective data passes and bottom row is evaluated in terms of computational time.}
\label{fig:compare SAN and SNM}
\end{figure}

In our second set of experiments, we compare those algorithms on large scale datasets.
As shown in Figure~\ref{fig:compare SAN and SNM}, we tested two datasets \texttt{webspam} and \texttt{epsilon}.  As we discuss below, both SNM and IQN are limited in this case, while SAN is able to efficiently solve the problem. 
IQN is disqualified in this large scale setting, because its memory cost of $\cO(nd^2)$ is prohibitive and makes it impossible to run on a laptop.
This cost comes from the fact that IQN maintains and updates $n$ approximations of the hessians $\nabla^2 f_i(w^k)$, each of size $d^2$, and these matrices are not low-rank even for a GLM, preventing from using GLM implementation tricks (as it is the case for SNM, see Remark \ref{R:SNM implementation} below).
We also did not run SANA, since we already know that it performs similarly to SAN in terms of effective passes, but suffers from a cost per iteration scaling with $n$, which is too large here.
It is possible to run SNM, but it is not efficient in terms of computational time due to its expensive cost per iteration.
For the dataset \texttt{epsilon}, just after one pass over the data, the running time of SNM exceeded our maximum allowed time while at the same time SAN has run $25$ data passes and reached a solution with a $10^{-6}$ precision.

Furthermore, note that we are running experiments in a setting which is favorable to SNM.
Indeed, its cost per iteration $\mathcal{O}(d^2)$ is only valid when using $L2$ regularization.
If we were to consider another separable regularizer, its cost per iteration would be $\cO(d^3)$, making SNM infeasible for large dimensional problems.
The next remark details those considerations about the complexity of SNM.

\begin{remark}[On the cost of SNM]\label{R:SNM implementation}
The updates of SNM can be written in closed form as 
\begin{equation}\label{eq: snm update}
w^{k+1} = \left( \frac{1}{n} \sum_{i=1}^n \nabla^2 f_i(\alpha_i^k)  \right)^{-1} \left(\frac{1}{n} \sum_{i=1}^n \nabla^2 f_i(\alpha_i^k)\alpha_i^k - \nabla f_i(\alpha_i^k) \right), 
\ \alpha_j^{k+1} = w^{k+1},
\ \alpha_i^{k+1} = \alpha_i^k \text{ for } i \neq j,
\end{equation}
where $w, \alpha_i \in \R^d$ for $i=1,\cdots, n$ are variables defined in~\eqref{eq:weqai} using a variable splitting trick.
The main cost of SNM is to update the following inverse matrix
$\left( \frac{1}{n} \sum_{i=1}^n \nabla^2 f_i(\alpha_i^k)  \right)^{-1}$ after updating a single $\alpha_j$.

For L2-regularized GLMs, by using the \textit{Sherman-Morrison formula}, the above term can be implemented efficiently in $\mathcal{O}(d^2)$~(See Algorithm~$3$ in~\citet{SNM}), exploiting rank one updates of the matrix. 

For other separable regularizers, such a formula is no longer available, as the perturbation becomes rank $d$ due to the diagonal Hessian of the regularizer derived by Lemma~\ref{L:GLM hessians}~\ref{itm:diagonal Hessian}. 
The inversion of the matrix, therefore costs $\cO(d^3)$ over all. 
Note that the memory cost is also impacted in this case: for general separable regularizers the memory cost will be $\cO(nd+d^2)$, instead of $\cO(n+d^2)$ as can be seen in Table~\ref{tab:complexity} for L2-regularized GLMs.
\end{remark}




\subsection{SAN vs SAN without the variable metric}
\label{sec:SANvsSANiI}

One of the main design features of SAN is that at every iteration we project our iterates onto an affine space with respect to a metric induced by the Hessian of one sampled function.
One could ask whether this is worth it, given that it makes the theoretical analysis much more difficult.
Let us consider again the problem introduced in~\eqref{eq:SAN implicit linearized 1..n} where the Hessian induced norm has been replaced by the L2 norm as following:
\begin{eqnarray}
\alpha^{k+1}_j, w^{k+1} & =& \arg\min_{\alpha_j \in \R^d, w \in \R^d} \norm{\alpha_j - \alpha_j^k}^2 + \norm{w-w^k}^2 \label{eq:SAN implicit linearized 1..n without variable metric}
 \\  
& &
 \mbox{subject to } \nabla f_j(w^k)+\nabla^2 f_j(w^k)(w-w^k) = \alpha_j .\nonumber 
\end{eqnarray} 
Using Lemma~\ref{lem:Hproject}, we can compute the closed form update of~\eqref{eq:SAN implicit linearized 1..n without variable metric} (the details are left to readers):
\begin{eqnarray}\label{eq: SAN without variable metric}
\alpha_j^{k+1} &=& \alpha_j^k - \left(\mI + (\nabla^2 f_j(w^k))^2 \right)^{-1} \left(\alpha_j^k - \nabla f_j (w^k) \right), \\
w^{k+1} &=& w^k - \nabla^2 f_j(w^k) \left(\alpha_j^{k+1} - \alpha_j^k \right).
\end{eqnarray}

\begin{figure}[!h]
\centering
\includegraphics[width=.24\linewidth]{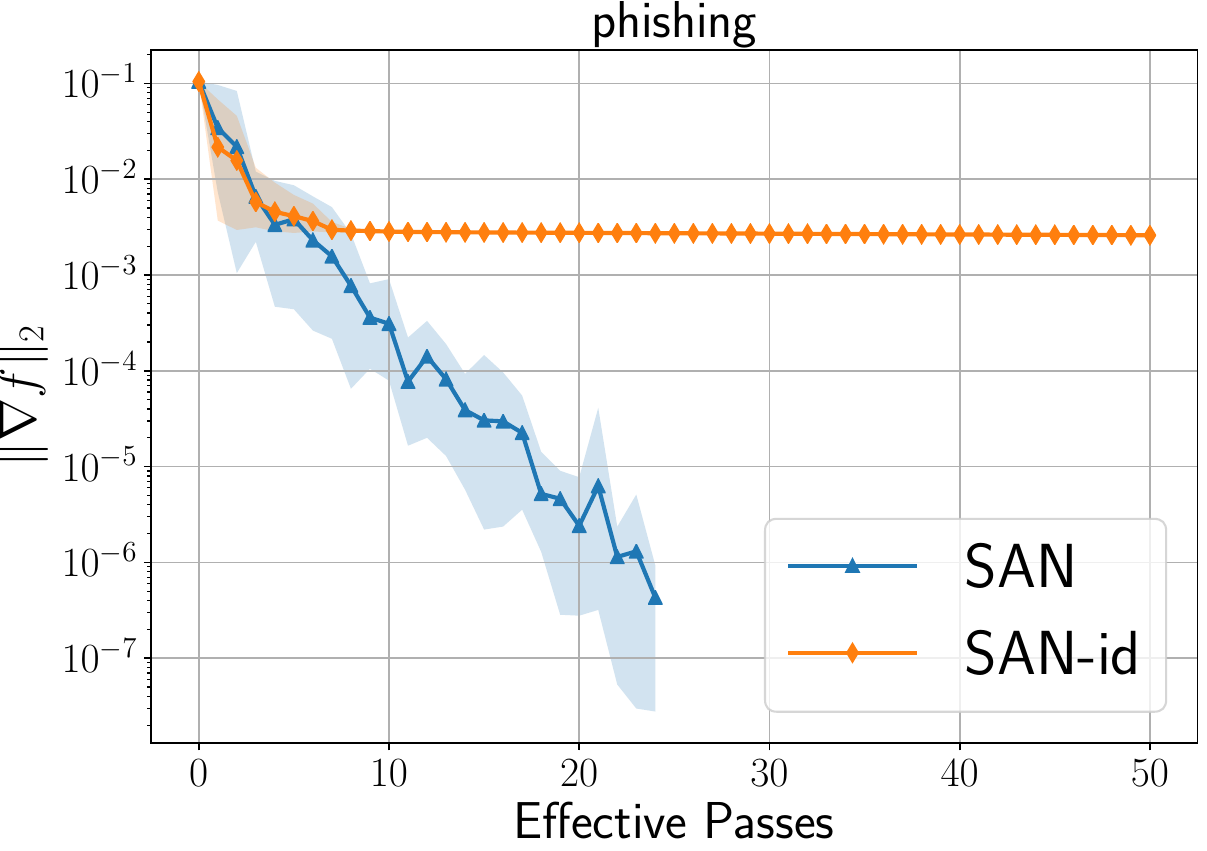}
\includegraphics[width=.24\linewidth]{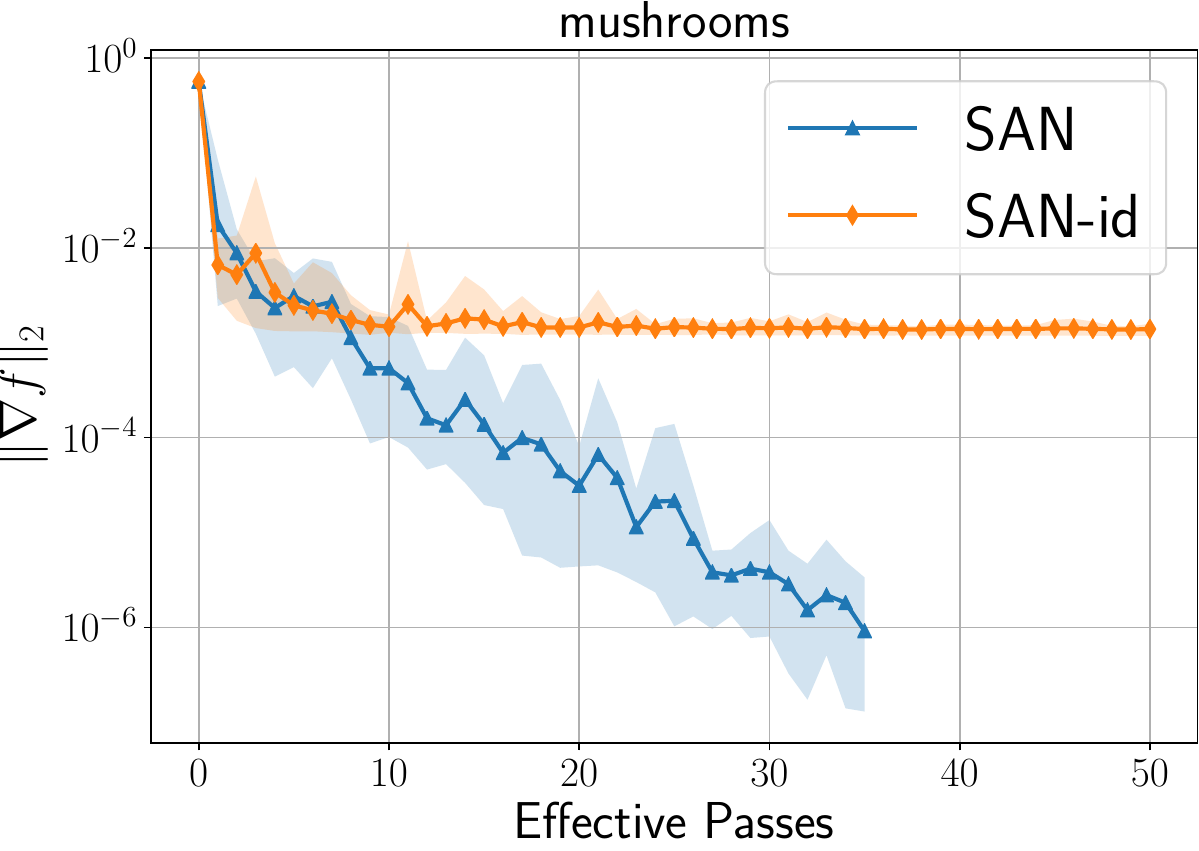}
\includegraphics[width=.24\linewidth]{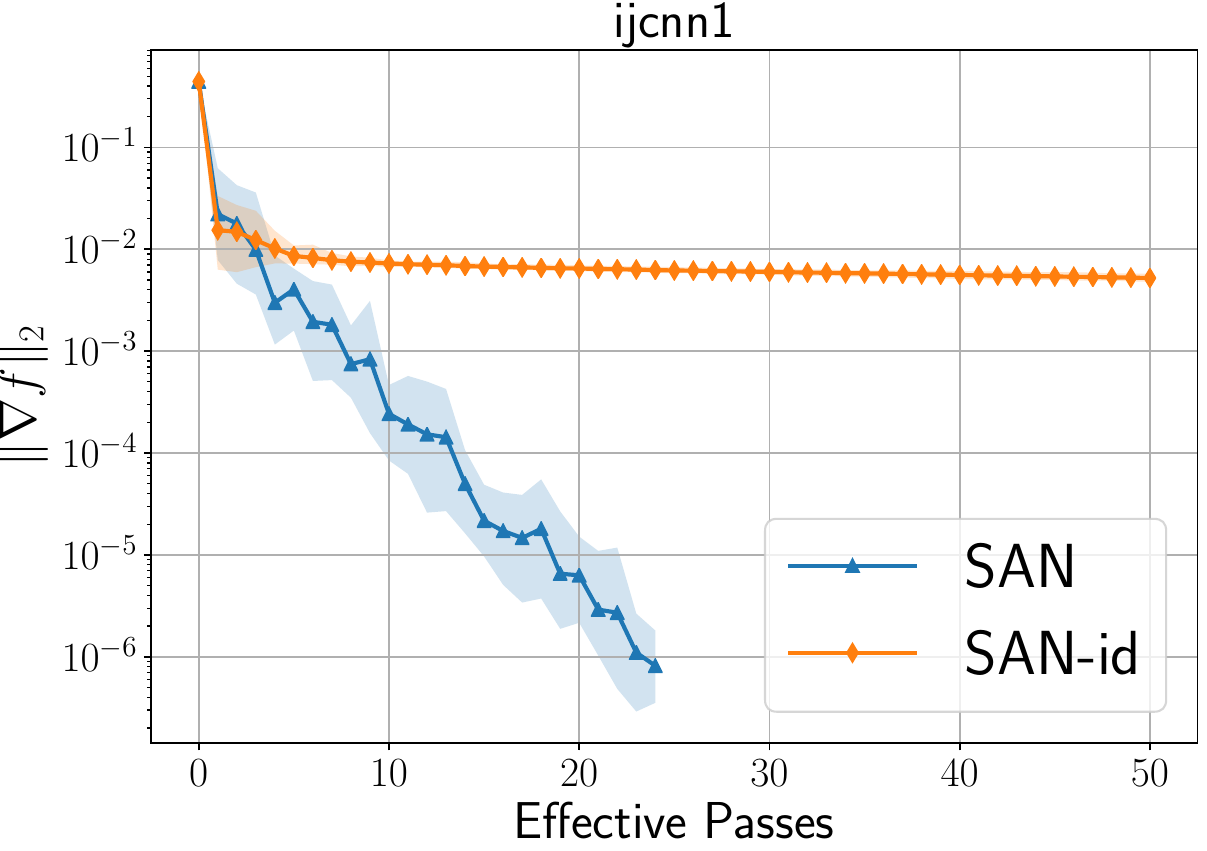}
\includegraphics[width=.24\linewidth]{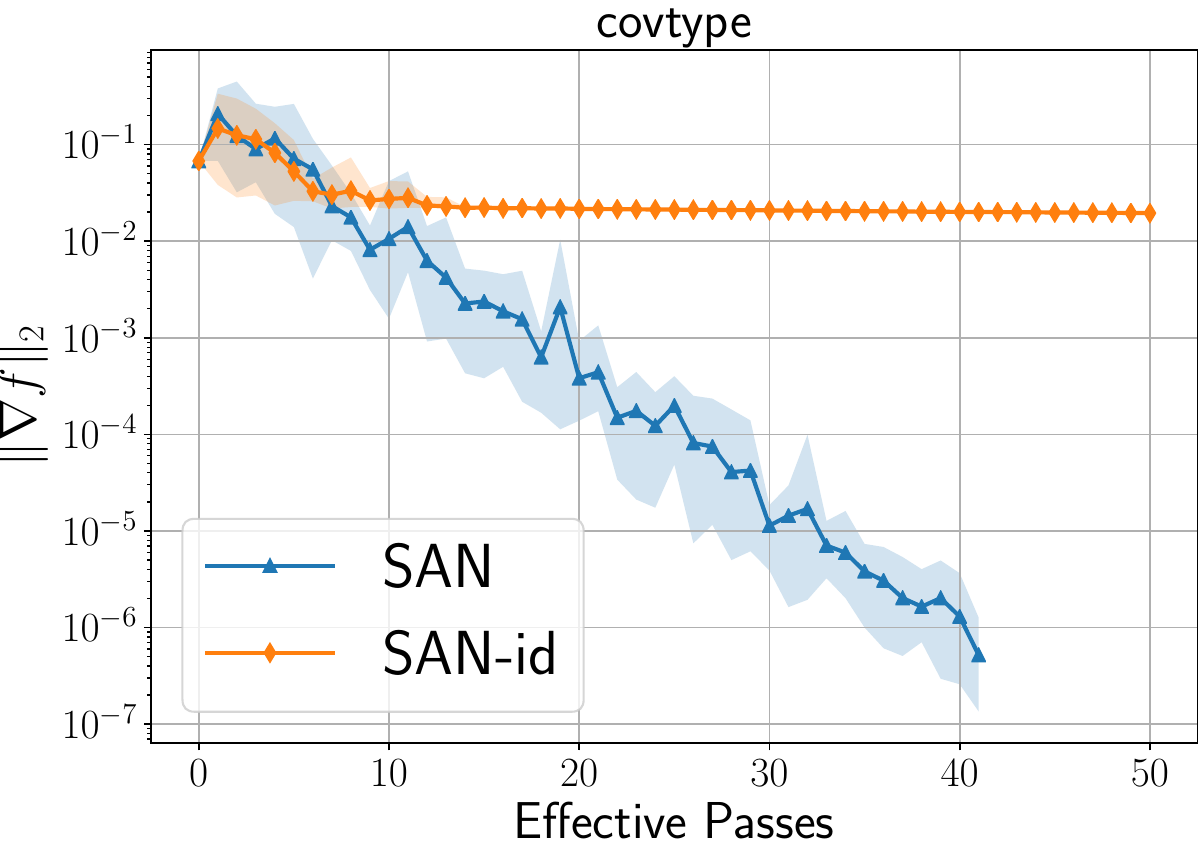}
\caption{L2-regularized logistic regression for SAN and SAN without the variable metric.}
\label{fig:san vs san-id}
\end{figure}

We call this algorithm \emph{SAN-id}\footnote{because this algorithm fits also in our SNRVM framework with $\mW_k \equiv \mI$ in~\eqref{algo:SNRVM}.} for short. However, we observe from Figure~\ref{fig:san vs san-id} that SAN-id only performs well at the early stage and stops converging to the optimum after the first few passes over the data. This motivated us to develop the version with the variable metric as introduced in the main text.

\section{SAN and SANA viewed as a sketched Newton Raphson method with variable metric } \label{sec:SNR viewpoint}

Here we provide a more detailed, step by step, introduction of the SAN and SANA methods. We also detail how SAN and SANA are particular instances of the Variable Metric Sketched Newton Raphson method introduced in the Section~\ref{sec:SNR}.

\subsection{A sketched Newton Raphson point of view}


Here we clarify  how the SAN and SANA methods are special cases of the sketched Newton Raphson method with a variable metric detailed in Section~\ref{sec:SNR}.

Let $x = \begin{bmatrix}
w \ ;
\alpha_1 \ ;
\cdots \ ;
\alpha_n 
 \end{bmatrix}\in \R^{(n+1)d}$ 
and $F:\R^{(n+1)d} \rightarrow \R^{(n+1)d}$ defined as
\begin{equation} 
F(x) \eqdef 
\begin{bmatrix}
\frac{1}{n}\sum\alpha_i \ ;
\nabla f_1(w) - \alpha_1 \ ;
\cdots \ ;
\nabla f_n(w) - \alpha_n
\end{bmatrix}.
\end{equation}
Therefore $w^* \in \mathbb{R}^d$ is a minimizer of~\eqref{eq:finite_sum} if and only if there 
exists $x^* = [w^*; \alpha_1^*; \dots; \alpha_n^*] \in \mathbb{R}^{(n+1)d}$ such that $F(x^*) = 0$. 
The Jacobian $\nabla F(x)$ is given by
\begin{eqnarray}
\nabla F(x) &=& \begin{bmatrix}
0 & \nabla^2f_1(w) & \cdots & \nabla^2f_n(w) \\
\frac{1}{n}\mI_d & & & \\
\vdots & & -\mI_{nd}& \\
\frac{1}{n}\mI_d & & & 
\end{bmatrix}
\in \R^{(n+1)d \times (n+1)d}, \label{eq: rsn_jac}
\end{eqnarray}
To find a zero of the function $F$, one could use the damped Newton Raphson method
\begin{equation}\label{D:SNR pseudoinverse}
x^{k+1} = x^k - \gamma \nabla F(x^k)^\top~^\dagger F(x^k), \quad \gamma \in (0,\,1].
\end{equation}
This can be equivalently rewritten 
 as a projection-and-relaxation step given by
\begin{equation}\label{D:SNR projection+relaxation gamma}
\begin{cases}
\bar x^{k+1} = {\rm{argmin}}~ \Vert x - x^k \Vert^2 \quad \text{ subject to } \quad \nabla F(x^k)^\top (x-x^k) = - F(x^k), \\
x^{k+1} = (1 - \gamma) x^k + \gamma \bar x^{k+1}.
\end{cases}
\end{equation}

Using the definition of our function $F$ in \eqref{eq: def_F}, we see that each iteration of the Newton Raphson method requires to project onto the following set of linear equations:
\begin{eqnarray}
& & 
\nabla F(x^k)^\top (x-x^k) = - F(x^k), \notag \\
& \Leftrightarrow & 
\begin{cases}
\frac{1}{n} \sum\limits_{i=1}^n (\alpha_i - \alpha_i^t) = - \frac{1}{n} \sum\limits_{i=1}^n \alpha_i^t, \notag \\
\nabla^2 f_i(w^t)(w - w^t) - (\alpha_i - \alpha_i^t) = \alpha_i^t- \nabla f_i(w^t) \text{ for } i \in \{1, \dots , n \},
\end{cases} \notag \\
& \Leftrightarrow &
\begin{cases}
\frac{1}{n} \sum\limits_{i=1}^n \alpha_i  = 0, \\
\nabla^2 f_i(w^t)(w - w^t) - \alpha_i = - \nabla f_i(w^t) \text{ for } i \in \{1, \dots , n \}.
\end{cases} \label{eq:NR linear eqt expanded}
\end{eqnarray}
Projecting onto \eqref{eq:NR linear eqt expanded} is challenging for two reasons: first it accesses all of the data (every function $f_i$ is involved) and second it requires solving a large linear system.

One approach to circumvent this bottleneck is to \textit{sketch} this linear system: at every iteration, instead of considering \eqref{eq:NR linear eqt expanded}, we will project onto a random row compression of this system.
Sketching can be for instance as simple as sampling one of the equations appearing in  \eqref{eq:NR linear eqt expanded}.
In its more general form, a sketch corresponds to any linear transformation of the equations.
In our context, this can be written as
\begin{equation}\label{eq:SNR linear eqt compact}
\mS^\top \nabla F(x^k)^\top (x-x^k) = - \mS^\top F(x^k),
\end{equation}
where $\mS \in \R^{(n+1)d \times \tau}$ is called the sketching matrix, and its number of columns $\tau$ is typically small.

This idea is at the core of the Sketched Newton Raphson method~\citep{Yuan2020sketched}, which aims at finding a zero of the function $F$ by iterating:
\begin{equation}\label{eq:SNR vanilla}
\begin{cases}
\bar x^{k+1} = {\rm{argmin}}~ \Vert x - x^k \Vert^2 \quad  \text{ subject to } \quad \mS_k^\top \nabla  F(x^k)^\top (x-x^k) = - \mS_k^\top F(x^k), \\
x^{k+1} = (1 - \gamma) x^k + \gamma \bar x^{k+1},
\end{cases}
\end{equation}
where $\mS_k$ is a sketching matrix randomly sampled at each iteration with respect to some distribution.

As we detailed in Section~\ref{sec:SNR}, the algorithms proposed in this paper can be seen as particular instances of a \textit{Variable Metric} Sketched Newton Raphson method.
This more general framework allows, at every iteration, to project the previous iterate onto  \eqref{eq:SNR linear eqt compact} with respect to some non-euclidean metric.
The algorithm writes as follows:
\begin{equation}\label{eq:SNR S&P projection+relaxation}
\begin{cases}
\bar x^{k+1} = {\rm{argmin}}~ \Vert x - x^k \Vert_{\mW_k}^2 \quad  \text{ subject to } \quad \mS_k^\top \nabla  F(x^k)^\top (x-x^k) = - \mS_k^\top F(x^k), \\
x^{k+1} = (1 - \gamma) x^k + \gamma \bar x^{k+1}.
\end{cases}
\end{equation}
Here, both $\mS_k$ and $\mW_k$ are randomly sampled with respect to a distribution which may depend on $x^k$. Besides, $\mW_k$ is positive-definite. The closed form solution to~\eqref{eq:SNR S&P projection+relaxation} is given in~\eqref{algo:SNRVM}, thanks to the following Lemma.

\begin{lemma}\label{L:SNRVM implicit to explicit}
If the iterates in \eqref{algo:SNRVM implicit} are well defined, then they are equivalent to \eqref{algo:SNRVM}.
\end{lemma}

\begin{proof}
Let $x^{k+1}$ be the iterate defined by \eqref{algo:SNRVM implicit}, where we assumed that the linear system $\mS_k^\top \nabla  F(x^k)^\top (x-x^k) = - \mS_k^\top F(x^k)$ has a solution.
Let us do a change of variable $u = x - x^k$, and write $\bar x^{k+1} = x^k + u^*$ where
\begin{equation*}
u^* = {\rm{argmin}}~ \Vert u \Vert_{\mW_k}^2 \quad  \text{ subject to } \quad \mS_k^\top \nabla  F(x^k)^\top u = - \mS_k^\top F(x^k).
\end{equation*}
We can call Lemma \ref{lem:Hproject} to obtain that 
\begin{equation*}
u^* = - \mW_k^{-1} \nabla F(x^k)\mS_k \left(\mS_k^\top  \nabla F(x^k)^\top \mW_k^{-1} \nabla F(x^k)\mS_k \right)^\dagger \mS_k^\top F(x^k).
\end{equation*}
The claim follows after writing that $x^{k+1} = (1- \gamma) x^k + \gamma (x^k + u^*) = x^k + \gamma u^*$.
\end{proof}



\subsection{SAN is a particular case of SNRVM}\label{sec:SNR viewpoint:SAN}

Let us consider SAN, described in Algorithm \ref{algo:SAN}, and rewrite it as an instance of the Variable Metric Sketched Newton Raphson method \eqref{eq:SNR S&P projection+relaxation}.
Given a probability $\pi \in (0,\,1)$, we define for all $x \in \mathbb{R}^{(n+1)d}$ a distribution $\mathcal{D}_x^{\mbox{\tiny SAN}}$ as follows: 
 $(\mS, \mW) \sim \mathcal{D}_x^{\mbox{\tiny SAN}}$ means that
\begin{itemize}
  \item with probability $\pi$ we have 
\begin{equation}\label{eq:rowsSWpo}
 \mS = 
\begin{bmatrix}
 \mI_d \\ \mo_d \\ \vdots \\ \mo_d
\end{bmatrix}
 \quad \mbox{and} \quad \mW = \mI_{(n+1)d},
\end{equation}
\item with probability  $1-\pi,$ we sample $j \in \{1, \cdots, n \}$ uniformly and set
\begin{equation}\label{eq:rowsSW1-po}
\mS = 
\begin{array}{cccc}
\begin{bmatrix}
\mo_d \\ \vdots \\ \mI_d \\ \vdots \\ \mo_d
\end{bmatrix}
&
\leftarrow j+1
&
\text{ and } 
&
\mW = \begin{bmatrix}
\nabla^2 f_j (w) & & & \\
      & \mI_d  & & \\
      & & \ddots & \\
      & & & \mI_d
\end{bmatrix}.
\end{array}
\end{equation}
\end{itemize}

\begin{lemma}\label{L:SAN from SNRVM}
Let $\pi \in (0,\,1 )$ and $\gamma \in (0,\,1]$ a step size.
Algorithm \ref{algo:SAN} (SAN) is equivalent to the Variable Metric Sketched Newton Raphson method \eqref{eq:SNR S&P projection+relaxation} applied to the function $F$ defined in \eqref{eq: def_F}, where at each iteration $(\mS_k, \mW_k)$ is sampled with respect to $\mathcal{D}_{x^k}^{\mbox{\tiny SAN}}$, as defined in \eqref{eq:rowsSWpo}-\eqref{eq:rowsSW1-po}.
\end{lemma}

\begin{proof}
Let us consider the Variable Metric Sketched Newton Raphson method  described in this Lemma.
We consider two cases, corresponding to the two classes of events described in \eqref{eq:rowsSWpo} and \eqref{eq:rowsSW1-po}.

Suppose that we are in the case (which holds with probability $\pi$)  given by~\eqref{eq:rowsSWpo}.
In this case we have
\begin{eqnarray*}
\mS_k^\top \nabla F(x^k)^\top (x - x^k) & \overset{\eqref{eq: rsn_jac}+\eqref{eq:rowsSWpo}}{=} &
\begin{bmatrix}
\mo_d & \frac{1}{n} \mI_d & \cdots & \frac{1}{n} \mI_d
\end{bmatrix}(x - x^k) \\
& = &
\frac{1}{n}\sum\limits_{i=1}^n (\alpha_i - \alpha_i^k) \\
\mS_k^\top F(x^k) &\overset{\eqref{eq: def_F}+\eqref{eq:rowsSWpo} }{ =}& \frac{1}{n}\sum\limits_{i=1}^n \alpha_i^k.
\end{eqnarray*}
Those expressions mean that the linearized equation  \eqref{eq:SNR linear eqt compact} is equivalent to $\frac{1}{n}\sum\limits_{i=1}^n \alpha_i = 0$.
So the update of the variables is exactly given by \eqref{e:sanie1}.

Let now $j$ be in $\{1, \dots, n\}$ sampled uniformly, and suppose that we are in the case given by~\eqref{eq:rowsSW1-po}.
We can then compute
\begin{eqnarray*}
\mS_k^\top \nabla F(x^k)^\top(x-x^k)
&  \overset{\eqref{eq: rsn_jac}+\eqref{eq:rowsSW1-po}}{=} &
\begin{bmatrix}
\nabla^2 f_j(w^k) & \mo_d & \cdots & -\mI_d& \cdots & \mo_d
\end{bmatrix}(x-x^k) \\
&=&
\nabla^2 f_j(w^k)(w - w^k) - (\alpha_j - \alpha_j^k) \\
\mS_k^\top  F(x^k)
&  \overset{\eqref{eq: def_F}+\eqref{eq:rowsSW1-po}}{=} &
\nabla f_j(w^k) - \alpha_j^k.
\end{eqnarray*}
Those expressions mean that the linearized equation  \eqref{eq:SNR linear eqt compact} is equivalent to $\nabla^2 f_j(w^k)(w - w^k) - \alpha_j = - \nabla f_j(w^k)$.
So the update of the variables is exactly given by \eqref{e:sanie2}.
The conclusion follows Lemma \ref{L:SAN implicit to explicit}.
\end{proof}

\subsection{SANA is a particular case of SNRVM}\label{sec:SNR viewpoint:SANA}

Let us consider SANA, described in Algorithm \ref{algo:SANA}, and rewrite it as an instance of the Variable Metric Sketched Newton Raphson method \eqref{eq:SNR S&P projection+relaxation}.
We define for all $x \in \mathbb{R}^{(n+1)d}$ a distribution $\mathcal{D}_x^{\mbox{\tiny SANA}}$ as follows: 
 $(\mS, \mW) \sim \mathcal{D}_x^{\mbox{\tiny SANA}}$ means that, with probability $1/n$ we sample $j \in \{1, \cdots, n \}$ and we have
\begin{equation}\label{eq:SANA distrib}
\mS = 
\begin{array}{cccc}
\begin{bmatrix}
\mI_d & \mo_d \\ 
\mo_d & \vdots \\ 
\vdots & \mI_d \\ 
\vdots & \vdots \\
\mo_d &  \mo_d
\end{bmatrix}
&
\leftarrow j+1
&
\text{ and } 
&
\mW = \begin{bmatrix}
\nabla^2 f_j (w) & & & \\
      & \mI_d  & & \\
      & & \ddots & \\
      & & & \mI_d
\end{bmatrix}.
\end{array}
\end{equation}

\begin{lemma}\label{L:SANA equivalence forms}
Let $\gamma \in (0,\,1]$ be a step size.
Algorithm \ref{algo:SANA} (SANA) is equivalent to the Variable Metric Sketched Newton Raphson method \eqref{eq:SNR S&P projection+relaxation} applied to the function $F$ defined in \eqref{eq: def_F}, where at each iteration $(\mS_k, \mW_k)$ is sampled with respect to $\mathcal{D}_{x^k}^{\mbox{\tiny SANA}}$, as defined in \eqref{eq:SANA distrib}.
\end{lemma}

\begin{proof}
Let $k \in \mathbb{N}$, and suppose that we have sampled $j\in \{1,\ldots, n\}$ and $\mS_k$ and $\mW_k$ according to~\eqref{eq:SANA distrib}.
Therefore,
\begin{eqnarray*}
\mS_k^\top \nabla F(x^k)^\top (x-x^k)
&  \overset{\eqref{eq: rsn_jac}+\eqref{eq:SANA distrib}}{=} &
\begin{bmatrix}
\mo_d & \frac{1}{n} \mI_d  & \cdots & \cdots & \cdots & \frac{1}{n} \mI_d  \\
\nabla^2 f_j(w^k) & \mo_d &\cdots & - \mI_d & \cdots & \mo_d
\end{bmatrix}(x-x^k) \\
&=&
\begin{bmatrix}
\frac{1}{n}\sum\limits_{i=1}^n (\alpha_i - \alpha_i^k) \\
\nabla^2 f_j(w^k) (w-w^k) - (\alpha_j - \alpha_j^k)
\end{bmatrix} \\
\mS_k^\top  F(x^k)
&  \overset{\eqref{eq: def_F}+\eqref{eq:SANA distrib}}{=} &
\begin{bmatrix}
\frac{1}{n}\sum\limits_{i=1}^n \alpha_i^k \\
\nabla f_j(w^k) - \alpha_j^k
\end{bmatrix}.
\end{eqnarray*}
Those expressions mean that the linearized equation  \eqref{eq:SNR linear eqt compact} is equivalent to the two equations $\sum\limits_{i=1}^n \alpha_i=0$ and $\nabla^2 f_j(w^k)(w - w^k) - \alpha_j = - \nabla f_j(w^k)$.
So the update of the variables is exactly given by \eqref{e:sanaie}.
The conclusion follows Lemma \ref{L:SANA implicit to explicit}.

\end{proof}

\section{Proofs for the results in Section \ref{sec:SNR}, including Theorems \ref{them: AdapSNR} and \ref{T:CV SAN explicit} }\label{sec:annex sec 4}

\subsection{Proof of Proposition \ref{P:Ass verified for SAN and SANA}}\label{SS:proof of Ass verified for SAN and SANA}

\begin{proof}
Let us start by showing that Assumption \ref{Ass:distribution D_x} is satisfied for SAN and SANA.
The distribution $\mathcal{D}_x^{\mbox{\tiny SAN}}$ (resp. $\mathcal{D}_x^{\mbox{\tiny SANA}}$) defined in the Section \ref{sec:SNR viewpoint:SAN} (resp. Section \ref{sec:SNR viewpoint:SANA}) is clearly finite and proper so long as $\pi \in (0,1)$.
It remains to compute $\mathbb{E}[\mS\mS^\top]$.
We can see that it is a block-diagonal matrix $\Diag{\pi   \mI_d,\frac{1-\pi}{n} \mI_d,  \cdots, \frac{1-\pi}{n} \mI_d}$ (resp. $\Diag{\mI_d, \frac{1}{n} \mI_d, \cdots, \frac{1}{n} \mI_d}$), which is invertible since  $\pi \in (0,1)$.

Now let us turn on Assumption \ref{ass:existence}.
To prove 
that $\nabla F(x)^\top \nabla F(x)$ is invertible, it is enough to show that $\nabla F(x)$ is injective.
Let $x =  ( w \ ; \alpha) \in \R^{d + dn}$, and let us first show that $\nabla F(x)$ is injective. 
Suppose there exists $\bar x = (\bar w \ ; \bar \alpha) \in \R^{d + dn}$ such that $\nabla F(x) \bar x =0.$ Consequently from~\eqref{eq: rsn_jac}  we have that 
\begin{align}
\sum_{i=1}^n \nabla^2 f_i(w) \bar \alpha_i &=0 \nonumber \\
\frac{1}{n} \bar w & = \bar \alpha_i, \quad \mbox{for }i=1,\ldots, n. \label{eq:tersmrojme9ps}
\end{align}
Substituting out the $\bar \alpha_i$'s we have that 
\[\nabla^2 f(w)\bar w = \frac{1}{n} \sum_{i=1}^n \nabla^2 f_i(w) \bar w =0.  \]
Consequently, since $\nabla^2 f(w)$ is positive definite (recall Assumption \ref{Ass:strict convexity of problem}), and in particular injective, we have that
$\bar w =0$.
Thus it follows from~\eqref{eq:tersmrojme9ps}  that $\bar \alpha_i =0$ for $i=1,\ldots, n$. 
This all shows that $\bar x =0$, and concludes the proof that $\nabla F(x) $ is injective.

Furthermore, $\nabla F(x)$ is a square matrix, thus invertible. We have $F(x) \; \in \; \Image{\nabla F(x)^\top}$.

Finally $\nabla F(x)^\top  \nabla F(x)$ is invertible since
\[\Null{\nabla F(x)^\top  \nabla F(x)} \; = \; \Null{\nabla F(x) } \; = \; \{0\}.\]

\end{proof}



\subsection{SNRVM is equivalent to minimizing a quadratic function over a random subspace}

\begin{lemma}\label{lem: lemma10_RSN}
(Lemma~$10$ in~\citet{RSN_nips}). 
For every matrix $\mM$ and symmetric positive semi-definite matrix $\mG$ such that $\Null{\mG} \subset 
\Null{\mM}$, 
we have that $\Null{\mM^\top} = \Null{\mM \mG \mM^\top}$.
\end{lemma}

\begin{lemma}\label{L:SNRVM as minimize quadratic over random subspace}
Let Assumptions~\ref{Ass:distribution D_x} and \ref{ass:existence} hold. 
Then the iterates of SNRVM are equivalent to
\begin{eqnarray}\label{algo:SNRVM newton random subspace}
x^{k+1} &=& \argmin_{x \in \R^p} \ \hat{f}_k(x^k) + \dotprod{\nabla \hat{f}_k(x^k) , x-x^k} + \frac{1}{2\gamma}\norm{x-x^k}^2_{\mW_k}  \\
& & \mbox{subject to } x \in x^k +\Image{ \mW_k^{-1} \nabla F(x^k) \mS_k},\nonumber
\end{eqnarray}
where $\hat{f}_k$ is defined in~\eqref{eq: def_hat_f}.
\end{lemma}
\begin{proof}
Start by observing that the problem in \eqref{algo:SNRVM newton random subspace} is strongly convex, and  therefore has a unique solution that we will note $x^*$. Let us prove that $x^*$ is exactly $x^{k+1}$ whose closed form expression is given in \eqref{algo:SNRVM}.
For this, let $\tau$ be the number of columns for $\mS_k$, and let $u \in \mathbb{R}^\tau$.
We can then write that $x^* = x^k + \mW_k^{-1} \nabla F(x^k) \mS_k u^*$, where $u^*$ is any solution of the following unconstrained optimization problem:
\begin{eqnarray*}
u^*  & \in & \argmin_{u \in \mathbb{R}^\tau} \dotprod{\nabla \hat{f}_k(x^k) , \mW_k^{-1} \nabla F(x^k) \mS_k u^*} + \frac{1}{2\gamma}\norm{\mW_k^{-1} \nabla F(x^k) \mS_k u^*}^2_{\mW_k} .
\end{eqnarray*}
Writing down the optimality conditions for this convex quadratic problem, we see that $u^*$ must verify:
\begin{equation*}
\gamma \mS_k^\top \nabla F(x^k)^\top \mW_k^{-1}  \nabla \hat{f}_k(x^k) + \mS_k^\top \nabla F(x^k)^\top \mW_k^{-1}\nabla F(x^k)\mS_k u^*.
\end{equation*}
Let us choose the pseudo inverse solution of this linear system:
\begin{equation}\label{e:snrvm rns0}
u^* = - \gamma \left( \mS_k^\top \nabla F(x^k)^\top \mW_k^{-1}\nabla F(x^k)\mS_k \right)^\dagger \mS_k^\top \nabla F(x^k)^\top \mW_k^{-1}  \nabla \hat{f}_k(x^k).
\end{equation}
Using the definition of $\hat f_k$, we can write
\[\nabla \hat{f}_k(x^k) = \nabla F(x^k) \left( \nabla F(x^k)^\top \mW_k^{-1} \nabla F(x^k) \right)^\dagger F(x^k).\]
All we need to prove now is that
\begin{equation}\label{e:snrvm rns1}
\nabla F(x^k)^\top \mW_k^{-1} \nabla F(x^k) \left( \nabla F(x^k)^\top \mW_k^{-1} \nabla F(x^k) \right)^\dagger F(x^k) = F(x^k).
\end{equation}
To see why \eqref{e:snrvm rns1} is true, first notice that $\nabla F(x^k)^\top \mW_k^{-1} \nabla F(x^k) \left( \nabla F(x^k)^\top \mW^{-1} \nabla F(x^k) \right)^\dagger$ is the orthogonal projector onto the range of $\nabla F(x^k)^\top \mW_k^{-1} \nabla F(x^k)$.
Moreover, using the fact that $\mW_k^{-1}$ is injective together with Lemma \ref{lem: lemma10_RSN}, we can write that
\begin{eqnarray*}
\Image{\nabla F(x^k)^\top \mW_k^{-1} \nabla F(x^k)} &=& (\Null{F(x^k)^\top \mW_k^{-1} \nabla F(x^k)}^\bot \\
&=&
(\Null{\nabla F(x^k)}^\bot \\
&=&
\Image{\nabla F(x^k)^\top}.
\end{eqnarray*}
Since we know from Assumption \ref{ass:existence} that $\nabla F(x^k)^\top$ is surjective, and so that $F(x^k)$ belongs in the range of $\nabla F(x^k)^\top$, we deduce that \eqref{e:snrvm rns1} is true.
We can now inject \eqref{e:snrvm rns1} into \eqref{e:snrvm rns0}, and obtain finally that
\begin{eqnarray*}
x^* & = & x^k +  \mW_k^{-1} \nabla F(x^k) \mS_k u^* \\
 & = & x^k - \gamma  \mW_k^{-1} \nabla F(x^k) \mS_k
 \left( \mS_k^\top \nabla F(x^k)^\top \mW_k^{-1}\nabla F(x^k)\mS_k \right)^\dagger \mS_k^\top
 F(x^k),
\end{eqnarray*}
which is exactly \eqref{algo:SNRVM}.
\end{proof}

\subsection{About \texorpdfstring{$\rho$}{r} in Theorem~\ref{them: AdapSNR} }

\begin{lemma}\label{L:null space sum}
If $A$, $B$ are two symmetric positive semi-definite matrices
then $\Null{A+B}~=~\Null{A} \cap \Null{B}$.
\end{lemma}

\begin{proof}
If $x \in \Null{A} \cap \Null{B}$ then it is trivial to see that $x \in \Null{A+B}$.
If $x \in \Null{A+B}$, then 
\begin{equation*}
0 = \langle (A+B)x,x \rangle = \langle Ax,x \rangle + \langle Bx,x \rangle,
\end{equation*}
where by positive semi-definiteness we have $\langle Ax,x \rangle \geq 0$ and $\langle Bx,x \rangle \geq 0$.
The sum of nonnegative numbers being nonegative, we deduce that 
\begin{equation*}
\langle Ax,x \rangle
=
\langle Bx,x \rangle
=
0.
\end{equation*}
Since $\langle Ax,x \rangle=0$, we deduce from the fact that $A$ is symmetric that $Ax=0$.
Similarly, $Bx=0$, which concludes the proof.
\end{proof}

The following Lemma will be needed in the proof of Theorem~\ref{them: AdapSNR}. 

\begin{lemma}\label{lem: sdjffskapp}
Recall the definition of $\mH(x)$ given by
\begin{align}\label{eq:H matrix}
\mH(x) \eqdef \mathbb{E}\left[\mS\left(\mS^\top \nabla F(x)^\top \mW^{-1} \nabla F(x)\mS\right)^\dagger\mS^\top\right],
\end{align}
If Assumption \ref{Ass:distribution D_x} and \ref{ass:existence}  hold, then $\mH(x)$  is invertible. Moreover, for every symmetric positive definite matrix $\mW$ and $x \in \mathbb{R}^p$ we have that 
\begin{equation}\label{eq:nasty smallest eigenvalue positive}
      \min_{v \in \Image{\mW^{-1/2}\nabla F(x)}\backslash\{0\}} \frac{\dotprod{\mW^{-1/2}\nabla F(x)\mH(x)\nabla F(x)^\top\mW^{-1/2}  v , v}}{\norm{v}^2}
\end{equation}
is exactly the smallest positive eigenvalue of $\mW^{-1/2}\nabla F(x) \mH(x) \nabla F(x)^\top \mW^{-1/2}$.
\end{lemma}
\begin{proof}
Let $x \in \R^{p}$ and $(\mS, \mW) \sim \cD_x$. Let $\mG = \nabla F(x)^{\top} \mW^{-1} \nabla F(x)$ which is symmetric positive semi-definite.  Since   $\nabla F(x)^\top \nabla F(x)$ and $\mW$ are invertible we have that $\mG$ is invertible. Consequently  $\Null{\mG}=\{0\} \subset \Null{\mS^\top}$. Thus
by Lemma~\ref{lem: lemma10_RSN}  (with  $\mM = \mS^{\top}$) we have that
\[
\Null{\left(\mS^{\top} \nabla F(x)^{\top} \mW^{-1} \nabla F(x) \mS \right)^{\dagger}}= 
\Null{\mS^{\top} \nabla F(x)^{\top} \mW^{-1} \nabla F(x) \mS } = \Null{\mS}.
\]
Using Lemma~\ref{lem: lemma10_RSN} once again with $\mG = \left(\mS^{\top} \nabla F(x)^{\top} \mW^{-1} \nabla F(x) \mS \right)^{\dagger}$ and $\mM = \mS$, we have that
\begin{equation} \label{eq: ker()=ker(SS^top)}
{\bf Null}\Big(\underbrace{\mS\left(\mS^{\top} \nabla F(x)^{\top} \mW^{-1} \nabla F(x) \mS \right)^{\dagger} \mS^{\top}}_{\eqdef \mH_{\mS, \mW}(x)}\Big) = 
\Null{\mS^{\top}} = \Null{\mS\mS^{\top}}.
\end{equation}
Observe that with our notations and from Assumption \ref{Ass:distribution D_x}, 
\begin{equation*}
\mH(x) = \EE{\mS, \mW \sim \cD_x}{\mH_{\mS, \mW}(x)}
=
\sum\limits_{i=1}^r p_i \mH_{\mS_i(x), \mW_i(x)}(x).
\end{equation*}
As $\mH_{\mS_i(x), \mW_i(x)}(x)$ is symmetric positive semi-definite, we can use Lemma \ref{L:null space sum} to write
\begin{eqnarray}
\Null{\mH(x) } = \Null{ \sum\limits_{i=1}^r p_i \mH_{\mS_i(x), \mW_i(x)}(x) } & = & \bigcap\limits_{i=1}^r
\Null{\mH_{\mS_i(x), \mW_i(x)}(x) } \nonumber \\
& \overset{\eqref{eq: ker()=ker(SS^top)}}{=} &  \bigcap\limits_{i=1}^r \Null{\mS_i(x)\mS_i(x)^{\top}}  \nonumber \\
& = & \Null{\EE{\mS \sim \cD_{x}}{\mS \mS^{\top}}} = \{0\} \nonumber
\end{eqnarray}

This means that $\mH(x)$ is invertible for all $x \in \R^p$.

Now,  take any $x \in \R^p$, and a symmetric positive definite matrix $\mW$. 
Then $\mW^{-1/2}\nabla F(x) \mH(x) \nabla F(x)^\top \mW^{-1/2}$ is symmetric, semi-definite positive.
 Since $\mH(x)$ and $\mW$ are invertible, we can apply
Lemma~\ref{lem: lemma10_RSN} again to obtain
\begin{eqnarray} \label{eq:Ker(DF^top) = Ker()}
\Null{\nabla F(x)^\top \mW^{-1/2}} &=& \Null{\mW^{-1/2} \nabla F(x) \mH(x) \nabla F(x)^\top \mW^{-1/2} }.
\end{eqnarray}
Consequently
\begin{align} \label{eq: range()=ker()^perp}
\Image{\mW^{-1/2} \nabla F(x)}& =  \left(\Null{\nabla F(x)^\top \mW^{-1/2}} \right)^{\perp}  \nonumber \\
&
\overset{\eqref{eq:Ker(DF^top) = Ker()}}{=} \left(\Null{ \mW^{-1/2} \nabla F(x) \mH(x) \nabla F(x)^\top \mW^{-1/2}} \right)^{\perp}.
\end{align}

Therefore, we conclude that \eqref{eq:nasty smallest eigenvalue positive} is equal to
\begin{align}
 &  \min_{v \in \left(\Null{\mW^{-1/2} \nabla F(x) \mH(x) \nabla F(x)^\top  \mW^{-1/2}} \right)^{\perp} \backslash\{0\}} \frac{\dotprod{\mW^{-1/2}\nabla F(x)\mH(x)\nabla F(x)^\top \mW^{-1/2}v , v}}{\norm{v}^2} \nonumber \\
& \qquad \qquad= \lambda^{+}_{\min} \left(\mW^{-1/2}\nabla F(x) \mH(x) \nabla F(x)^\top \mW^{-1/2} \right) > 0 \nonumber.
\end{align}
\end{proof}

\if{ 
\subsection{Proof of Proposition \ref{P:rhopositive SAN}}


\begin{proof}
Let $x \in \mathbb{R}^{p}$ 
and $i \in \{1, \dots, r \}$.
By definition, 
\begin{equation*}
\rho_i(x) \eqdef \lambda^{+}_{\min} \left(\mW_{i}(x)^{-1/2}\nabla F(x) \mH(x) \nabla F(x)^\top \mW_{i}(x)^{-1/2} \right)
\end{equation*}
is positive.
Moreover, Assumptions \ref{Ass:distribution D_x} and \ref{ass:existence}, together with Lemma \ref{lem: sdjffskapp}, guarantee that the matrix $\mW_{i}(x)^{-1/2}\nabla F(x) \mH(x) \nabla F(x)^\top \mW_{i}(x)^{-1/2}$ is invertible.
So we actually have that 
\begin{equation*}
\rho_i(x) = \lambda_{\min} \left(\mW_{i}(x)^{-1/2}\nabla F(x) \mH(x) \nabla F(x)^\top \mW_{i}(x)^{-1/2} \right),
\end{equation*}
where $\lambda_{\min}$ denotes the smallest eigenvalue.
This means that $\rho_i(x)$ depends continuously on $x$, provided that $\mW_{i}(x)^{-1/2}\nabla F(x) \mH(x) \nabla F(x)^\top \mW_{i}(x)^{-1/2}$ depends continuously on $x$. This is what we are going to show next. 

For SAN and SANA, we have that $F$ is of class $\cC^1$, and so is $\mW_i(x)$ 
(see Sections \ref{sec:SNR viewpoint:SAN} and \ref{sec:SNR viewpoint:SANA}).
It remains to verify that $\mH(x)$ is continuous.
From its definition \eqref{eq:H matrix}, we see that it is enough to show that $\mS^\top \nabla F(x)^\top \mW^{-1} \nabla F(x)\mS$ is invertible, where $(\mS,\mW)$ is drawn with respect to $\mathcal{D}_x^{\mbox{\tiny SAN}}$ (resp. $\mathcal{D}_x^{\mbox{\tiny SANA}}$).
And this is true due to Lemma \ref{lem: lemma10_RSN}, using  the fact that $\nabla F(x)$ is injective (Proposition \ref{P:Ass verified for SAN and SANA}), that $\mW$ is invertible, and that $\mS$ is injective (see Sections \ref{sec:SNR viewpoint:SAN} and \ref{sec:SNR viewpoint:SANA}).

We have proved that $\rho_i$ is continuous and positive, so we deduce that $\rho(x) = \min\limits_{i=1, \dots, s} \rho_i(x)$ is also continuous and positive.
Therefore, over any compact set $\rho$ is bounded away from zero.
The conclusion follows directly from our boundedness assumption on the sequence of iterates.
\end{proof}

}\fi 

\subsection{Proof of Theorem~\ref{them: AdapSNR} }

\label{asec:proofTheoremAdapSNR}

\begin{proof}

Let $k \in \mathbb{N}$.
In this proof, we will write $\nabla F_k$ as a shorthand for $\nabla F(x^k)$, and we introduce the notation $ \nabla^{\mW} F_k \eqdef  \mW_k^{-1/2} \nabla F(x^k)$.
First we aim to establish a relationship between  $\hat{f}_k(x^k) = \norm{F(x^k)}^2_{(\nabla F_k^\top \mW_k^{-1} \nabla F_k)^{\dagger}}$ and $\norm{F(x^k)}^2_{\mH(x^k)}$. 
Observe that Assumption \ref{ass:existence} allows us to write that
\begin{equation}\label{eq:DFxDF^dagxF=F}
(\forall x \in \mathbb{R}^p) \quad
F(x) = \nabla F(x)^\top  \mW_k^{-1/2}(\nabla F(x)^\top  \mW_k^{-1/2})^{\dagger} F(x).
\end{equation}
This is due to the fact that
$ \nabla F(x)^\top  \mW_k^{-1/2}(\nabla F(x)^\top  \mW_k^{-1/2})^{\dagger} $ is the projection matrix onto $\Image{\nabla F(x)^\top  \mW_k^{-1/2} }$, where $\mW_k^{-1/2}$ is surjective, meaning that $\Image{\nabla F(x)^\top }  = \Image{\nabla F(x)^\top  \mW_k^{-1/2} }$.
From~\eqref{eq:DFxDF^dagxF=F} we have that
\begin{eqnarray}
\norm{F(x^k)}^2_{\mH(x^k)} & = & \dotprod{F(x^k), \mH(x^k) F(x^k)} \nonumber \\
&  \overset{\eqref{eq:DFxDF^dagxF=F}}{=} & \dotprod{\nabla^{\mW} F_k ^\top(\nabla^{\mW} F_k^\top )^{\dagger} F(x^k), \mH(x^k)\nabla^{\mW} F_k^\top (\nabla^{\mW} F_k^\top)^{\dagger} F(x^k)} \nonumber \\
&= & \dotprod{(\nabla^{\mW} F_k^\top)^{\dagger} F(x^k), \nabla^{\mW} F_k \mH(x^k)\nabla^{\mW} F_k^\top \left( (\nabla^{\mW} F_k^\top)^{\dagger} F(x^k)\right)} \nonumber \\
 & \overset{\text{Lemma \ref{lem: sdjffskapp}}} {\geq} & \rho \ \Vert (\nabla^{\mW} F_k^\top)^{\dagger} F(x^k) \Vert^2 \label{eq:tempzmoemze}\\
& = & \rho \norm{F(x^k)}_{\left(\nabla^{\mW} F_k^\top \nabla^{\mW} F_k\right)^\dagger}^2   \label{eq: norm_H_jac-temp2} \\
 & \overset{\eqref{eq: def_hat_f}}{ =}  &2 \rho \hat{f}_k(x^k),  \label{eq: norm_H_jac2}
\end{eqnarray}
where in~\eqref{eq:tempzmoemze} we used that $\Image{(\nabla^{\mW} F_k^\top)^{\dagger}}  =\Image{\nabla^{\mW} F_k}$ together with Lemma \ref{lem: sdjffskapp}, and in~\eqref{eq: norm_H_jac-temp2} we used that  $(\mM)^{\dagger\top}(\mM)^{\dagger} = (\mM^\top)^{\dagger}(\mM)^{\dagger}  =(\mM \mM^\top)^\dagger$ for every matrix $\mM.$
Now, we turn onto the study the term $\hat{f}_k(x^k)$.
Compute its gradient with respect to the metric induced by $\mW_k$ at $x^k$:
\begin{eqnarray} \label{eq:gradfhatisNewton}
\nabla^{\mW_k} \hat{f}_k(x^k) & =&  \mW_k^{-1}\nabla F_k(\nabla F_k^\top  \mW_k^{-1} \nabla F_k)^\dagger F(x^k) \nonumber \\
& =& \mW_k^{-1/2}(\nabla F_k^\top  \mW_k^{-1/2})^{\dagger} F(x^k).
\end{eqnarray}
This, together with the Assumption~\ref{ass:AdapUpperBnd}, allows us to write that
\begin{eqnarray}
\hat f_{k+1}(x^{k+1}) & \leq &    \hat{f}_k(x^k) + \dotprod{\nabla \hat{f}_k(x^k) , x^{k+1}-x^k} + \frac{L}{2}\norm{x^{k+1}-x^k}^2_{\mW_k} \nonumber \\
& \overset{\eqref{algo:SNRVM}+\eqref{eq:gradfhatisNewton}} {=}& \hat{f}_k(x^k)   -\gamma\dotprod{ n^{\mW_k}(x^k)  , \mW_k^{-1} \nabla F_k\mS_k\left(\mS_k^\top \nabla F_k^\top \mW_k^{-1} \nabla F_k\mS_k\right)^\dagger\mS_k^\top F(x^k)}_{\mW_k}  \nonumber \\ & & \quad + \frac{\gamma^2 L}{2}\norm{\mW_k^{-1} \nabla F_k\mS_k\left(\mS_k^\top \nabla F_k^\top \mW_k^{-1} \nabla F_k\mS_k\right)^\dagger\mS_k^\top F(x^k)}^2_{\mW_k} \nonumber \\
& \overset{\eqref{eq:DFxDF^dagxF=F}}{=} & \hat{f}_k(x^k)   - \gamma  \dotprod{ F(x^k)  , \mS_k\left(\mS_k^\top \nabla F_k^\top \mW_k^{-1} \nabla F_k\mS_k\right)^\dagger\mS_k^\top F(x^k)}  \nonumber \\ & & \quad + \frac{\gamma^2L}{2}\norm{\mW_k^{-1} \nabla F_k\mS_k\left(\mS_k^\top \nabla F_k^\top \mW_k^{-1} \nabla F_k\mS_k\right)^\dagger \mS_k^\top F(x^k)}^2_{\mW_k} \nonumber\\
& = & \hat{f}_k(x^k)   -\gamma\left(1- \frac{\gamma L}{2} \right) \norm{ F(x^k) }^2_{\mS_k\left(\mS_k^\top \nabla F_k^\top \mW_k^{-1} \nabla F_k\mS_k\right)^\dagger\mS_k^\top} \nonumber \\
& \overset{\gamma = 1 / L}{=} & \hat{f}_k(x^k)   - \frac{\gamma}{2}  \norm{ F(x^k) }^2_{\mS_k\left(\mS_k^\top \nabla F_k^\top \mW_k^{-1} \nabla F_k\mS_k\right)^\dagger\mS_k^\top} \label{eq: adap_hat_f_k_decreasing} 
\end{eqnarray}
where in~\eqref{eq: adap_hat_f_k_decreasing} we use the identity $\mM^\dagger \mM \mM^\dagger = \mM^\dagger$ with $\mM = \mS_k^\top \nabla F_k^\top \mW_k^{-1} \nabla F_k\mS_k$.
Taking the expectation conditioned on $x_k$ in the inequality~\eqref{eq: adap_hat_f_k_decreasing} gives
\begin{eqnarray}
\EE{}{\hat{f}_{k+1}(x^{k+1}) \; \vert \; x^k} & \leq & \EE{}{\hat{f}_{k}(x^{k}) \; \vert \; x^k}  -  \frac{\gamma}{2}  \norm{ F(x^k) }^2_{\mH(x^k)} \nonumber\\
& \leq & \EE{}{\hat{f}_{k}(x^{k}) \; \vert \; x^k}  - \rho \gamma  \hat{f}_{k}(x^{k}). \nonumber
\end{eqnarray}
Taking full expectation and expanding the recurrence gives finally
\begin{equation*}
\EE{}{\hat{f}_{k+1}(x^{k+1})} \leq (1- \rho \gamma) \EE{}{\hat{f}_{k}(x^{k})}.
\end{equation*}
\end{proof}

\subsection{SNRVM for solving linear systems}
\label{asec:linearexample}

Here we consider the simplified case in which  our objective function~\eqref{eq:finite_sum} is a quadratic function.
In this case, the stationarity condition \eqref{eq:weqai} is a linear system.
To simplify the notation, let us denote in this section the resulting linear system as 
\begin{equation} \label{eqn:problin}
\mA x = b, \quad \mbox{where } \mA \in \mathcal{M}_{p}(\mathbb{R}), \ p=d(n+1).
\end{equation}
In other words, our nonlinear map is given by $F(x) =\mA x - b$.
Because $\nabla F(x) = \mA^\top$, where $\mA$ is a square matrix, we see that the assumptions on $F$ in Assumption \ref{ass:existence} are verified if and only if $\mA$ is invertible.
%
%
%

In this setting the SNR method~\eqref{eq:linearFxsketchedproj} is known as the sketch-and-project method~\citep{Gower2015}. The sketch-and-project method has been shown to converge linearly at a fast rate~\citep{Gower2015,Gower2015c}. 
Thus this quadratic case serves as a good sanity check to verify if our rate of convergence in Theorem~\ref{them: AdapSNR} recovers the well known fast linear rate of  the sketch-and-project method. This is precisely what we investigate in the next lemma.
It only remains to reformulate Assumption~\ref{ass:AdapUpperBnd}, which we do in the following lemma.

\begin{lemma}\label{L:nasty ass for linear case}If $\{\mW_{k}\}_{k \in \mathbb{N}}$ is a sequence of invertible matrices such that
\begin{equation}\label{eq:metricdescentquad}
 \mW_{k+1} \preceq  \mW_k,
\end{equation}
and  $\mA$ is invertible, 
then Assumption \ref{ass:AdapUpperBnd} holds with $L=1$.
\end{lemma}

\begin{proof}
Using the fact that $F(x) = \mA x -b$ and $\nabla F(x) = \mA^\top$, we can rewrite definition \eqref{eq: def_hat_f} as
 \begin{equation}\label{eq:fhatlinear}
   \hat{f}_k(x)  \; = \; \frac{1}{2}\norm{\mA x- b}_{(\mA \mW_k^{-1} \mA^\top)^\dagger}^2.
\end{equation}
Now, since $\mA$ is invertible, we have $(\mA \mW_k^{-1} \mA^\top)^\dagger = {\mA^\top}^{-1} \mW_k {\mA}^{-1}$.
So, if $x^*$ is the unique solution to~\eqref{eqn:problin}, then we obtain
 \begin{equation}\label{eq:fhatlinear2}
   \hat{f}_k(x)  \;=\;  \frac{1}{2} \norm{x-x^*}_{\mW_k}^2.
\end{equation}
Using  Assumption~\ref{eq:metricdescentquad} together with the fact that $\hat f_k$ is quadratic to conclude that
\begin{eqnarray*}
\hat f_{k+1}(x^{k+1})
& \overset{\eqref{eq:metricdescentquad}}{\leq } &
\hat f_k(x^{k+1}) \\
& = &\hat f_k(x^{k}) + \langle \nabla \hat f_k(x^{k}), x^{k+1} - x^k \rangle + \frac{1}{2} \Vert x^{k+1} - x^k \Vert^2_{\nabla^2 \hat f_k(x^{k})} \\
& = &\hat f_k(x^{k}) + \langle \nabla \hat f_k(x^{k}), x^{k+1} - x^k \rangle + \frac{1}{2} \Vert x^{k+1} - x^k \Vert^2_{\mW_k}.
\end{eqnarray*}
\end{proof}

\begin{proposition}\label{lem:linearcase}
Let $\mA \in \mathcal{M}_p(\mathbb{R})$ be invertible, $b \in \mathbb{R}^p$, and $x^*$ be the solution to~\eqref{eqn:problin}.
Let $(x^k)_{k \in \mathbb{N}}$ be a sequence generated from SNRVM \eqref{algo:SNRVM}, with $F(x)= \mA x-b$, and $\gamma=1$.
We assume that, at every iteration $k \in \mathbb{N}$, the matrices $(\mS_k,\mW_k)$ are sampled from a finite proper distribution (see Assumption \ref{Ass:distribution D_x}) such that for all $x$, $\mathcal{D}_{x}$ is independant of $x$, that $\mathbb{E}_{\mathcal{D}_{x}}[\mS \mS^\top]$ is invertible and $\mW_{k}$ is constant and equal to some invertible matrix $\mW \in \mathcal{M}_p(\mathbb{R})$.

Let \[\rho \; \eqdef \; \lambda_{\min}\left( \mW^{-1/2} \mA^\top \mathbb{E}_{\mathcal{D}_{x^0}}[\mS(\mS^\top \mA \mW^{-1} \mA^\top \mS)^\dagger \mS^\top] \mA \mW^{-1/2} \right).\]
It follows that $\rho \in (0,1)$, and
\begin{eqnarray}\label{eq:linear case cv rate}
\E{\Vert x^k - x^* \Vert^2_{\mW}}
&\leq& (1- \rho) \Vert x^2 - x^* \Vert^2_\mW.
\end{eqnarray}
\end{proposition}

\begin{proof}
We are going to apply the result in Theorem \ref{them: AdapSNR}, so we start by checking its assumptions.
First, our assumptions on the sampling ensure that Assumption \ref{Ass:distribution D_x} is verified.
Second, as discussed earlier in this section, the fact that $\mA$ is invertible ensures that Assumption \ref{ass:existence} holds true.
Third, our assumption that $\mW_k \equiv \mW$ together with Lemma \ref{L:nasty ass for linear case} tells us that Assumption \ref{ass:AdapUpperBnd} holds with $L=1$, meaning that we take a stepsize $\gamma = 1$.
Let $\mH \eqdef \mathbb{E}_{\mathcal{D}_{x^0}}[\mS(\mS^\top \mA \mW^{-1} \mA^\top \mS)^\dagger \mS^\top]$.
Note that this matrix is independant of $k$, because we assumed the distribution $\mathcal{D}_x$ to be independant of $x$.
We also know that $\mH$  is invertible, thanks to Lemma \ref{lem: sdjffskapp}.
Therefore, $\rho = \lambda_{min}\left( \mW^{-1/2} \mA^\top \mH \mA \mW^{-1/2} \right) >0$.

To prove that $\rho \leq 1$, observe that
\begin{eqnarray*}
&  & \hspace{-2cm}
\mW^{-1/2} \mA^\top \mathbb{E}[\mS(\mS^\top \mA \mW^{-1} \mA^\top \mS)^\dagger \mS^\top] \mA \mW^{-1/2} \\
& = &
\mathbb{E}[\mW^{-1/2}   \mA^\top \mS(\mS^\top \mA \mW^{-1/2}\mW^{-1/2} \mA^\top \mS)^\dagger \mS^\top \mA \mW^{-1/2} ] \\
& = & 
\mathbb{E}[ (\mS^\top \mA \mW^{-1/2} )^\top ((\mS^\top \mA \mW^{-1/2})(\mS^\top \mA \mW^{-1/2} ))^\dagger (\mS^\top \mA \mW^{-1/2} ) ] \\
& = &
\mathbb{E}[ (\mS^\top \mA \mW^{-1/2})^\dagger (\mS^\top \mA \mW^{-1/2} ) ],
\end{eqnarray*}
where $(\mS^\top \mA \mW^{-1/2})^\dagger (\mS^\top \mA \mW^{-1/2} )$ is the orthogonal projection onto the range of $\mW^{-1/2} \mA^\top \mS$.
Consequently $\lambda_{max}((\mS^\top \mA \mW^{-1/2})^\dagger (\mS^\top \mA \mW^{-1/2} )) \leq 1$, and from Jensen's inequality, we deduce that the eigenvalues of its expectation also are in $(0,1]$.
Whence $\rho \in (0,1]$.

To conclude the proof, see that under our assumptions, the quantity $\rho(x)$ defined in Theorem \ref{them: AdapSNR} is independent of $x$, and equal to $\rho$.
We have verified all the assumptions needed to call Theorem \ref{them: AdapSNR}, which proves the claim.
\end{proof}

The rate of convergence given in~\eqref{eq:linear case cv rate} is exactly the well known linear rate of convergence given in Theorem 4.6 in~\citep{Gower2015c}. 
For example, if $\mA$ is symmetric positive definite, and we can set $\mW_k \equiv \mA$ and sample the sketching matrix $\mS \in \mathbb{R}^{p \times 1}$ according to
\[ \Prob{\mS = e_i} \; = \; \frac{\mA_{ii}}{\trace{\mA}}, \quad \mbox{for }i=1,\ldots, m.\footnote{Here $e_i \in \R^{m}$ is the $i$-th unit coordinate vector and $\trace{\mA} = \sum_{i=1}^m\mA_{ii}$ is the trace of $\mA.$}\]
With this choice  of sketch and metric, the resulting method~\eqref{algo:SNRVM} is known as coordinate descent~\citep{Leventhal2010,Gower2015}. 
In this case, our resulting rate in~\eqref{eq:linear case cv rate} is  controlled by
\begin{equation*}
\rho = \lambda_{\min} \left(\mA^{1/2}\EE{}{ e_i\left(\mA_{ii} \right)^\dagger e_i^\top} \mA^{1/2} \right) = \frac{\lambda_{\min} \left(\mA \right)}{\trace{\mA}},
\end{equation*}

which is exactly the celebrated linear convergence rate of coordinate descent first given in~\citep{Leventhal2010}.

Proposition~\ref{lem:linearcase} shows that our main convergence theory in Theorem~\ref{them: AdapSNR} is \emph{tight} in this quadratic setting. That is, when specialized to a linear mapping $F(x)$ and a fixed metric $\mW \equiv \mW_k$, our Theorem~\ref{them: AdapSNR} recovers the best known convergences results as a special case.

\subsection{Proof of Theorem \ref{T:SNRVM explicit rates}}


\begin{lemma}\label{L:rho explicit lower bound}
Let Assumption \ref{Ass:explicit rates for SNRVM} hold.
Let $x \in \Omega$, let $(\hat\mS,\hat\mW)$ be in the domain of $\mathcal{D}_x$, and consider $\mA := \hat\mW^{-1/2}\nabla F(x) \mH(x) \nabla F(x)^\top \hat\mW^{-1/2}$, where $\mH(x)$ is defined in \eqref{eq:H matrix}.
Then 
\begin{equation*}
\lambda_{min}(\mA) \geq 
\frac{\mu_{\nabla F}^2}{L_{\nabla F}^{2}} 
\frac{\mu_{W}}{L_{W}}
\frac{\bar \mu_S}{L_S}
.
\end{equation*}
\end{lemma}

\begin{proof}
Let us write $\mH := \mH(x)$,  $J :=  \nabla F(x)$ and $U:= J^\top \hat\mW^{-1/2}$, so that $\mA = U^\top \mH U$.
Therefore,
\begin{equation*}
\lambda_{min}(\mA) \geq \lambda_{min}(U^\top U) \lambda_{min}(\mH).
\end{equation*}
From Assumption \ref{ass:existence} we know that $\hat\mW$ is invertible, and also that $J$ is injective, so from Assumption \ref{Ass:explicit rates for SNRVM} and the fact that $J$ is a square matrix, we deduce that $J$ is invertible, and therefore deduce that $U$ is invertible as well.
This means that $\lambda_{min}(U^\top U) >0$ and that 
\begin{equation*}
\lambda_{min}(U^\top U) =
\lambda_{min}(UU^\top)
=
\lambda_{min}(J^\top \hat\mW^{-1} J)
\geq
\lambda_{min}(J^\top J)\lambda_{min}(\hat\hat\mW^{-1})
=
\frac{\sigma_{min}(J)^2}{\lambda_{max}(\mW)}
\geq
\frac{\mu_{\nabla F}^2}{L_W}.
\end{equation*}
Now we turn to $\mH$, and write $\mH = \mathbb{E}\left[ \mS B \mS^\top \right]$, where $B=(\mS^\top G \mS)^\dagger$, with $G = J^\top \mW ^{-1} J$.
From the same arguments as above, we know that that $G$ is invertible under our assumptions.
So, using properties of the pseudo inverse with the fact that $G$ is symmetric and Lemma \ref{lem: lemma10_RSN}, we can write that
\begin{equation*}
\Null{B} = \Null{(\mS^\top G \mS)^\top} = \Null{\mS^\top G \mS}
=
\Null{\mS}.
\end{equation*}
Therefore, for all $x \in \mathbb{R}^p$ we have $\mS^\top x \in \Null{B}^\bot$.
So, by noting $\lambda_{min}^*(B)$ the smallest nonzero eigenvalue of $B$, we can write that
\begin{equation*}
\langle \mS B\mS^\top x,x \rangle
=
\langle B\mS^\top x, B\mS^\top \rangle 
\geq
\lambda_{min}^*(B) \Vert \mS^\top x \Vert^2
=
\lambda_{min}^*(B) \langle \mS\mS^\top x,x \rangle.
\end{equation*}
Here 
\begin{equation*}
\lambda_{min}^*(B) 
= \Vert B^\dagger \Vert^{-1} 
= \Vert \mS^\top G \mS \Vert^{-1}
\geq \Vert \mS\mS^\top \Vert^{-1} \Vert G \Vert^{-1}
\geq\Vert \mS\mS^\top \Vert^{-1} \Vert J^\top J \Vert^{-1} \Vert \mW ^{-1} \Vert^{-1}
\geq
L_S^{-1} L_{\nabla F}^{-2} \mu_W,
\end{equation*}
where we used the fact that $\Vert \mW ^{-1} \Vert^{-1} = \lambda_{min}(\mW )$
and $\Vert J^\top J \Vert^{-1} = \sigma_{max}(J)^{-2}$.
By combining those last inequalities we obtain that
\begin{eqnarray*}
\langle \mH x,x \rangle
&=&
\mathbb{E}\left[ \langle \mS B\mS^\top x,x \rangle \right] \\
&\geq &
L_S^{-1} L_{\nabla F}^{-2} \mu_W
\mathbb{E}\left[\langle {\mS\mS^\top} x,x \rangle \right] \\
&= &
L_S^{-1} L_{\nabla F}^{-2} \mu_W
\langle \mathbb{E}\left[{\mS\mS^\top}\right] x,x \rangle  \\
& \geq & L_S^{-1} L_{\nabla F}^{-2} \mu_W \bar \mu_S
\Vert x \Vert^2.
\end{eqnarray*}
This means that $\lambda_{min}(\mH) \geq L_S^{-1} L_{\nabla F}^{-2} \mu_W \bar \mu_S$.
If we recombine all our inequalities, we ultimately obtain that
\begin{equation*}
\lambda_{min}(\mA)
\geq
\frac{\mu_{\nabla F}^2}{L_W}
L_S^{-1} L_{\nabla F}^{-2} \mu_W \bar \mu_S,
\end{equation*}
which is what we needed.
\end{proof}

\begin{lemma}\label{L:fhat lower bound}
Let Assumption \ref{Ass:explicit rates for SNRVM} hold.
Let $x \in \Omega$, let $(\hat\mS,\hat\mW)$ be in the domain of $\mathcal{D}_x$.
Then
\begin{equation*}
\lambda_{min}\left( (\nabla F(x)^\top \mW^{-1} \nabla F(x))^\dagger \right)
\geq
L_{\nabla F}^{-2} \mu_W > 0.
\end{equation*}
In particular, for all $k \in \mathbb{N}$, if $x^k \in \Omega$ almost surely, then
	\begin{equation*}
	\mathbb{E}\left[  \hat f_k(x^k) \right]
	\geq 
	\frac{\mu_W}{2L_{\nabla F}^{2}} \mathbb{E}\left[ \Vert F(x^k) \Vert^2 \right]
	\quad 
	\text{ almost surely}.
	\end{equation*}
\end{lemma}

\begin{proof}
We have $\lambda_{min}\left( (\nabla F(x)^\top \mW^{-1} \nabla F(x))^\dagger \right) 
=
\Vert \nabla F(x)^\top \mW^{-1} \nabla F(x) \Vert^{-1}$
where
\begin{equation*}
\Vert \nabla F(x)^\top \mW^{-1} \nabla F(x) \Vert
\leq
\sigma_{max}(\nabla F(x))^2 \lambda_{min}(\mW)^{-1} 
\leq
L_{\nabla F}^2 \mu_W^{-1},
\end{equation*}
which gives the desired lower bound on the eigenvalues.
Now, given $x^k \in \Omega$ we immediately deduce that 
\begin{equation}\label{flb1}
 \hat f_k(x^k)
\geq 
\frac{\mu_W}{2L_{\nabla F}^{2}} \Vert F(x^k) \Vert^2.
\end{equation}
The conclusion follows by taking the expectation on this inequality.
\end{proof}

	\begin{proof}[Proof of Theorem \ref{T:SNRVM explicit rates}]
Keeping the notations of Theorem \ref{them: AdapSNR}, we see from Lemma \ref{L:rho explicit lower bound} that we can take
\begin{equation*}
\rho = \frac{\mu_{\nabla F}^2}{L_{\nabla F}^{2}} 
\frac{\mu_{W}}{L_{W}}
\frac{\bar \mu_S}{L_S} >0,
\end{equation*}
from which we obtain that
   \begin{eqnarray*}
(\forall k \in \mathbb{N}) \quad
 \E{\hat{f}_{k}(x^{k})} &\leq&   \left(1  -  \rho\gamma \right)^k \E{\hat{f}_0(x^0)}
 \quad 
 \text{ almost surely.}
\end{eqnarray*}
We now can use Lemma \ref{L:fhat lower bound} to lower bound the left member of that inequality, and obtain
\begin{equation}\label{snmvrer0}
(\forall k \in \mathbb{N}) \quad
 \frac{\mu_W}{2L_{\nabla F}^{2}} \mathbb{E}\left[\Vert F(x^k) \Vert^2\right] \leq  \left(1  -  \rho\gamma \right)^k \E{\hat{f}_0(x^0)}
 \quad 
 \text{ almost surely.}
\end{equation}
	The conclusion follows by taking 
	\begin{equation*}
	C=
	\E{\hat{f}_0(x^0)}
	\frac{L_{\nabla F}^{2}}{\mu_W}.
	\end{equation*}
\end{proof}

\subsection{Proof of convergence for SAN and SAN for bounded sequences}\label{SS:SAN bounded sequence}

\begin{proposition}\label{P:SAN constants on compact}
Let Assumption \ref{Ass:strict convexity of problem} hold.
For SAN and SANA, Assumption \ref{Ass:explicit rates for SNRVM} holds on every compact set $\Omega$.
\end{proposition}

\begin{proof}
First, remember that Assumption \ref{ass:existence} holds for SAN and SANA (see Proposition \ref{P:Ass verified for SAN and SANA}), and that $m=p=(n+1)d$.
Now, let $\Omega$ be a compact set, and verify that the bounds in Assumption \ref{Ass:explicit rates for SNRVM} hold.

We start with the sketching matrices $\mS$, for which we know (see the proof of Proposition \ref{P:Ass verified for SAN and SANA} in Section \ref{SS:proof of Ass verified for SAN and SANA}) that
\begin{equation*}
\Vert \mS \mS^\top \Vert =1 
\quad \text{ and } \quad 
\mathbb{E}\left[ \mS \mS^\top \right]
=
\Diag{\pi   \mI_d,\frac{1-\pi}{n} \mI_d,  \cdots, \frac{1-\pi}{n} \mI_d} \text{ or }\Diag{\mI_d, \frac{1}{n} \mI_d, \cdots, \frac{1}{n} \mI_d}.
\end{equation*}
In both cases, we see that we can take $L_S=1$ and $\bar \mu_L = \min\{\frac{1}{n}, \frac{1-\pi}{n}, \pi\}$.

Second, let $\mW_i$ be in the domain of $\mathcal{D}_x$.
According to their definition in (\ref{eq:rowsSW1-po},~\ref{eq:SANA distrib}), and because the $f_i$ is of class $C^2$ (see Assumption \ref{Ass:strict convexity of problem}),
we know that each $\mW_i$ is continuous with respect to $x$.
Moreover, we know (again from Assumption \ref{Ass:strict convexity of problem}) that $\mW_i$ is definite positive : 
$\lambda_{min}(\mW_i)>0$.
This is true for every $x \in \Omega$, so by continuity of $\lambda_{min}$ and the compactness of $\Omega$, we deduce that $\inf \limits_{x \in \Omega} \lambda_{min}(\mW_i)>0$.
Similarly, $\sup \limits_{x \in \Omega} \lambda_{max}(\mW_i)< + \infty$.
This means that the constants $\mu_W$ and $L_W$ are well defined in $(0, + \infty)$.

Finally, we need to control the singular values of $\nabla F(x)$ over $\Omega$.
We use here the same arguments that we used for $\mW_i$ : $\nabla F(x)$ is continuous with respect to $x$, and it is invertible (because it is square and injective, see Proposition \ref{P:Ass verified for SAN and SANA}).
\end{proof}

\begin{theorem}\label{T:CV SAN bounded sequence}
Let Assumptions \ref{Ass:strict convexity of problem} and \ref{ass:AdapUpperBnd} hold.
Let $\{x^k\}_{k \in \mathbb{N}}$ be a sequence generated by SAN with $\pi = 1/(n+1)$, or by SANA, with $\gamma=1/L$.
Suppose that   $\{x^k\}_{k \in \mathbb{N}}$ is bounded almost surely.
Then there exists $C, \rho >0$ such that, for every $k \in \mathbb{N}$, 
\begin{equation*}
\mathbb{E}\left[ \Vert F(x^k) \Vert^2 \right] \leq C(1-\gamma \rho)^k \text{ a.s.}
\end{equation*}
\end{theorem}

\begin{proof}
There exists a compact set $\Omega$ containing almost surely the sequence.
So it remains to combine Theorem \ref{T:SNRVM explicit rates} together with Proposition \ref{P:SAN constants on compact} and Proposition \ref{P:Ass verified for SAN and SANA}.
\end{proof}

\subsection{Proof of Theorem \ref{T:CV SAN explicit}}

The proof of Theorem \ref{T:CV SAN explicit}, which can be found at the end of this section, will combine  Theorem \ref{T:SNRVM explicit rates} with the forthcoming Propositions \ref{P:explicit constants from SAN to SNRVM} and \ref{P:SAN properties of F bijection}.

%
%
%

\begin{lemma}\label{L:varphi properties explicit rates}
Let $\varphi : [0,+\infty) \longrightarrow [1,+\infty)$ be defined as
\begin{equation}\label{D:varphi rates explicit}
\varphi(t):=\sqrt{1 + \frac{1}{2} \left( t + \sqrt{4t + t^2} \right)}.
\end{equation}
\begin{enumerate}
	\item\label{L:varphi properties explicit rates:monotone} $\varphi(t)$ is well defined and increasing on $[0,+\infty[$.
	\item\label{L:varphi properties explicit rates:inverse} $\varphi(t)^{-1}=\sqrt{1 + \frac{1}{2} \left( t - \sqrt{4t + t^2} \right)}$.
	\item\label{L:varphi properties explicit rates:asymptotic} For all $a \in (0,+\infty)$, $\varphi(at)\varphi(t^{-1})t^{-1/2}$ is decreasing on $(0,+\infty)$. 
	\item\label{L:varphi properties explicit rates:explicit bound} For all $a \in (0,+\infty)$, and all $t \geq 1$, $\varphi(at)\varphi(t^{-1}) \leq \varphi(1) \sqrt{t} \sqrt{2+a}$.
\end{enumerate}
\end{lemma}

\begin{proof}
\ref{L:varphi properties explicit rates:monotone} : It is well defined because $t + \sqrt{4t + t^2} \geq 0$.  
It is increasing because it is the composition, sum and product of increasing functions on $[0,+\infty[$. 
Point \ref{L:varphi properties explicit rates:inverse} is a simple exercise.
Point \ref{L:varphi properties explicit rates:asymptotic} is a bit more technical.
Let $\phi(t) = \varphi(at)^2\varphi(t^{-1})^2 t^{-1}$, which is the square of the quantity of interest.
We can compute its derivative, and a some effort we obtain that
\begin{equation*}
t^3 \phi'(t) = 
-\left[ t+\frac{1}{2} \left( 1+ \sqrt{1+4t} \right) \right]\left[ 1+ \frac{at}{\sqrt{4at+a^2t^2}} \right]
-
\frac{1}{2}
\left[ 1+\frac{1+2t}{\sqrt{1+4t}} \right]
\left[ 1+\frac{1}{2}\left( at+\sqrt{4at+a^2t^2} \right) \right].
\end{equation*}
It is clear that the above expression is the sum of two negative terms, implying that $\phi$ is decreasing.
For item \ref{L:varphi properties explicit rates:explicit bound},
we use the monotonicity of item \ref{L:varphi properties explicit rates:asymptotic} to get
\begin{equation*}
\varphi(at)\varphi(t^{-1}) \leq
\varphi(a\cdot 1)\varphi(1) \sqrt{t}.
\end{equation*}
From the inequality $a + \sqrt{a^2+4a} \leq 2a+2$  (it is easy to prove it by rearranging the terms and taking the square), we deduce that
\begin{equation*}
\varphi(a) \leq \sqrt{1+\frac{1}{2}\left( 2a+2 \right)} = \sqrt{2+a}.
\end{equation*}
\end{proof}

\begin{lemma}\label{L:singular values triangular block matrix}
Let $\mA \in \mathbb{R}^{m\times d}$ be an injective matrix.
Let $\varphi$ be defined as in \eqref{D:varphi rates explicit} and consider:
\begin{equation*}
A := 
\begin{bmatrix}
\mI_d & \mo_{d,m} \\
\mA & \mI_{m}
\end{bmatrix}.
\end{equation*}
Then $\Vert A \Vert = \varphi\left(\Vert \mA^\top \mA \Vert\right).$
\end{lemma}

\begin{proof}
We start by remembering that $\Vert A \vert$ is the largest singular value of $A$.
The singular values of $A$ are exactly the square root of the eigenvalues of $A^\top A$, which is given by
\begin{equation*}
A^\top A := 
\begin{bmatrix}
\mI_d + \mA^\top \mA & \mA^\top \\
\mA & \mI_{m}
\end{bmatrix}.
\end{equation*}
We compute its eigenvalues by finding the roots of its characteristic polynomial, that we note $P \in \mathbb{R}[X]$.
Using a simple formula for computing the determinant of a $2\times2$ block matrix, we can write, for all $X \neq -1$ : 
\begin{eqnarray*}
P(X) 
&= &
\det \begin{bmatrix}
(1-X)\mI_d + \mA^\top \mA & \mA^\top \\
\mA & (1-X)\mI_{m}
\end{bmatrix} \\ 
&=&
\det\left( (1-X)\mI_{m} \right)
\det\left( (1-X)\mI_d + \mA^\top \mA - \mA^\top((1-X)\mI_{m})^{-1}\mA \right) \\
&=&
(1-X)^m
\det\left( (1-X)\mI_d + \mA^\top \mA - \frac{1}{1-X}\mA^\top \mA \right)\\
&=&
(1-X)^{m-d}
\det\left( (1-X)^2\mI_d + (1-X)\mA^\top \mA - \mA^\top \mA \right)\\
&=&
(1-X)^{m-d}
\det\left( (1-X)^2\mI_d  -X\mA^\top \mA \right).
\end{eqnarray*}
The right member of this equality is polynomial in $X$, since the determinant of a matrix is polynomial in its coefficients, and our assumption that $\mA$ is injective implies that $m-d \geq 0$.
In particular this right member is well defined and continuous at $X = -1$, which means that the equality holds true for every $X \in \mathbb{R}$.

We see that $1$ is a root of $P$, with multiplicity $m-d$.
The other roots are the zeroes of $\det\left( (1-X)^2\mI_d  -X\mA^\top \mA \right)$, for which we see that
\begin{eqnarray*}
\det\left( (1-X)^2\mI_d  -X\mA^\top \mA \right) = 0 
& \Leftrightarrow &
\text{$(1-X)^2$ is an eigenvalue of $X \mA^\top \mA$} \\
& \Leftrightarrow &
\text{$(1-X)^2 = X \lambda$ for $\lambda \in {\rm{spec}}(\mA^\top \mA)$} \\
& \Leftrightarrow &
\text{$X = 1 + \frac{1}{2} \left( \lambda \pm \sqrt{4\lambda + \lambda^2} \right)$ for $\lambda \in {\rm{spec}}(\mA^\top \mA)$},
\end{eqnarray*}
which gives us the remaining $2d$ roots (counted with multiplicity).
This proves that the singular values of $A$ are $1$ (with multiplicity $m-d$) and (see Lemma \ref{L:varphi properties explicit rates}.\ref{L:varphi properties explicit rates:inverse})
\begin{equation*}
\{ \varphi(\lambda), \varphi(\lambda)^{-1} \ | \ \lambda \in {\rm spec}\left( \mA^\top \mA \right)\}.
\end{equation*}
Since $\varphi$ is increasing (Lemma \ref{L:varphi properties explicit rates}.\ref{L:varphi properties explicit rates:monotone}), and $\varphi(\lambda) \geq 1$, we conclude that the largest singular value of $A$ is $\varphi\left(\Vert \mA^\top \mA \Vert\right)$.

%
\end{proof}

\begin{proposition}\label{P:explicit constants from SAN to SNRVM}
Let Assumption \ref{Ass:explicit rates for SAN} hold, and consider the SAN (resp. SANA) algorithm.
Let $c = \sqrt{\frac{3+\sqrt{5}}{2}}$.
Then Assumption \ref{Ass:explicit rates for SNRVM} is verified, with $\Omega = \mathbb{R}^p$ and:
\begin{eqnarray*}
&\mu_W = \min\{1,\mu_f\}, \
&L_W = \max\{1, L_f\},\\
&\bar \mu_S = \min\left\{\frac{1-\pi}{n}, \pi\right\} \left(\text{resp. } \bar \mu_S =\frac{1}{n}\right),\
&L_S = 1, \\
&\mu_{\nabla F} = \frac{\mu_W}{c \sqrt{n} \sqrt{2+L_f^2}},\
&L_{\nabla F} = L_W c \sqrt{n} \sqrt{2+L_f^2}.
\end{eqnarray*}
\end{proposition}

%

\begin{proof}
Let $x \in \mathbb{R}^p$ be fixed, $J:= \nabla F(x) \in \mathbb{R}^{p \times p}$, and $(\mS,\mW)$ in the domain of $\mathcal{D}_x$.
We need to find uniform spectral bounds on those three matrices.
For $\mS$, we have seen already in the proof of Proposition \ref{P:SAN constants on compact} that we can take $L_S = 1$, and $\bar \mu_L = \frac{1}{n}$ for SAN, or $\min\{\frac{1-\pi}{n}, \pi\}$ for SANA.
For $\mW$, we see directly from (\ref{eq:rowsSW1-po},~\ref{eq:SANA distrib}) that it is a block-diagonal matrix, whose eigenvalues are included in $[\mu_W,L_W]$, with $\mu_W = \min\{1,\mu_f\}$ and $L_W=\max\{1,L_f\}$.
The rest of this proof is dedicated to the study of $J$, which requires more work.

Remember that the expression for $J$ is given in \eqref{eq: rsn_jac}. 
We  write for convenience that
\begin{equation*}
J = 
\begin{bmatrix}
\mo_d & \mH \\ \frac{1}{n}\mE^\top & - \mI_{nd}
\end{bmatrix}
\quad \text{ with } \quad
\mE:=\left[ \mI_d, \cdots, \mI_d \right] \in \mathbb{R}^{d \times nd},
\quad 
\mH:=\left[ H_1, \cdots, H_n \right] \in \mathbb{R}^{d \times nd}, 
\quad
H_i := \nabla^2 f_i(w).
\end{equation*}
Let us now introduce a few more matrices.
Let $\bar H:= \nabla^2 f(w)$, which can equivalently be written as $\bar H = \frac{1}{n}\sum_i H_i = \frac{1}{n} \mH \mE^\top$.
Now consider
\begin{equation*}
U := 
\begin{bmatrix}
\mI_d & \mo_{d,nd} \\ -\mH^\top & \mI_{nd}
\end{bmatrix},
\quad
D :=
\begin{bmatrix}
\bar H & \mo_{d,nd} \\ \mo_{nd,d} & -\mI_{nd}
\end{bmatrix},
\quad
V := 
\begin{bmatrix}
\mI_d & \mo_{d,nd} \\ \frac{-1}{n} \mE^\top & \mI_{nd}
\end{bmatrix}.
\end{equation*}
Note that those three matrices are triangular, and invertible because $\bar H$ is invertible.
It is easy to see that $J=U^\top DV$.
Indeed,
\begin{equation*}
DV
=
\begin{bmatrix}
\bar H & \mo_{d,nd} \\
\frac{1}{n}\mE^\top & -\mI_{nd}
\end{bmatrix},
\quad
U^\top DV
=
\begin{bmatrix}
\bar H - \mH \frac{1}{n}\mE^\top & \mH \\
\frac{1}{n}\mE^\top & - \mI_{nd}
\end{bmatrix}
=
\begin{bmatrix}
\mo_d & \mH \\ \frac{1}{n}\mE^\top & - \mI_{nd}
\end{bmatrix},
\end{equation*}
where the last equality comes from the fact that $\mH\frac{1}{n}\mE^\top = \bar H$.
Therefore, it remains to upper bound the right member of
\begin{equation*}
\sigma_{max}(J) = \sqrt{\lambda_{max}(J^\top J)} = \Vert J \Vert \leq \Vert D \Vert \Vert U \Vert \Vert V \Vert.
\end{equation*}
Bounding the smallest singular value will follow the same argument.
Indeed, from our assumptions, $J^T J$ is invertible (see Proposition \ref{P:Ass verified for SAN and SANA}), but $J$ is a square matrix, therefore $J$ itself is invertible.
In consequence, we can write that $J^{-1} = {V^{-1}}^\top D^{-1} U^{-1}$, so that
\begin{equation*}
\sigma_{min}(J)=\frac{1}{\Vert J^{-1}\Vert}
\geq
\frac{1}{\Vert D^{-1} \Vert \Vert U^{-1} \Vert \Vert V^{-1} \Vert},
\end{equation*}
where one easily computes that
\begin{equation*}
U^{-1} = 
\begin{bmatrix}
\mI_d & \mo_{d,nd} \\ \mH^\top & \mI_{nd}
\end{bmatrix},
\quad
D^{-1} =
\begin{bmatrix}
\bar H^{-1} & \mo_{d,nd} \\ \mo_{nd,d} & -\mI_{nd}
\end{bmatrix},
\quad
V^{-1} = 
\begin{bmatrix}
\mI_d & \mo_{d,nd} \\ \frac{1}{n} \mE^\top & \mI_{nd}
\end{bmatrix}.
\end{equation*}
It is easy to see, given our smoothness and strong convexity assumptions, that 
\begin{equation*}
\Vert D \Vert = \max\{1,L\} = L_W
\quad \text{ and } \quad 
\Vert D^{-1} \Vert = \max\{1, \mu^{-1}\} = \mu_W^{-1}.
\end{equation*}
Now, observe that $V$ and $V^{-1}$ share the same structure, so we can call Lemma \ref{L:singular values triangular block matrix} with $\mA=\frac{1}{n}\mE^T$ or $-\frac{1}{n}\mE^T$.
In both cases $\mA^\top \mA = n^{-1}\mI_{d}$, meaning that  $\Vert\mA^\top \mA \Vert ={n}^{-1}$, and so we deduce that 
\begin{equation*}
\Vert V \Vert = \Vert V^{-1} \Vert = \varphi(n^{-1}).
\end{equation*}
Finally, we do the same for $U$ and $U^{-1}$ : we use  Lemma \ref{L:singular values triangular block matrix} with $\mA = \pm \mH^\top$.
In both cases $\mA^\top \mA = \mH \mH^\top = \sum_{i=1}^n H_i^2$, whose eigenvalues belong to $[n\mu^2, nL^2]$.
Due to the monotonicity of $\varphi$, we deduce that
\begin{equation*}
\Vert U \Vert 
=
\Vert U^{-1} \Vert \leq \varphi(nL^2).
\end{equation*}
As a result, we conclude that
\begin{equation*}
\sigma_{min}(J)
\geq
\frac{1}{\mu_W^{-1}\varphi({n}^{-1})\varphi(nL^2)}
=
\frac{\mu_W}{\varphi({n}^{-1})\varphi(nL^2)},
\end{equation*}
while
\begin{equation*}
\sigma_{max}(J)
\leq
L_W\varphi({n}^{-1})\varphi(nL^2).
\end{equation*}
We can then conclude by using Lemma \ref{L:varphi properties explicit rates}.\ref{L:varphi properties explicit rates:explicit bound} and noting $c=\varphi(1)$.
\end{proof}

\begin{proposition}\label{P:SAN properties of F bijection}
Let Assumption \ref{Ass:explicit rates for SAN} hold, and consider the SAN (resp. SANA) algorithm.
Then 
\begin{enumerate}
	\item $F : \mathbb{R}^p \longrightarrow \mathbb{R}^p$ is a diffeomorphism
	\item $F^{-1} : \mathbb{R}^p \longrightarrow \mathbb{R}^p$ is Lipschitz continuous : 
	\begin{equation*}
	(\forall x,y \in \mathbb{R}^p)
	\quad
	\Vert x - y \Vert \leq \frac{c \sqrt{n} \sqrt{2+L_f^2}}{\min\{1,\mu_f\}} \Vert F(x)- F(y) \Vert
	\quad \text{ with } \quad
	c = \sqrt{\frac{3+\sqrt{5}}{2}}.
	\end{equation*}
	\item $F^{-1}(0) = \left[ w^*,\ \nabla f_1(w^*),\ \dots,\ \nabla f_n(w^*) \right]$, with $w^* = {\rm{argmin}}~f$.
\end{enumerate}
\end{proposition}

\begin{proof}
Let us start by showing that $F$ is injective.
For $x, \hat x \in \mathbb{R}^p$, we have
\begin{eqnarray}\label{spof1}
F(x) = F(\hat x) 
& \Rightarrow &
\frac{1}{n}\sum_i \alpha_i = \frac{1}{n}\sum_i \hat \alpha_i
\quad \text{ and } \quad 
\nabla f_i(w) - \alpha_i = \nabla f_i(\hat w) - \hat \alpha_i, 
\quad \forall i=1,\dots, n.
\end{eqnarray}
Summing the right member over $i$, we obtain that 
$\nabla f(w) - \nabla f(\hat w) = \frac{1}{n}\sum_i\alpha_i - \frac{1}{n}\sum_i \hat \alpha_i =0$.
In other words, we obtained that $\nabla f(w) = \nabla f(\hat w)$.
Now, we assumed that $f$ is strongly convex, therefore $\nabla f$ is injective (this can be seen from the fact that $\nabla f$ is strongly monotone).
So we deduce that $w = \hat w$.
Going back to \eqref{spof1}, we see now that $\alpha_i = \hat \alpha_i$, from which we conclude that $x= \hat x$, and that $F$ is indeed injective.

We know that $\nabla F(x)$ is invertible for all $x \in \mathbb{R}^p$ (see Proposition \ref{P:Ass verified for SAN and SANA}), so we can use the global inversion theorem to deduce that $F$ is a diffeomorphism between $\mathbb{R}^p$ and $F(\mathbb{R}^p)$.
Let us prove now that $F(\mathbb{R}^p) = \mathbb{R}^p$.
For this we will use an argument analog to what we used in the proof of Proposition \ref{P:explicit constants from SAN to SNRVM}.
Let $u,d,v^\top : \mathbb{R}^p \longrightarrow \mathbb{R}^p$ be defined as
\begin{eqnarray*}
u(x) := (w, -\nabla f_1(w)+\alpha_1, \dots, -\nabla f_n(w) + \alpha_n),
\quad
v^\top(x) := (w - \frac{1}{n}\sum\limits_{i=1}^n \alpha_i, \alpha_1, \dots, \alpha_n), \\
d(x) := (\nabla f(w), -\alpha_1, \dots, -\alpha_n).
\end{eqnarray*}
It is a simple exercise to verify that $F = v^\top \circ d \circ u $.
We will see now that those three functions are invertible and that their inverses are given by
\begin{eqnarray*}
u^{-1}(x) := (w, \nabla f_1(w)+\alpha_1, \dots, \nabla f_n(w) + \alpha_n),
\quad
{v^\top}^{-1}(x) := (w + \frac{1}{n}\sum\limits_{i=1}^n \alpha_i, \alpha_1, \dots, \alpha_n), \\
d^{-1}(x) := (\nabla f^*(w), -\alpha_1, \dots, -\alpha_n).
\end{eqnarray*}
For $u$, it is again a simple exercise to verify that $u \circ u^{-1} = u^{-1} \circ u = id_{\mathbb{R}^p}$.
Same for $v$ (which actually is linear).
For $d$, we used the notation $f^*$ which refers to the Fenchel transform of $f$,
\begin{equation*}
f^*(u) := \sup\limits_{w \in \mathbb{R}^d}  \langle u,w \rangle - f(w).
\end{equation*}
We assumed $f$ to be strongly convex and smooth, which means that $f^*$ is well-defined and differentiable on $\mathbb{R}^p$, and that  $(\nabla f)^{-1} = \nabla f^*$  \cite[Theorems 13.37, 16.29 \& 18.15]{BauCom}.
This proves that our expression for $d^{-1}$ is correct.
Now we conclude, by seeing that $u^{-1}, {v^\top}^{-1},d^{-1}$ are well-defined on $\mathbb{R}^p$, that $F^{-1}$ is well-defined on $\mathbb{R}^p$, and so that $F(\mathbb{R}^p)=\mathbb{R}^p$.

Now, we focus on $F^{-1} : \mathbb{R}^p \longrightarrow \mathbb{R}^p$.
It is differentiable everywhere, so  we can use the mean value theorem to deduce that $F^{-1}$ is Lipschitz continuous on $\mathbb{R}^p$, with a Lipschitz constant being bounded by:
	\begin{equation*}
		\sup\limits_{y \in \mathbb{R}^p} \Vert \nabla F^{-1}(y) \Vert 
		=
		\sup\limits_{x \in \mathbb{R}^p} \Vert \nabla F^{-1}(F(x)) \Vert 
		=
		\sup\limits_{x \in \mathbb{R}^p} \Vert \nabla F(x)^{-1} \Vert 		
		= 
		\sup\limits_{x \in \mathbb{R}^p} \frac{1}{\sigma_{min}(\nabla F(x))} 
		\leq \frac{1}{\mu_{\nabla F}},
	\end{equation*}
	where $\mu_{\nabla F}$ was computed in Proposition \ref{P:explicit constants from SAN to SNRVM}.
	We obtain that
	\begin{equation*}
		\Vert x - y \Vert
		=
		\Vert F^{-1}(F(x)) - F^{-1}(y) \Vert
		\leq 
		\frac{1}{\mu_{\nabla F}} \Vert F(x) - F(y) \Vert.
	\end{equation*}
	To conclude about the expression of $F^{-1}(0)$, compute 
	\begin{equation*}
	F^{-1}(x) = (u^{-1} \circ d^{-1} \circ {v^\top}^{-1})(x)
	=
	\left[ \hat w,\ 
	\nabla f_1 \left( \hat w \right) - \alpha_1,\
	\cdots, \
	\nabla f_n \left( \hat w \right) - \alpha_n	
	 \right], 
	 \text{ with } 
	 \hat w = \nabla f^* \left( w+\frac{1}{n}\sum\limits_{i=1}^n \alpha_i \right),
	\end{equation*}
	and use the fact that $\nabla f^*(0) = w^*$.
\end{proof}


\begin{proof}[Proof of Theorem \ref{T:CV SAN explicit}]
Let $w^* = {\rm{argmin}}~f$, and $x^* = \left[ w^*,\ \nabla f_1(w^*),\ \dots,\ \nabla f_n(w^*) \right]$, such that $F(x^*) = 0$ according to Proposition \ref{P:SAN properties of F bijection}.
Using again Proposition \ref{P:SAN properties of F bijection}, we obtain for all $k \in \mathbb{N}$ that
\begin{equation*}
	\Vert x^k - x^* \Vert^2 \leq \frac{c^2 {n} ({2+L_f^2})}{\min\{1,\mu_f^2\}} \Vert F(x^k) \Vert^2.
\end{equation*}
Here we have
\begin{equation*}
\Vert x^k - x^* \Vert^2 = 
\Vert w^k - w^* \Vert^2 
+
\sum\limits_{i=1}^n \Vert \alpha_i^k - \nabla f_i(w^*) \Vert^2.
\end{equation*}
Taking the expectation on the above expressions, and using the fact that $c^2 \leq 3$, we obtain that 
\begin{equation*}
\mathbb{E}\left[ \Vert w^k - w^* \Vert^2  \right]
+ \sum\limits_{i=1}^n \mathbb{E}\left[   \Vert \alpha_i^k - \nabla f_i(w^*) \Vert^2 \right]
\leq
\frac{3 {n} ({2+L_f^2})}{\min\{1,\mu_f^2\}} \mathbb{E}\left[ \Vert F(x^k) \Vert^2 \right].
\end{equation*}
Now we can use Theorem \ref{T:SNRVM explicit rates} together with Proposition \ref{P:explicit constants from SAN to SNRVM}, and combine the constants into play.
For the sake of the presentation, we assume $\pi = 1/(n+1)$, so that we can have a unique lower bound for SAN and SANA : $\bar \mu_S = \min\{\frac{1}{n},\frac{1}{n+1}\} \geq \frac{1}{2n}$.
We also simplify the expression of $c$, by again using bounds  like $c^2 \leq 3$ or $c^4 \leq 7$.
This allows us to write that
\begin{equation*}
	\mathbb{E}\left[ \Vert F(x^k) \Vert^2 \right] \leq C'(1-\gamma \rho)^k 
	\quad \text{ almost surely},
	\end{equation*}
	with
	$\rho = \frac{\min\{1,\mu_f^3\}}{14n^3(2+L_f^2)^2\max\{1,L_f^3\}}$, 
	and
$C'=
6n\E{\hat{f}_0(x^0)}
\frac{\max\{1,L_f^2\}(2+L_f^2)}{\min\{1,\mu_f\}}.
$	
	The conclusion follows by taking $C = C' \frac{3 {n} ({2+L_f^2})}{\min\{1,\mu_f^2\}}$:
\begin{equation*}
\mathbb{E}\left[ \Vert w^k - w^* \Vert^2  \right]
+ \sum\limits_{i=1}^n \mathbb{E}\left[   \Vert \alpha_i^k - \nabla f_i(w^*) \Vert^2 \right]
\leq
18n^2\E{\hat{f}_0(x^0)}
\frac{\max\{1,L_f^2\}(2+L_f^2)^2}{\min\{1,\mu_f^3\}} (1-\gamma \rho)^k.
\end{equation*}
\end{proof}

\end{document}